\documentclass[11pt]{article}

\usepackage[margin=1in]{geometry}
\usepackage{lipsum}
\usepackage{amsfonts}
\usepackage{graphicx}
\graphicspath{ {/} }
\usepackage{epstopdf}
\usepackage[normalem]{ulem}

\usepackage{mathtools}
\usepackage{empheq}
\usepackage{amsfonts}
\usepackage{amsgen,amsmath,amstext,amsbsy,amsopn,amssymb,subfigure,stmaryrd,mathabx}
\usepackage{comment}
\usepackage{enumitem}
\usepackage{natbib}
\usepackage{booktabs}
\usepackage{blkarray}
\usepackage{algorithm}
\usepackage{algorithmic}
\usepackage{url}
\usepackage{longtable}
\usepackage{multirow}
\usepackage{xcolor}
\usepackage{hyperref}

\allowdisplaybreaks

\newtheorem{Definition}{Definition}
\newtheorem{Theorem}{Theorem}

\newtheorem{Lemma}{Lemma}
\newtheorem{Remark}{Remark}

\newtheorem{Proposition}{Proposition}

\newcommand{\bea}{\begin{eqnarray}}
\newcommand{\eea}{\end{eqnarray}}
\newcommand{\beas}{\begin{eqnarray*}}
\newcommand{\eeas}{\end{eqnarray*}}

\newcommand{\Range}{{\rm Ran}}

\newcommand{\cA}{{\mathcal{A}}}

\newcommand{\cL}{{\mathcal{L}}}

\newcommand{\0}{{\mathbf{0}}}

\renewcommand{\a}{{\mathbf{a}}}

\newcommand{\e}{{\mathbf{e}}}
\renewcommand{\d}{{\mathbf{d}}}

\renewcommand{\u}{{\mathbf{u}}}
\renewcommand{\v}{{\mathbf{v}}}
\newcommand{\w}{{\mathbf{w}}}
\newcommand{\m}{{\mathbf{m}}}
\newcommand{\n}{{\mathbf{n}}}

\newcommand{\x}{{\mathbf{x}}}
\newcommand{\y}{{\mathbf{y}}}
\newcommand{\z}{{\mathbf{z}}}
\newcommand{\X}{{\mathbf{X}}}
\newcommand{\B}{{\mathbf{B}}}

\newcommand{\D}{{\mathbf{D}}}
\newcommand{\I}{{\mathbf{I}}}
\newcommand{\M}{{\mathbf{M}}}
\newcommand{\N}{{\mathbf{N}}}
\newcommand{\bF}{{\mathbf{F}}}

\newcommand{\G}{{\mathbf{G}}}
\renewcommand{\L}{{\mathbf{L}}}

\newcommand{\Q}{{\mathbf{Q}}}
\newcommand{\R}{{\mathbf{R}}}
\renewcommand{\S}{{\mathbf{S}}}
\newcommand{\U}{{\mathbf{U}}}
\newcommand{\V}{{\mathbf{V}}}
\renewcommand{\H}{{\mathbf{H}}}

\newcommand{\F}{{\rm F}}
\newcommand{\grad}{{\rm grad}}
\newcommand{\W}{{\mathbf{W}}}
\newcommand{\A}{{\mathbf{A}}}
\newcommand{\Y}{{\mathbf{Y}}}
\newcommand{\Z}{{\mathbf{Z}}}

\newcommand{\bepsilon}{{\boldsymbol{\epsilon}}}
\newcommand{\bSigma}{{\boldsymbol{\Sigma}}}

\newcommand{\bvarepsilon}{\boldsymbol{\varepsilon}}

\newcommand{\bDelta}{\boldsymbol{\Delta}}

\newcommand{\rank}{{\rm rank}}

\newcommand{\tr}{{\rm tr}}

\newcommand{\diag}{{\rm diag}}

\newcommand{\QR}{{\rm QR}}
\newcommand{\rmvec}{{\rm vec}}

\newcommand{\Hess}{{\rm Hess}}

\newcommand{\Newton}{{\rm Newton}}

\newcommand{\argmin}{\mathop{\rm arg\min}}
\newcommand{\argmax}{\mathop{\rm arg\max}}
\newcommand{\bbR}{\mathbb{R}}

\newcommand{\bbO}{\mathbb{O}}
\newcommand{\bbP}{\mathbb{P}}
\newcommand{\bbE}{\mathbb{E}}
\newcommand{\cM}{\mathcal{M}}

\newcommand{\inprod}[2]{\langle #1 , #2 \rangle}

\allowdisplaybreaks

\definecolor{lred}{rgb}{1,0.8,0.8}
\definecolor{lblue}{rgb}{0.8,0.8,1}
\definecolor{dred}{rgb}{0.8,0,0}
\definecolor{dblue}{rgb}{0,0,0.5}
\definecolor{dgreen}{rgb}{0,0.5,0}

\begin{document}


%
\title{Recursive Importance Sketching for Rank Constrained Least Squares: Algorithms and High-order Convergence}
\author{Yuetian Luo$^1$, ~ Wen Huang$^2$, ~ Xudong Li$^3$, ~ and ~ Anru R. Zhang$^{4}$}
\date{}
\maketitle
%
%
\footnotetext[1]{Department of Statistics, University of Wisconsin-Madison \texttt{yluo86@wisc.edu} . Y. Luo would like to thank RAship from Institute for Foundations of Data Science at UW-Madison.  }
\footnotetext[2]{School of Mathematical Sciences, Xiamen University \texttt{wen.huang@xmu.edu.cn} }
\footnotetext[3]{School of Data Science, Fudan University  \texttt{lixudong@fudan.edu.cn}}
\footnotetext[4]{Department of Biostatistics \& Bioinformatics, Computer Science, Mathematics, and Statistical Science, Duke University \texttt{anru.zhang@duke.edu}}


\maketitle

\begin{abstract}
In this paper, we propose {\it \underline{R}ecursive} {\it \underline{I}mportance} {\it \underline{S}ketching} algorithm for {\it \underline{R}ank} constrained least squares {\it \underline{O}ptimization} (RISRO). 
The key step of RISRO is recursive importance sketching, a new sketching framework based on deterministically designed recursive projections, which significantly differs from the randomized sketching in the literature \citep{mahoney2011randomized,woodruff2014sketching}. Several existing algorithms in the literature can be reinterpreted under this new sketching framework and RISRO offers clear advantages over them. RISRO is easy to implement and computationally efficient, where the core procedure in each iteration is to solve a dimension-reduced least squares problem. We establish the local quadratic-linear and quadratic rate of convergence for RISRO under some mild conditions. We also discover a deep connection of RISRO to the Riemannian Gauss-Newton algorithm on fixed rank matrices. The effectiveness of RISRO is demonstrated in two applications in machine learning and statistics: low-rank matrix trace regression and phase retrieval. Simulation studies demonstrate the superior numerical performance of RISRO.
\end{abstract}


%


\section{Introduction}\label{sec:intro}

The focus of this paper is on the rank constrained least squares:
\begin{equation}\label{eq:minimization}
\min_{\X \in \mathbb{R}^{p_1\times p_2}} f(\X) := \frac{1}{2} \left\|\y - \mathcal{A}(\X)\right\|_2^2, \quad \text{subject to}\quad \rank(\X)= r.
\end{equation}
Here, $\y\in\mathbb{R}^n$ is the given response and $\mathcal{A}\in \mathbb{R}^{p_1\times p_2} \to \mathbb{R}^n$ is a known linear map that can be explicitly represented as
\begin{equation} \label{eq: affine operator}
\mathcal{A}(\X) = \left[
\langle\A_1, \X\rangle, \ldots,
\langle\A_n, \X\rangle \right]^\top, \quad \langle \A_i, \X\rangle = \sum_{1\leq j\leq p_1, 1\leq k\leq p_2}(\A_{i})_{[j,k]} \X_{[j,k]}
\end{equation}
with given measurement matrices $\A_i \in \mathbb{R}^{p_1\times p_2}$, $i=1,\ldots,n$. 

The rank constrained least squares \eqref{eq:minimization} is motivated by the widely studied low-rank matrix recovery problem, where the goal is to recover a low-rank matrix $\X^*$ from the observation $\y = \cA(\X^*) + \bepsilon$ ($\bepsilon$ is the noise). This problem is of fundamental importance in a variety of fields such as optimization, machine learning, signal processing, scientific computation, and statistics. With different realizations of $\cA$, \eqref{eq:minimization} covers many applications, such as matrix trace regression \citep{candes2011tight,davenport2016overview}, matrix completion \citep{candes2010power,keshavan2009matrix,koltchinskii2011nuclear,miao2016rank}, phase retrieval \citep{candes2013phaselift,shechtman2015phase}, blind deconvolution \citep{ahmed2013blind}, and matrix recovery via rank-one projections \citep{cai2015rop,chen2015exact}. To overcome the non-convexity and NP-hardness of directly solving \eqref{eq:minimization} \citep{recht2010guaranteed}, various computational feasible schemes have been developed in the past decade, including the prominent convex relaxation \citep{recht2010guaranteed,candes2011tight}:
\begin{equation} \label{eq: convex relaxation}
	\min_{\X \in \mathbb{R}^{p_1\times p_2} } \frac{1}{2}\|\y - \cA (\X)\|_2^2 + \lambda \|\X\|_{*},
\end{equation}
where $\|\X\|_* = \sum_{i=1}^{\min(p_1, p_2)} \sigma_i(\X)$ is the nuclear norm of $\X$ and $\lambda > 0$ is a tuning parameter. Nevertheless, the convex relaxation technique has one well-documented limitation: the parameter space after relaxation is usually much larger than that of the target problem. Also, algorithms for solving the convex program often require the singular value decomposition as the stepping stone and can be prohibitively time-consuming for large-scale instances. 

In addition, non-convex optimization, which directly enforces the rank $r$ constraint on the iterates, renders another important class of algorithms for solving \eqref{eq:minimization}. Since each iterate lies in a low dimensional space, the computation cost of the non-convex approach can be much smaller than the convex regularized approach. Over the last a few years, there is a flurry of research on non-convex methods in solving \eqref{eq:minimization} \citep{chen2015fast,hardt2014understanding,jain2013low,miao2016rank,sun2015guaranteed,tran2016extended,tu2016low,wen2012solving,zhao2015nonconvex,zheng2015convergent}, and many of the algorithms such as gradient descent and alternating minimization are shown to have nice convergence results under proper model assumptions \citep{hardt2014understanding,jain2013low,sun2015guaranteed,tong2020accelerating,tu2016low,zhao2015nonconvex}. We refer readers to Section \ref{sec: related literature} for more review of recent works. 

In the existing literature, many algorithms for solving \eqref{eq:minimization} either require careful tuning of hyper-parameters or have a convergence rate no faster than linear. Thus, we raise the following question: 

\vskip.5cm
\begin{center}
\noindent \parbox{0.92\textwidth}{
{\it Can we develop an easy-to-compute and efficient (hopefully has comparable per-iteration computational complexity as the first-order methods) algorithm with provable high-order convergence guarantees (possibly converge to a stationary point due to the non-convexity) for solving \eqref{eq:minimization}?}
}
\end{center}
\vskip.5cm

In this paper, we give an affirmative answer to this question by making contributions as outlined next. 

\subsection{Our Contributions} \label{sec: constribution}
We introduce an easy-to-implement and computationally efficient algorithm, {\it \underline{R}ecursive} {\it \underline{I}mportance} {\it \underline{S}ketching} for {\it \underline{R}ank} constrained least squares {\it \underline{O}ptimization} (RISRO), for solving \eqref{eq:minimization} in this paper. The proposed algorithm is tuning free and has the same per-iteration computational complexity as Alternating Minimization \citep{jain2013low}, as well as comparable complexity to many popular first-order methods such as iterative hard thresholding \citep{jain2010guaranteed} and gradient descent \citep{tu2016low} when $r\ll p_1, p_2, n$. We then illustrate the key idea of RISRO under a general framework of recursive importance sketching. This framework also renders a platform to compare RISRO and several existing algorithms for rank constrained least squares. 

\begin{figure}[h]
	\centering
	\subfigure[Noiseless low-rank matrix trace regression. Here, $\y_i = \langle \A_i, \X^* \rangle$ for $1 \leq i \leq n$, $\X^* \in \bbR^{p \times p}$ with $p = 100$, $\sigma_1(\X^*) = \cdots = \sigma_3(\X^*) = 3, \sigma_k(\X^*) = 0$ for $4 \leq k \leq 100$ and $\A_i$ has independently identically distributed (i.i.d.) standard Gaussian entries]{\includegraphics[width=0.4\textwidth]{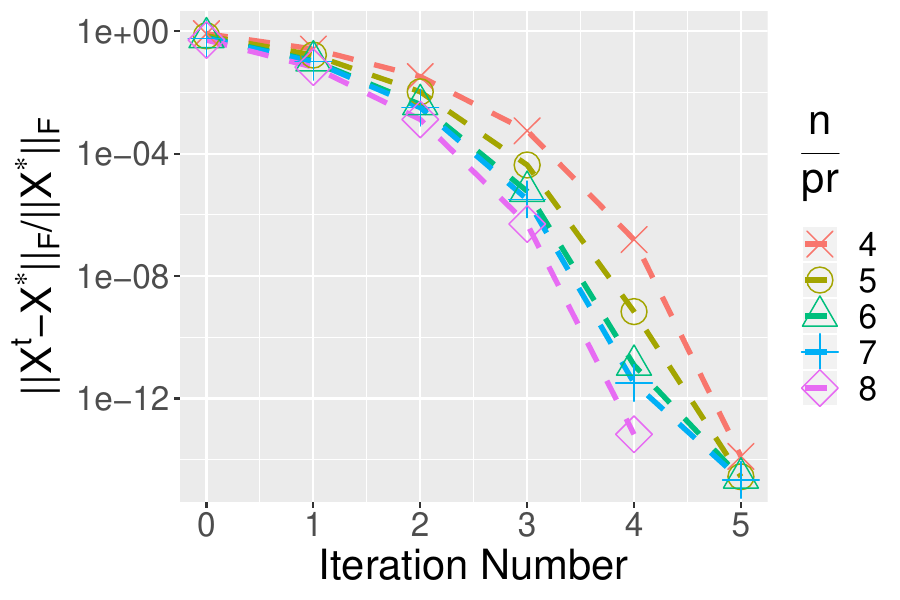}}\qquad 
	\subfigure[Phase Retrieval. Here, $\y_i = \langle \a_i \a_i^\top, \x^* \x^{*\top} \rangle$ for $1 \leq i \leq n$, $\x^* \in \bbR^{p}$ with $p = 1200$, $\a_i \overset{i.i.d.}\sim N(0,\I_{p})$ ]{\includegraphics[width=0.4\textwidth]{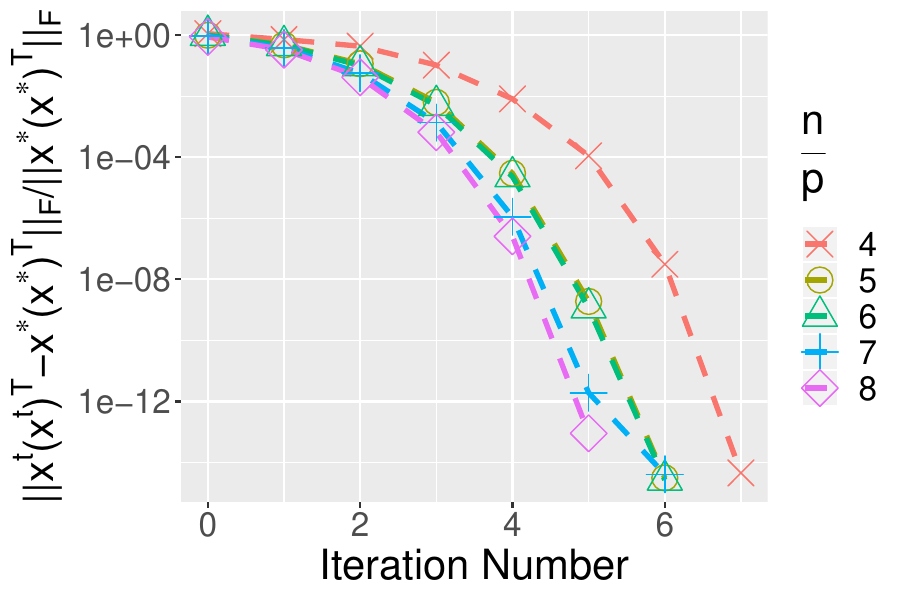}}
	\caption{RISRO achieves a quadratic rate of convergence (spectral initialization is used in each setting and more details about the simulation setup are given in Section \ref{sec:numerics})}\label{fig: Simul Mini performance illustration}
\end{figure}

 Assuming $\cA$ satisfies the restricted isometry property (RIP), we prove RISRO is local quadratic-linearly convergent in general and quadratically convergent {\it to a stationary point} under some extra conditions. Figure \ref{fig: Simul Mini performance illustration} provides a numerical example of the performance of RISRO in the noiseless low-rank matrix trace regression (left panel) and phase retrieval (right panel). In both problems, RISRO converges to the underlying parameter quadratically and reaches a highly accurate solution within five iterations. 

In addition, we discover a deep connection between RISRO and the Riemannian Gauss-Newton optimization algorithm on fixed rank matrices manifold. The least squares step in RISRO implicitly solves a \emph{Fisher Scoring} or \emph{Riemannian Gauss-Newton equation} on the Riemannian optimization of low-rank matrices and the updating rule in RISRO can be seen as a \emph{retraction map}. With this connection, our theory on RISRO also improves the existing convergence results on the Riemannian Gauss-Newton method for the rank constrained least squares problem.

Next, we further apply RISRO to two prominent problems in machine learning and statistics: low-rank matrix trace regression and phase retrieval. In the noisy low-rank matrix trace regression, we prove the {\it estimation error rate} of RISRO converges quadratically to the information-theoretical limit with only a {\it double-logarithmic} number of iterations under the Gaussian ensemble design. To the best of our knowledge, RISRO is the first algorithm that provably achieves the minimax rate-optimal estimation error in matrix trace regression with only a double-logarithmic number of iterations, which offers an exponential improvement over the existing results of first-order methods \citep{jain2010guaranteed,jain2013low,chen2015fast}. We also discover a new ``quadratic + one-iteration optimality" phenomenon for RISRO on low-rank matrix recovery (Remark \ref{rem: trace regression remark 1}). In phase retrieval, where $\cA$ does not satisfy the RIP condition, we can still establish the local convergence of RISRO given a proper initialization. We also develop RISRO in the matrix completion and robust PCA applications, where the restricted isometry property completely fails. We find RISRO still has similar empirical performance as in the setting where the RIP condition holds.

Finally, we conduct simulation studies to support our theoretical findings and compare RISRO with many existing algorithms. The numerical results show RISRO not only offers faster and more robust convergence but also requires a smaller sample size requirement for low-rank matrix recovery, compared to existing approaches. 

\subsection{Related Literature} \label{sec: related literature}

This work is related to a range of literature on low-rank matrix recovery, convex/non-convex optimization, and sketching arising from several communities, including optimization, machine learning, statistics, and applied mathematics. We make an attempt to review the related literature without claiming the survey is exhaustive.

One class of the most popular approaches to solve \eqref{eq:minimization} is the nuclear norm minimization (NNM) \eqref{eq: convex relaxation}. Many algorithms have been proposed to solve NNM, such as proximal gradient descent \citep{toh2010accelerated}, fixed-point continuation (FPC) \citep{goldfarb2011convergence}, and proximal point methods \citep{jiang2014partial}. It has been shown that the solution of NNM has desirable properties under proper models, such as matrix trace regression and matrix completion  \citep{cai2013sharp,cai2014sparse,cai2015rop,candes2011tight,recht2010guaranteed}. In addition to NNM, the max norm minimization is another widely considered convex realization for the rank constrained optimization \citep{lee2010practical,cai2013max}. However, these convex programs are usually computationally intensive to solve, which motivates a line of work on using non-convex approaches. Since \cite{burer2003nonlinear}, one of the most popular non-convex methods for solving \eqref{eq:minimization} is to first factor the low-rank matrix $\X$ to $\R \L^\top$ with two factor matrices $\R \in \bbR^{p_1 \times r}, \L \in \bbR^{p_2 \times r}$, then run either gradient descent or alternating minimization on $\R$ and $\L$ \citep{candes2015phase,li2019rapid,ma2019implicit,park2018finding,sanghavi2017local,sun2015guaranteed,tu2016low,wang2017unified,zhao2015nonconvex,zheng2015convergent,tong2020low}. Other methods, such as singular value projection or iterative hard thresholding \citep{goldfarb2011convergence,jain2010guaranteed,tanner2013normalized}, Grassmann manifold optimization \citep{boumal2011rtrmc,keshavan2009matrix}, Riemannian manifold optimization \citep{huang2018blind,meyer2011linear,mishra2014fixed,vandereycken2013low,wei2016guarantees} have also been proposed and studied. We refer readers to the recent survey paper \cite{chi2019nonconvex} for a comprehensive overview of existing literature on convex and non-convex approaches on solving \eqref{eq:minimization}. Most of the convergence analyses in the literature were conducted under certain statistical models (e.g., noisy/noiseless matrix trace regression, matrix completion, and phase retrieval) and the goal was to recover the underlying parameter matrix. Here, we study \eqref{eq:minimization} from both an optimization perspective (how the algorithm {\it converges to a stationary point}) and a statistical perspective (how the iterates {\it estimate the underlying true parameter}). These two perspectives overlap in the noiseless settings as the parameter becomes a stationary point then, while disjoint in the more general noisy settings.

There are a few recent attempts in connecting the geometric structures of different approaches \citep{ha2020equivalence,li2019non}, and the landscape of problem \eqref{eq:minimization} has also been studied in various settings \citep{bhojanapalli2016global,ge2017no,uschmajew2018critical,zhang2019sharp,zhu2018global}. 

Our work is also related to the idea of sketching in numerical linear algebra. Performing sketching to speed up the computation via dimension reduction has been explored extensively in recent years \citep{mahoney2011randomized,woodruff2014sketching}. Sketching methods have been applied to solve a number of problems including but not limited to matrix approximation \citep{song2017low,zheng2012practical,drineas2012fast}, linear regression  \citep{clarkson2017low,dobriban2019asymptotics,pilanci2016iterative,raskutti2016statistical}, ridge regression \citep{wang2017sketching}, etc. In most of the sketching literature, the sketching matrices are randomly constructed  \citep{mahoney2011randomized,woodruff2014sketching}. Randomized sketching 
matrices are easy to generate and require little storage for sparse sketching. However, randomized sketching can be suboptimal in statistical settings \citep{raskutti2016statistical}. To overcome this, \cite{zhang2020islet} introduced an idea of importance sketching in the context of low-rank tensor regression. In contrast to the randomized sketching, importance sketching matrices are constructed deterministically with the supervision of the data and are shown capable of achieving better statistical efficiency. However, the method developed is \cite{zhang2020islet} is essentially a ``one-time" importance sketching, which yield a sub-optimal outcome when the noise level is small or moderate. This paper proposes a more powerful recursive importance sketching algorithm that iteratively refines the sketching matrices. We also provide a comprehensive convergence analysis for the proposed algorithm without the sample-splitting assumption used in \cite{zhang2020islet}; our theory demonstrates the optimality of the proposed algorithm at all different noise levels and advantages over other algorithms for the rank constrained least squares problem.

\subsection{Organization of the Paper}\label{sec: organization}
The rest of this article is organized as follows. After a brief introduction of notation in Section \ref{sec: notation}, we present our main algorithm RISRO with an interpretation from the recursive importance sketching perspective in Section \ref{sec: main algorithm}. The theoretical results of RISRO are given in Section \ref{sec:theory}. In Section \ref{sec: Riemannian manifold interpre}, we present another interpretation for RISRO from Riemannian manifold optimization. The computational complexity of RISRO and its applications to low-rank matrix trace regression and phase retrieval are discussed in Sections \ref{sec: computation complexity} and \ref{sec: statistics applications}, respectively. Numerical studies of RISRO and the comparison with existing algorithms in the literature are presented in Section \ref{sec:numerics}. Conclusion and future work are given in Section~\ref{sec: conclusion}.

\subsection{Notation} \label{sec: notation}

The following notation will be used throughout this article. Upper and lowercase letters (e.g., $A, B, a, b$), lowercase boldface letters (e.g. $\u, \v$), uppercase boldface letters (e.g., $\U, \V$) are used to denote scalars, vectors, matrices, respectively. For any two series of numbers, say $\{a_n\}$ and $\{b_n\}$, denote $a = O(b)$ if there exists uniform constants $ C>0$ such that $a_n \leq Cb_n, \forall n$. For any $a, b \in \bbR$, let $a \wedge b := \min\{a,b\}, a \vee b = \max\{a,b\}$. For any matrix $\X \in \mathbb{R}^{p_1\times p_2}$ with singular value decomposition $\sum_{i=1}^{p_1 \land p_2} \sigma_i(\X)\u_i \v_i^\top$, where $\sigma_1(\X) \geq \sigma_2(\X) \geq \cdots \geq \sigma_{p_1 \wedge p_2} (\X)$, let $\X_{\max(r)}= \sum_{i=1}^{r} \sigma_i(\X)\u_i \v_i^\top$ be the best rank-$r$ approximation of $\X$ and denote $\|\X\|_\F = \sqrt{\sum_{i} \sigma^2_i(\X)}$ and $\|\X\| = \sigma_1(\X)$ as the Frobenius norm and spectral norm, respectively. Let $\QR(\X)$ be the $Q$ part of the QR decomposition outcome of $\X$. $\rmvec (\X) \in \bbR^{p_1 p_2}$ represents the vectorization of $\X$ by its columns. In addition, $\I_r$ is the $r$-by-$r$ identity matrix. Let $\mathbb{O}_{p, r} = \{\U: \U^\top \U=\I_r\}$ be the set of all $p$-by-$r$ matrices with orthonormal columns. For any $\U\in \mathbb{O}_{p, r}$, $P_{\U} = \U\U^\top$ represents the orthogonal projector onto the column space of $\U$; we also note $\U_\perp\in \mathbb{O}_{p, p-r}$ as the orthonormal complement of $\U$. We use bracket subscripts to denote sub-matrices. For example, $\X_{[i_1,i_2]}$ is the entry of $\X$ on the $i_1$-th row and $i_2$-th column; $\X_{[(r+1):p_1, :]}$ contains the $(r+1)$-th to the $p_1$-th rows of $\X$. For any matrix $\X$, we use $\X^\dagger$ to denote its Moore-Penrose inverse.  For matrices $\U\in \mathbb{R}^{p_1\times p_2}, \V\in \mathbb{R}^{m_1\times m_2}$, let
$$\U\otimes \V = \begin{bmatrix}
 \U_{[1,1]}\cdot \V & \cdots & \U_{[1, p_2]}\cdot \V\\
 \vdots & & \vdots\\
 \U_{[p_1,1]}\cdot \V & \cdots & \U_{[p_1, p_2]}\cdot \V\\
 \end{bmatrix}\in \mathbb{R}^{(p_1m_1)\times(p_2m_2)}$$ 
be their Kronecker product. Finally, for any given linear operator ${\cal L}$, we use $\cL^*$ to denote its adjoint, and use $\Range(\cL)$ to denote its range space.

\section{Recursive Importance Sketching for Rank Constrained Least Squares} \label{sec: main algorithm}

In this section, we discuss the procedure and interpretations of RISRO, then compare it with existing algorithms from a sketching perspective. The pseudocode of RISRO is summarized in Algorithm \ref{alg: recursive IS alg 2}. 

\subsection{RISRO Procedure and Recursive Importance Sketching} \label{sec: sketching for Simul Mini}

\begin{figure}
	\centering
	\subfigure{\includegraphics[width = 0.80\textwidth]{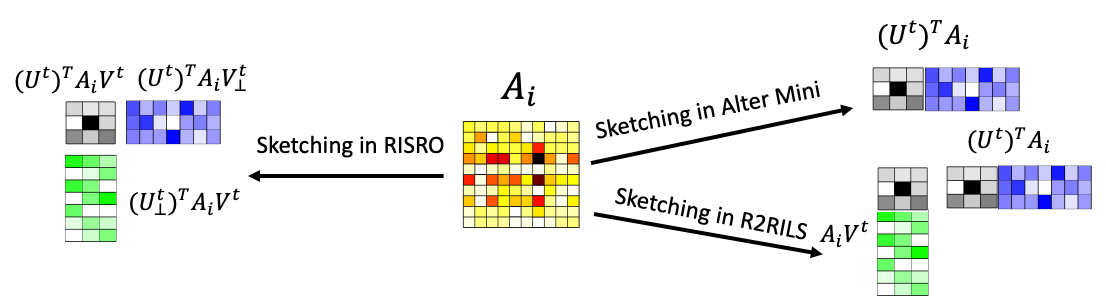}}
	\caption{ Illustration of sketching strategies of RISRO (this work), Alter Mini \citep{hardt2014understanding,jain2013low}, and R2RILS \citep{bauch2020rank}. Here, $\A_i$ denotes the covariate matrix of the $i$th observation; $\U^t$ and $\V^t$ span the column and row subspaces of $\X^t$, respectively. Covariate matrices colored in gray represent the sketching of $\A_i$ onto the column and row subspaces of $\X^t$, covariate matrices colored in green represent the sketching of $\A_i$ onto the perpendicular column subspace and row subspace of $\X^t$ and covariate matrices colored in blue represent the sketching of $\A_i$ onto the column subspace and perpendicular row subspace of $\X^t$. In Alter Mini and R2RILS, the sketched covariates colored in gray and blue (or green) are combined to represent the actual algorithmic implementation. }
	\label{fig: sketching illutstration}
\end{figure}

In each iteration $t=1,2,\ldots$, RISRO includes three steps. 
\begin{itemize}
\item[Step 1] We sketch each $\A_i$ ($i = 1,\ldots, n$) onto the subspace spanned by $[\U^t\otimes\V^t, \U_{\perp}^t\otimes\V^t, \U^t\otimes\V_{\perp}^t]$, where $\U^t$ and $\V^t$ span the column and row subspaces of $\X^t$, respectively. This yields the sketched importance covariates $\U^{t\top} \A_i \V^t, \U^{t\top}_\perp \A_i \V^t, \U^{t\top} \A_i \V^t_\perp$. See Figure \ref{fig: sketching illutstration} left panel for an illustration of the sketching scheme of RISRO. Then we construct the covariates maps $\cA_B: \bbR^{r \times r} \to \bbR^n$, $\cA_{D_1}: \bbR^{(p_1 - r)\times r} \to \bbR^n$ and $\cA_{D_2}: \bbR^{r \times (p_2 -r)} \to \bbR^n$: for matrix ``$\cdot$", let
		\begin{equation} \label{eq: importance sketches}
		[\cA_B(\cdot)]_i = \langle \cdot, \U^{t\top} \A_i \V^{t}\rangle,\quad 
		[\cA_{D_1}(\cdot)]_i = \langle \cdot, \U_{\perp}^{t\top} \A_{i} \V^{t}\rangle,\quad
		[\cA_{D_2}(\cdot)]_i = \langle \cdot, \U^{t\top} \A_{i} \V_{\perp}^{t}\rangle,		\quad i=1,\ldots, n.
		\end{equation} 
\item[Step 2] We solve a dimension reduced least squares problem \eqref{eq: alg2 least square} (provided in the box of Algorithm \ref{alg: recursive IS alg 2}) where the number of parameters is reduced to $(p_1 + p_2 -r)r$ while the sample size remains $n$. 
\item[Step 3] We update the sketching matrices $\U^{t+1}, \V^{t+1}$ and $\X^{t+1}$ in Steps~ \ref{alg:RISRO:st6} and \ref{alg:RISRO:st7}. By construction, $\U^{t+1}, \V^{t+1}$ contain both the column and row spans of $\X^{t+1}$. 
\end{itemize}

\begin{algorithm}[h]
\caption{Recursive Importance Sketching for Rank Constrained Least Squares (RISRO)}
	\begin{algorithmic}[1]
		\STATE Input: $\mathcal{A}(\cdot): \mathbb{R}^{p_1\times p_2}\to \mathbb{R}^n$, $\y\in \mathbb{R}^n$, rank $r$, initialization $\X^0$ which admits singular value decomposition $\U^0 \bSigma^0 \V^{0\top}$, where $\U^0 \in \bbO_{p_1,r}, \V^0 \in \bbO_{p_2, r}, \bSigma^0 \in \bbR^{r \times r}$
		\FOR{$t=0, 1, \ldots$} 
		\STATE Perform importance sketching on $\cA$ and construct the covariates maps $\cA_B: \bbR^{r \times r} \to \bbR^n$, $\cA_{D_1}: \bbR^{(p_1 - r)\times r} \to \bbR^n$ and $\cA_{D_2}: \bbR^{r \times (p_2 -r)} \to \bbR^n$: for matrix ``$\cdot$", let
		\begin{equation*} 
		[\cA_B(\cdot)]_i = \langle \cdot, \U^{t\top} \A_i \V^{t}\rangle,\quad 
		[\cA_{D_1}(\cdot)]_i = \langle \cdot, \U_{\perp}^{t\top} \A_{i} \V^{t}\rangle,\quad
		[\cA_{D_2}(\cdot)]_i = \langle \cdot, \U^{t\top} \A_{i} \V_{\perp}^{t}\rangle,		\quad i=1,\ldots, n.
		\end{equation*}  
		\STATE Solve the unconstrained least squares problem
		\begin{equation} \label{eq: alg2 least square}
		(\B^{t+1}, \D_1^{t+1}, \D_2^{t+1}) = \argmin_{\substack{\B\in \bbR^{r\times r}, \D_i \in \bbR^{(p_i -r) \times r}, i=1,2}} \left\|\y - \cA_B(\B) - \cA_{D_1}(\D_1) - \cA_{D_2} (\D_2^\top) \right\|_2^2
		\end{equation}
		\STATE Compute $\X^{t+1}_{U} = \left(\U^{t}\B^{t+1} + \U^{t}_{\perp}\D_1^{t+1}\right)$ and $\X^{t+1}_V = \left(\V^{t}\B^{t+1\top} + \V^{t}_{\perp}\D_2^{t+1}\right)$.
		\STATE \label{alg:RISRO:st6} Perform QR orthogonalization: $\U^{t+1} = \QR(\X^{t+1}_U),\quad \V^{t+1} = \QR(\X^{t+1}_V).$
		\STATE \label{alg:RISRO:st7} Update $\X^{t+1} = \X^{t+1}_U \left(\B^{t+1}\right)^{\dagger} \X_V^{t+1\top}$.
		\ENDFOR
	\end{algorithmic} \label{alg: recursive IS alg 2}
\end{algorithm}

We give a high-level explanation of RISRO through a decomposition of $\y_i$. Suppose $\y_i = \langle \A_i, \widebar{\X} \rangle + \widebar{\bepsilon}_i$ where $\widebar{\X}$ is a rank $r$ target matrix with singular value decomposition $\widebar{\U}\widebar{\bSigma} \widebar{\V}^\top$ with $\widebar \U \in \bbO_{p_1, r}$, $\widebar \bSigma \in \bbR^{r\times r}$ and $\widebar \V \in \bbO_{p_2,r}$. Then
\begin{equation} \label{eq: sketching interp of Simul Mini}
\begin{split}
	\y_i = & \langle \U^{t\top} \A_i \V^t, \U^{t\top} \widebar{\X} \V^t \rangle + \langle \U^{t\top}_\perp \A_i \V^t, \U^{t\top}_\perp \widebar{\X} \V^t \rangle  + \langle \U^{t\top} \A_i \V^t_\perp, \U^{t\top} \widebar{\X} \V^t_\perp \rangle + \langle \U^{t\top}_{\perp}\A_i\V^t_{\perp}, \U^{t\top}_{\perp}\widebar{\X}\V^t_{\perp}\rangle + \widebar{\bepsilon}_i \\
	:= &\langle \U^{t\top} \A_i \V^t, \U^{t\top} \widebar{\X} \V^t \rangle + \langle \U^{t\top}_\perp \A_i \V^t, \U^{t\top}_\perp \widebar{\X} \V^t \rangle + \langle \U^{t\top} \A_i \V^t_\perp, \U^{t\top} \widebar{\X} \V^t_\perp \rangle + \bepsilon^t_i.
\end{split}
\end{equation}
Here, $\bepsilon^t:= \cA(P_{\U^t_\perp} \widebar{\X} P_{\V^t_\perp})+\widebar{\bepsilon}\in \bbR^n$ can be seen as the residual of the new regression model \eqref{eq: sketching interp of Simul Mini}, and $\U^{t\top}\A_i\V^t, \U^{t\top}_\perp \A_i \V^t, \U^{t\top}\A_i\V^t_\perp$ are exactly the importance covariates constructed in \eqref{eq: importance sketches}.
Let 
\begin{equation} \label{eq: tildeB,D}
  \widetilde{\B}^t := \U^{t\top} \widebar{\X} \V^t, \widetilde{\D}_1^t := \U^{t\top}_\perp \widebar{\X} \V^t, \widetilde{\D}_2^{t\top} := \U^{t\top} \widebar{\X} \V^{t}_\perp.
\end{equation}
If $\bepsilon^t = 0$, we have $(\widetilde{\B}^t,\widetilde{\D}_1^t,\widetilde{\D}_2^{t})$ is a solution of the least squares in \eqref{eq: alg2 least square}. Hence, we could set $\B^{t+1} = \widetilde{\B}^t, \D_1^{t+1} = \widetilde{\D}_1^t, \D_2^{t+1} = \widetilde{\D}_2^{t}$ and thus $ \X_U^{t+1} = \widebar{\X}\V^t, \X_V^{t+1} = \widebar{\X}^{\top} \U^t$. 
Furthermore, if $\B^{t+1}$ is invertible, then it holds that
\begin{equation}\label{eq: B, D1, D2 update equality equation}
\begin{split}
 \X_U^{t+1} (\B^{t+1})^{-1} \X_V^{t+1 \top} = \widebar{\X} \V^t (\U^{t \top} \widebar{\X} \V^t)^{-1} (\widebar{\X}^{\top} \U^t)^\top = \widebar{\X},
\end{split}
\end{equation}
which means $\widebar{\X}$ can be exactly recovered by one iteration of RISRO.

In general, $\bepsilon^t \neq 0$. When the column spans of $\U^t, \V^t$ well approximate the ones of $\widebar{\U}, \widebar{\V}$, i.e., the column and row subspaces that the target parameter $\widebar{\X}$ lie on, we expect $\U^{t\top}_{\perp}\widebar{\X}\V^t_{\perp}$ and $\bepsilon^t_i = \langle\U_{\perp}^{t\top}\A_i\V_{\perp}^t, \U_{\perp}^{t\top}\widebar{\X}\V_{\perp}^t\rangle + \widebar{\bepsilon}_i$ to have a small amplitude, then $\B^{t+1}, \D_1^{t+1}, \D_2^{t+1}$, the outcome of the least squares problem \eqref{eq: alg2 least square}, can well approximate $\widetilde{\B}^t, \widetilde{\D}_1^t, \widetilde{\D}_2^{t}$. In Lemma \ref{lm: partial least square for Simul Mini}, we give a precise characterization for this approximation. Before that, let us introduce a convenient notation so that \eqref{eq: alg2 least square} can be written in a more compact way. 

Define the linear operator ${\cal L}_t$ as
\begin{equation} \label{eq: linear operator L_t}
\begin{split}
		&{\cal L}_t: \W = \begin{bmatrix}
\W_0 \in \bbR^{r \times r} & \W_2 \in \bbR^{r \times (p_2 - r)}\\
\W_1 \in \bbR^{(p_1-r) \times r} & \0_{(p_1-r)\times (p_2 - r)} 
\end{bmatrix} \to [\U^t \quad \U^t_\perp] \begin{bmatrix}
\W_0  & \W_2 \\
\W_1  & \0
\end{bmatrix} [\V^t \quad \V^t_\perp]^\top,
\end{split}
\end{equation}
and it is easy to compute its adjoint ${\cal L}^*_t: \M \in \bbR^{p_1 \times p_2} \to \begin{bmatrix}
\U^{t\top} \M \V^t  & \U^{t\top} \M \V^t_\perp\\
(\U^t_\perp)^\top \M \V^t  & \0
\end{bmatrix}$. Then, the least squares problem in \eqref{eq: alg2 least square} can be written as
\begin{equation} \label{eq: compact repre of least square}
	(\B^{t+1}, \D_1^{t+1}, \D_2^{t+1}) = \argmin_{\substack{\B\in \bbR^{r\times r}, \D_i \in \bbR^{(p_i -r) \times r}, i=1,2}} \left\|\y - \cA \cL_t \left( \begin{bmatrix}
	\B  & \D_2^\top \\
	\D_1  & \0
	\end{bmatrix} \right) \right\|_2^2.
\end{equation}
\begin{Lemma}[Iteration Error Analysis for RISRO] \label{lm: partial least square for Simul Mini}
	Let $\widebar{\X}$ be any given target matrix. Recall the definition of $\bepsilon^t = \widebar{\bepsilon} + \cA(P_{\U^t_\perp} \widebar\X P_{\V^t_\perp})$ from \eqref{eq: sketching interp of Simul Mini}. If the operator $\cL_t^* \cA^* \cA \cL_t$ is invertible over $\Range(\cL_t^*)$, then $\B^{t+1}, \D_1^{t+1}, \D_2^{t+1}$ in \eqref{eq: alg2 least square} 
	satisfy
\begin{equation}\label{eq: Bt+1 - Btilde}
\begin{split}
    	\begin{bmatrix}
		\B^{t+1} - \widetilde{\B}^t & \D_2^{t+1 \top} - \widetilde{\D}_2^{t\top}\\
	\D_1^{t+1} - \widetilde{\D}_1^t & \0
	\end{bmatrix} = (\cL^*_t \cA^* \cA \cL_t)^{-1} \cL_t^* \cA^* \bepsilon^t,
\end{split}
\end{equation} and 
\begin{equation} \label{eq: B and D part bound}
	\|\B^{t+1} - \widetilde{\B}^t \|_\F^2 + \sum_{k=1}^2 \|\D_k^{t+1} - \widetilde{\D}_k^t\|_\F^2 = \left\|(\cL^*_t \cA^* \cA \cL_t)^{-1} \cL_t^* \cA^* \bepsilon^t\right\|_\F^2.
\end{equation}
\end{Lemma}

In view of Lemma \ref{lm: partial least square for Simul Mini}, the approximation errors of $\B^{t+1}, \D_1^{t+1}, \D_2^{t+1}$ to $\widetilde{\B}^t, \widetilde{\D}_1^t, \widetilde{\D}_2^t$ are driven by the least squares residual $\|(\cL_t^* \cA^* \cA \cL_t)^{-1} \cL_t^*\cA^* \bepsilon^t \|_\F^2$. This fact plays a key role in the proof for the high-order convergence theory of RISRO, see later in Remark \ref{rm: proof-sketch}.

\begin{Remark}[Comparison with Randomized Sketching]
The importance sketching in RISRO is significantly different from the randomized sketching in the literature (see surveys \cite{mahoney2011randomized,woodruff2014sketching} and the references therein). The randomized sketching matrices are often randomly generated and reduce the sample size ($n$), the importance sketching matrices are deterministically constructed under the supervision of $\y$ and reduce the dimension of parameter space ($p_1p_2$). See \cite[Section 1.3 and 2]{zhang2020islet} for more comparison of randomized and importance sketchings. 
\end{Remark}

\subsection{Comparison with More Algorithms in the View of Sketching} \label{sec: sketching interpre for more algorithms}

In addition to RISRO, several classic algorithms for rank constrained least squares can be interpreted from the recursive importance sketching perspective. Through the lens of the sketching, RISRO exhibits advantages over these existing algorithms.

We first focus on Alternating Minimization (Alter Mini) proposed and studied in \cite{hardt2014understanding,jain2013low,zhao2015nonconvex}. Suppose $\U^t$ is the left singular vectors of $\X^t$, the outcome of the $t$-th iteration, Alter Mini solves the following least squares problems to update $\U$ and $\V$,
\begin{equation}\label{eq: sketching point view of Alter Mini}
\begin{split}
	& \widecheck{\V}^{t+1} = \argmin_{\V \in \bbR^{p_2 \times r}} \sum_{i=1}^n \left( \y_i - \langle \A_i, \U^t \V^\top \rangle \right)^2  = \argmin_{\V\in \bbR^{p_2 \times r}} \sum_{i=1}^n \left( \y_i - \langle \U^{t \top} \A_i,  \V^\top \rangle \right)^2,\\
	& \widecheck{\U}^{t+1} = \argmin_{\U \in \bbR^{p_1 \times r}} \sum_{i=1}^n\left(\y - \langle\A_i, \U(\V^{t+1})^\top\rangle\right)^2 = \argmin_{\U \in \bbR^{p_1 \times r}} \sum_{i=1}^n\left(\y - \langle\A_i\V^{t+1}, \U\rangle\right)^2,\\
	& \V^{t+1} = \QR(\widecheck{\V}^{t+1}), \quad \U^{t+1} = \QR(\widecheck{\U}^{t+1}).
\end{split}
\end{equation}
Then, Alter Mini essentially solves least squares problems with sketched covariates $\U^{t\top}\A_i, \A_i\V^{t+1}$ to update $\widecheck{\V}^{t+1}, \widecheck{\U}^{t+1}$ alternatively and iteratively. The number of parameters of the least squares in \eqref{eq: sketching point view of Alter Mini} are $rp_2$ and $rp_1$ as opposed to $p_1p_2$, the number of parameters in the original least squares problem. See Figure \ref{fig: sketching illutstration} upper right panel for an illustration of the sketching scheme in Alter Mini.
Consider the following decomposition of $\y_i$,
\begin{equation}\label{eq: intuition of Alter Mini}
\begin{split}
	\y_i & = \langle \A_i, P_{\U^t} \widebar{\X} \rangle + \langle \A_i, P_{\U^t_\perp} \widebar{\X} \rangle +\widebar{\bepsilon}_i = \langle \U^{t\top} \A_i,  \U^{t\top} \widebar{\X} \rangle + \langle \A_i, P_{\U^t_\perp} \widebar{\X} \rangle + \widebar{\bepsilon}_i :=  \langle \U^{t\top} \A_i,  \U^{t\top} \widebar{\X} \rangle + \widecheck{\bepsilon}^t_i,
\end{split}
\end{equation}
where $\widecheck{\bepsilon}^t:= \cA(P_{\U^t_\perp}\widebar{\X})+\bar\bepsilon\in \bbR^n$. Define $\widecheck{\A}^t \in \bbR^{n \times p_2 r}$ with $\widecheck{\A}_{[i,:]} = \rmvec(\U^{t\top} \A_i)$. Similar to how Lemma \ref{lm: partial least square for Simul Mini} is proved, we can show $\|\widecheck{\V}^{t+1 \top} - \U^{t\top} \widebar{\X}\|^2_\F = \|(\widecheck{\A}^{t\top} \widecheck{\A}^t)^{-1} \widecheck{\A}^{t\top} \widecheck{\bepsilon}^t\|_2^2,$
which implies the approximation error of $\V^{t+1}=\QR(\widecheck{\V}^{t+1})$ (i.e., the outcome of one iteration Alter Mini) to $\widebar{\V}$ (i.e., true row span of the target matrix $\widebar{\X}$) is driven by $\widecheck{\bepsilon}^t = \cA(P_{\U^t_\perp}\widebar{\X})+\bar\bepsilon$, i.e., the residual of least squares problem \eqref{eq: intuition of Alter Mini}. Recall for RISRO, Lemma \ref{lm: partial least square for Simul Mini} shows the approximation error of $\V^{t+1}$ is driven by $\bepsilon^t = \cA(P_{\U^t_\perp} \widebar{\X} P_{\V^t_\perp})+\bar\bepsilon$. Since $\|P_{\U_{\perp}^t}\widebar{\X} P_{\V_{\perp}^t}\|_\F \leq \|P_{\U_{\perp}^t}\widebar{\X}\|_\F$, the approximation error in per iteration of RISRO can be smaller than the one of Alter Mini. Such a difference between RISRO and Alter Mini is due to the following fact: in Alter Mini, the sketching captures the importance covariates correspond to only the row (or column) span of $\X^t$ in updating $\V^{t+1}$ (or $\U^{t+1}$), while the importance sketching of RISRO in \eqref{eq: importance sketches} catches the importance covariates from both the row span and column span of $\X^t$. As a consequence, Alter Mini iterations yield first-order convergence while RISRO iterations render high-order convergence as will be established in Section \ref{sec:theory}.

\begin{Remark}
	Recently, \cite{kummerle2018harmonic} proposed a harmonic mean iterative reweighted least squares (HM-IRLS) method for low-rank matrix recovery: they specifically solve $\min_{\X \in \bbR^{p_1 \times p_2}} \|\X\|_q^q$ subject to $\y = \cA(\X)$, where $\|\X\|_q = (\sum_i \sigma^q_i(\X))^{1/q}$ is the Schatten-$q$ norm of the matrix $\X$. Compared to the original iterative reweighted least squares (IRLS) \citep{fornasier2011low,mohan2012iterative}, which only involves either the column span or the row span of $\X^t$ in constructing the reweighting matrix, HM-IRLS leverages both the column and row spans of $\X^t$ in constructing the reweighting matrix per-iteration and performs better. Such a comparison of HM-IRLS versus IRLS shares the same spirit as RISRO versus Alter Mini: the importance sketching of RISRO simultaneously captures the information of both column and row spans of $\X^t$ per iteration and achieves a better performance. Utilizing both row and column spans of $\X^t$ simultaneously is the key to achieve high-order convergence performance by RISRO. 
\end{Remark}

Another example is the rank $2r$ iterative least squares (R2RILS) proposed in \cite{bauch2020rank} for solving ill-conditioned matrix completion problems. In particular, at the $t$-th iteration, Step 1 of R2RILS solves the following least squares problem
\begin{equation} \label{eq: R2RILS least square step}
	\min_{\M \in \bbR^{p_1 \times r}, \N \in \bbR^{p_2 \times r}} \sum_{(i,j) \in \Omega} \left\{\left(\U^t \N^\top + \M \V^{t \top} - \X\right)_{[i,j]} \right\}^2,
\end{equation} 
where $\Omega$ is the set of index pairs of the observed entries. In the matrix completion setting, it turns out the following equivalence holds (proof given in Appendix)
\begin{equation} \label{eq: sketching in R2RILS}
	\begin{split}
		&\argmin_{\M \in \bbR^{p_1 \times r}, \N \in \bbR^{p_2 \times r}} \sum_{(i,j) \in \Omega} \left\{\left(\U^t \N^\top + \M \V^{t \top} - \X\right)_{[i,j]} \right\}^2\\
		 = &\argmin_{\M \in \bbR^{p_1 \times r}, \N \in \bbR^{p_2 \times r}} \sum_{(i,j) \in \Omega} \left( \langle \U^{t\top} \A^{ij} , \N^\top \rangle +  \langle \M, \A^{ij} \V^t \rangle  - \X_{[i,j]} \right)^2,
	\end{split}
\end{equation} 
where $\A^{ij} \in \bbR^{p_1 \times p_2}$ is the special covariate in matrix completion satisfying $(\A^{ij})_{[k,l]}=1$ if $(i,j)=(k,l)$ and $(\A^{ij})_{[k,l]}=0$ otherwise. This equivalence reveals that the least squares step \eqref{eq: R2RILS least square step} in R2RILS can be seen as an implicit sketched least squares problem similar to \eqref{eq: alg2 least square} and \eqref{eq: sketching point view of Alter Mini} with covariates $\U^{t\top} \A^{ij}$ and $\A^{ij} \V^t$ for $(i,j) \in \Omega$. 

We give a pictorial illustration for the sketching interpretation of R2RILS on the bottom right part of Figure \ref{fig: sketching illutstration}. Different from the sketching in RISRO, R2RILS incorporates the core sketch $\U^{t\top} \A_i \V^t$ twice, which results in the rank deficiency in the least squares problem \eqref{eq: R2RILS least square step} and brings difficulties in both implementation and theoretical analysis. RISRO overcomes this issue by performing a better-designed sketching and covers more general low-rank matrix recovery settings than R2RILS. With the new sketching scheme, we are able to give a new and solid theory for RISRO with high-order convergence. 

\section{Theoretical Analysis}\label{sec:theory}

In this section, we provide convergence analysis for the proposed algorithm. For technical convenience, we assume $\cA$ satisfies the Restricted Isometry Property (RIP) \citep{candes2008restricted}. The RIP condition, first introduced in compressed sensing, has been widely used as one of the most standard assumptions in the low-rank matrix recovery literature \citep{cai2013sharp,cai2014sparse,candes2011tight,chen2015fast,jain2010guaranteed,recht2010guaranteed,tu2016low,zhao2015nonconvex}. It also plays a critical role in analyzing the landscape of the rank constrained optimization problem \eqref{eq:minimization}  \citep{bhojanapalli2016global,ge2017no,uschmajew2018critical,zhang2019sharp,zhu2018global}. On the other hand, RIP is only a sufficient but not necessary condition for the convergence of RISRO. We will illustrate later in several examples that RISRO converges quadratically while RIP completely fails. 
\begin{Definition}[Restricted Isometry Property (RIP)] \label{def: RIP} Let $\cA : \bbR^{p_1 \times p_2} \to \bbR^n$ be a linear map. For every integer $r$ with $1 \leq r \leq \min(p_1, p_2)$, define the $r$-restricted isometry constant to be the smallest number $R_r$ such that $(1-R_r) \|\Z\|^2_\F \leq \|\cA(\Z)\|_2^2 \leq (1+R_r) \|\Z\|_\F^2$ holds for all $\Z$ of rank at most $r$. And $\cA$ is said to satisfy the $r$-restricted isometry property ($r-$RIP) if $ 0\leq R_r < 1$.
\end{Definition}

 The RIP condition provably holds when $\cA$ has independent random sub-Gaussian design or $\cA$ is a random projection \citep{candes2011tight,recht2010guaranteed}. In addition, the definition of RIP above can be equivalently stated in a matrix format: define $\tilde{\A} = [\rmvec(\A_1), \cdots, \rmvec(\A_n) ]^\top$ and $\cA(\Z) = \tilde{\A}\rmvec(\Z)$. Then $\cA$ satisfies the RIP condition is equivalent to $(1-R_r)\|\rmvec(\Z)\|_2^2 \leq \|\tilde{\A}(\rmvec(\Z))\|_2^2 \leq (1+R_r)\|\rmvec(\Z)\|_2^2$ for all matrices $\Z$ of rank at most $r$. By definition, $R_r \leq R_{r'}$ for any $r \leq r'$.

By assuming RIP for $\cA$, we can show the linear operator $\cL_t^* \cA^* \cA \cL_t$ mentioned in Lemma \ref{lm: partial least square for Simul Mini} is always invertible over $\Range({\cal L}_t^*)$ (i.e. the least squares \eqref{eq: alg2 least square} has a unique solution). The following lemma gives explicit lower and upper bounds for the spectrum of this operator.
\begin{Lemma}[Bounds for Spectrum of $\cL_t^* \cA^* \cA \cL_t$] \label{lm: spectral norm bound of Atop A}
	Recall the definition of $\cL_t$ in \eqref{eq: linear operator L_t}. It holds that 
	\begin{equation}
	    \label{eq: norm equal Range Lt}
	    \|\cL_t(\M)\|_\F = \|\M\|_\F, \quad \forall \, \M\in \Range(\cL^*_t).
	\end{equation}
	Suppose the linear map $\cA$ satisfies the 2r-RIP. Then, it holds that for any matrix $\M\in \Range(\cL_t^*)$,
	\begin{equation*}
		(1 - R_{2r}) \|\M\|_\F \leq \|\cL_t^* \cA^* \cA \cL_t(\M) \|_\F \leq (1+R_{2r})\|\M\|_\F.
	\end{equation*}
\end{Lemma}
\begin{Remark}[Bounds for spectrum of $(\cL_t^* \cA^* \cA \cL_t)^{-1}$] By the relationship of the spectrum of an operator and its inverse, from Lemma \ref{lm: spectral norm bound of Atop A}, we also have the spectrum of $(\cL_t^* \cA^* \cA \cL_t)^{-1}$ is lower and upper bounded by $\frac{1}{(1+R_{2r})}$ and $\frac{1}{(1-R_{2r})}$, respectively.
\end{Remark}

 In the following Proposition \ref{prop: bound for iter approx error}, we bound the iteration approximation error given in Lemma \ref{lm: partial least square for Simul Mini}. 

\begin{Proposition}[Upper Bound for Iteration Approximation Error] \label{prop: bound for iter approx error}
	Let $\widebar{\X}$ be a given target rank r matrix and $\widebar{\bepsilon} = \y - \cA(\widebar{\X})$. Suppose that $\cA$ satisfies the $2r$-RIP. Then at $t$-th iteration of RISRO, the approximation error \eqref{eq: B and D part bound} has the following upper bound:
	\begin{equation} \label{ineq: prop least square residual bound}
	\begin{split}
		&\left\|(\cL^*_t \cA^* \cA \cL_t)^{-1} \cL_t^* \cA^* \bepsilon^t\right\|_\F^2\\
		 \leq& \frac{R_{3r}^2 \|\X^t - \widebar{\X} \|^2 \|\X^t - \widebar{\X} \|_\F^2 }{(1-R_{2r})^2\sigma_r^2(\widebar{\X})} + \frac{\| \cL_t^* \cA^* (\widebar{\bepsilon})\|_\F^2}{(1-R_{2r})^2} + \| \cL_t^* \cA^* (\widebar{\bepsilon})\|_\F \frac{2R_{3r} \|\X^t - \widebar{\X} \| \|\X^t - \widebar{\X} \|_\F }{\sigma_r(\widebar{\X})(1-R_{2r})^2}.
	\end{split}
\end{equation}   
\end{Proposition}

Note that Proposition \ref{prop: bound for iter approx error} is rather general in the sense that it applies to any $\widebar{\X}$ of rank $r$ and we will pick different choices of $\widebar{\X}$ depending on our purposes. For example, in studying the convergence of RISRO, e.g., the upcoming Theorem \ref{th: local contraction general setting}, we treat $\widebar{\X}$ as a stationary point and in the setting of estimating the model parameter in matrix trace regression, we take $\widebar{\X}$ to be the ground truth (see Theorem \ref{th: local convergence in local rank matrix recovery}).

Now, we are ready to establish the deterministic convergence theory for RISRO.
For problem \eqref{eq:minimization}, we use the following definition of stationary points: a rank $r$ matrix $\widebar\X$ is said to be a stationary point of \eqref{eq:minimization} if $ \nabla f(\widebar{\X})^\top \widebar{\U} =0$ and $\nabla f(\widebar{\X}) \widebar{\V} = 0$ where $\nabla f(\widebar{\X}) = \cA^*( \cA(\widebar{\X}) - \y)$, and $\widebar{\U}, \widebar{\V}$ are the left and right singular vectors of $\widebar{\X}$. See also \cite{ha2020equivalence}. In Theorem \ref{th: local contraction general setting}, we show that given any target stationary point $\widebar{\X}$ and proper initialization, RISRO has a local quadratic-linear convergence rate in general and quadratic convergence rate if $\y = \mathcal{A}(\bar{\X})$. 

\begin{Theorem}[Local Quadratic-Linear and Quadratic Convergence of RISRO]\label{th: local contraction general setting}
Let $\widebar{\X}$ be a stationary point to problem \eqref{eq:minimization} and $\widebar{\bepsilon} = \y - \cA(\widebar{\X})$. Suppose that $\cA$ satisfies the $2r$-RIP, and the initialization $\X^0$ satisfies
\begin{equation} \label{ineq: optimization theory initialization condition}
	\|\X^0 - \widebar{\X}\|_\F \leq \left( \frac{1}{4} \wedge \frac{1-R_{2r}}{4\sqrt{5}R_{3r}}  \right) \sigma_r(\widebar{\X}), 
\end{equation} 
and $\|\cA^*(\widebar{\bepsilon})\|_\F \leq \frac{1-R_{2r}}{4\sqrt{5}}\sigma_r(\widebar{\X})$. Then, we have $\{\X^t\}$, the sequence generated by RISRO (Algorithm \ref{alg: recursive IS alg 2}), converges linearly to $\widebar{\X}$: $\|\X^{t+1} - \widebar{\X} \|_\F \leq \frac{3}{4} \|\X^{t} - \widebar{\X}\|_\F, \quad \forall\, t\ge 0.$

More precisely, it holds that $\forall \, t \ge 0$:
\begin{equation}\label{ineq: optimization theory quadratic and linear converge}
\begin{split}
	\|\X^{t+1} - & \widebar{\X}\|^2_\F \leq \frac{5\|\X^t - \widebar{\X}\|^2}{(1-R_{2r})^2 \sigma^2_r(\widebar{\X})}\cdot \left( R_{3r}^2 \|\X^t - \widebar{\X}\|_\F^2 + 4 R_{3r} \|\cA^*(\widebar{\bepsilon})\|_\F\|\X^t - \widebar{\X}\|_\F + 4\|\cA^*(\widebar{\bepsilon})\|^2_\F   \right).
\end{split}
\end{equation}

In particular, if $\widebar{\bepsilon} = 0$, then $\{\X^t\}$ converges quadratically to $\widebar{\X}$ as 
\begin{equation*}
	\|\X^{t+1} - \widebar{\X} \|_\F \leq \frac{\sqrt{5}R_{3r}}{(1-R_{2r})\sigma_r(\widebar{\X}) } \|\X^t - \widebar{\X}\|_\F^2, \quad \forall\, t \ge 0.
\end{equation*}
\end{Theorem}

\begin{Remark}[Quadratic-linear and Quadratic Convergence of RISRO]\label{rem: quadratic-linear convergence}
	We call the convergence in \eqref{ineq: optimization theory quadratic and linear converge} quadratic-linear since the sequence $\{\X^t\}$ generated by RISRO exhibits a phase transition from quadratic to linear convergence: when $\|\X^t - \widebar{\X}\|_\F \gg \|\cA^*(\widebar{\bepsilon})\|_\F$, the algorithm has a quadratic convergence rate; when $\X^{t}$ becomes close to $\widebar{\X}$ such that $\|\X^t - \widebar{\X}\|_\F \leq c \|\cA^*(\widebar{\bepsilon})\|_\F$ for some $c > 0$, the convergence rate becomes linear. Even though the ultimate convergence of RISRO is linear {\rm to a stationary point} in the noisy setting, we will show later in Section \ref{sec: low rank matrix regression} that RISRO achieves quadratic convergence {\rm in estimating the underlying parameter matrix} in statistical applications. Moreover, as $\widebar{\bepsilon}$ becomes smaller, the stage of quadratic convergence becomes longer (see Section \ref{sec: numerical Simul Mini property} for a numerical illustration of this convergence pattern). In the extreme case $\widebar{\bepsilon} = 0$, Theorem \ref{th: local contraction general setting} covers the widely studied matrix sensing problem under the RIP framework \citep{chen2015fast,jain2010guaranteed,park2018finding,recht2010guaranteed,tu2016low,zhao2015nonconvex,zheng2015convergent}. It shows as long as the initialization error is within a constant factor of $\sigma_r(\widebar{\X})$, RISRO enjoys quadratic convergence to the target matrix $\widebar{\X}$. To the best of our knowledge, we are among the first to give quadratic-linear algorithmic convergence guarantees for general rank constrained least squares and quadratic convergence for matrix sensing. Recently, \cite{charisopoulos2019low} formulated \eqref{eq:minimization} as a non-convex composite optimization problem based on $\X = \R \L^\top$ factorization and showed that the prox-linear algorithm \citep{burke1985descent,lewis2016proximal} achieves local quadratic convergence when $\bar{\bepsilon}=0$. In each iteration therein, a carefully tuned convex program needs to be solved exactly and the tuning parameter relies on the unknown weakly convexity parameter of the composite objective function. In contrast, the proposed RISRO is tuning-free, only solves a dimension-reduced least squares in each step, and can be as cheap as many first-order methods. See Section \ref{sec: computation complexity} for a detailed discussion on the computational complexity of RISRO.

	Moreover, a quadratic-linear convergence rate also appears in several other methods under different settings: \cite{pilanci2017newton} studied the local convergence of the randomized Newton Sketch for objectives with strong convexity and smooth properties; \cite{erdogdu2015convergence} considered the sub-sampled Newton method to optimize an objective function in the form of a sum of convex functions and established their convergence theory with the well-conditioned sub-sampled Hessian. We consider the non-convex matrix optimization problem \eqref{eq:minimization} and use the recursive importance sketching method. Our quadratic-linear convergence result can be boosted to quadratic when $\bar{\bepsilon}=0$.
	
\end{Remark}

\begin{Remark}[Initialization] \label{rem: globalization and initialization}
	The convergence theory in Theorem \ref{th: local contraction general setting} requires a good initialization condition. Practically, the spectral method often provides a sufficiently good initialization that meets the requirement in \eqref{ineq: optimization theory initialization condition} in many statistical applications. In Section \ref{sec: statistics applications} and \ref{sec:numerics}, we will illustrate this point from two applications: matrix trace regression and phase retrieval. 
\end{Remark}

\begin{Remark}[Small residual condition in Theorem \ref{th: local contraction general setting}]
    In addition to the initialization condition, the small residual condition $\|\cA^*(\widebar{\bepsilon})\|_\F \leq \frac{1-R_{2r}}{4\sqrt{5}}\sigma_r(\widebar{\X})$ is also needed in Theorem \ref{th: local contraction general setting}. This condition essentially means that the signal strength at point $\widebar{\X}$ needs to dominate the noise. If $\widebar{\bepsilon} = \y - \cA(\widebar{\X})= 0$, then the aforementioned small residual condition holds automatically.
\end{Remark}

\begin{Remark}\label{rm: proof-sketch}We provide a proof sketch of Theorem \ref{th: local contraction general setting} and discuss our technical contributions therein. 

\noindent{\bf Step 1}. 
We bound $\| \cL_t^* \cA^* (\widebar{\bepsilon})\|_\F \leq \frac{4\|\X^t - \widebar{\X} \|^2}{\sigma_r^2(\widebar{\X})} \|\cA^*(\widebar{\bepsilon})\|_\F^2$, and then apply Proposition \ref{prop: bound for iter approx error} to obtain an upper bound for the approximation error in \eqref{eq: B and D part bound}: 
\begin{equation} \label{ineq: least square residual upper bound2}
	\begin{split}
	\|(\cL^*_t \cA^* & \cA \cL_t)^{-1} \cL_t^* \cA^* \bepsilon^t\|_\F^2 \leq  \frac{\|\X^t - \widebar{\X}\|^2}{(1-R_{2r})^2 \sigma^2_r(\widebar{\X})}\\
	&\quad \quad \quad \cdot \left( R_{3r}^2 \|\X^t - \widebar{\X}\|_\F^2 + 4 \|\cA^*(\widebar{\bepsilon})\|_\F^2 + 4 R_{3r}\|\cA^*(\widebar{\bepsilon})\|_\F \|\X^t - \widebar{\X}\|_\F \right).
	\end{split}
\end{equation}

\noindent{\bf Step 2}. We use induction to show the following three claims,
\begin{equation*}
	\begin{split}
		 &(C1)~ \max\{ \left\| \sin \Theta(\U^{t}, \widebar{\U})  \right\|, \left\| \sin \Theta(\V^{t},\widebar{\V})  \right\| \} \leq \frac{1}{2}; \quad (C2)~  \B^{t+1} \text{ in \eqref{eq: alg2 least square} is invertible} ;\\
		& (C3)~ 	\|\X^{t+1} - \widebar{\X}\|_\F^2 \leq \frac{5\|\X^t - \widebar{\X}\|^2}{(1-R_{2r}^2) \sigma^2_r(\widebar{\X})} \left( R_{3r}^2 \|\X^t - \widebar{\X}\|_\F^2 + 4 R_{3r} \|\cA^*(\widebar{\bepsilon})\|_\F\|\X^t - \widebar{\X}\|_\F + 4\|\cA^*(\widebar{\bepsilon})\|^2_\F   \right)\\
	& \quad \quad \quad \quad \quad \quad \quad  \leq \frac{9}{16} \|\X^t - \widebar{\X} \|_\F^2.
	\end{split}
\end{equation*}
Here, (C2) means the iterates $\X^t$ are always rank $r$. This fact is useful in Section \ref{sec: Riemannian manifold interpre} in connecting RISRO to Riemannian optimization on fixed rank matrix manifolds. (C2) is proved by (C1) and Lemma \ref{lm: partial least square for Simul Mini}. In proving (C3), we introduce an intermediate quantity $ \rho_{t+1} = \max \{ \left\| \D_1^{t+1} (\B^{t+1})^{-1} \right\|, \left\| (\B^{t+1})^{-1} \D_2^{t+1\top} \right\|\}$ and obtain
\begin{equation} \label{ineq: th1-key-argument}
	\begin{split}
		\|\X^{t+1} - \widebar{\X}\|_\F^2 = \left\| \begin{array}{c c}
	\B^{t+1} - \widetilde{\B}^t & \D_2^{t+1 \top} - \widetilde{\D}_2^{t\top}\\
	\D_1^{t+1} - \widetilde{\D}_1^t & 
	~~ \D_1^{t+1} (\B^{t+1})^{-1} \D_2^{t+1\top}- \widetilde{\D}_1^{t} (\widetilde{\B}^{t})^{-1} \widetilde{\D}_2^{t\top}
	\end{array} \right\|_\F^2 \overset{(a)} \leq 5  \left\|(\cL^*_t \cA^* \cA \cL_t)^{-1} \cL_t^* \cA^* \bepsilon^t\right\|_\F^2,
	\end{split}
\end{equation} 
Here (a) is by the induction assumptions, Lemma \ref{lm: partial least square for Simul Mini}, and Lemma \ref{lm:FGH}. Finally, (C3) follows by plugging \eqref{ineq: least square residual upper bound2} into \eqref{ineq: th1-key-argument} and the induction assumptions and this proves the main result of Theorem \ref{th: local contraction general setting}.

\end{Remark}

\section{A Riemannian Manifold Optimization Interpretation of RISRO}\label{sec: Riemannian manifold interpre}

The superior performance of RISRO yields the following question: 

{\it Is there a connection of RISRO to any class of optimization algorithms in the literature?}

In this section, we give an affirmative answer to this question. We show RISRO can be viewed as a Riemannian optimization algorithm on the manifold $\cM_r:=\left\{ 
\X\in \bbR^{p_1\times p_2}\mid {\rm rank}(\X) = r
\right\}$. We find the sketched least squares in~\eqref{eq: alg2 least square} of RISRO actually solves the \emph{Fisher Scoring} or \emph{Riemannian Gauss-Newton equation} and Step~\ref{alg:RISRO:st7} in RISRO performs a type of \emph{retraction} under the framework of Riemannian optimization. 

Riemannian optimization concerns optimizing a real-valued function $f$ defined on a Riemannian manifold $\mathcal{M}$. One commonly-encountered manifold is a submanifold of $\mathbb{R}^n$. Under such circumstances, a manifold can be viewed as a smooth subset of $\mathbb{R}^n$. When a smooth-varying inner product is further defined on the subset, the subset together with the inner product is called a Riemannian manifold. We refer to~\cite{absil2009optimization} for the rigorous definition of Riemannian manifolds. Optimization on a Riemannian manifold often relies on the notion of Riemannian gradient/Riemannian Hessian (which are used to find a search direction) and the notion of retraction (which is defined for the motion of iterates on the manifold). The remainder of this section describes the required Riemannian optimization tools and the connection of RISRO to Riemannian optimization.

It has been shown in~\cite[Example 8.14]{lee2013smooth} that the set $\mathcal{M}_r$ is a smooth submanifold of $\mathbb{R}^{p_1 \times p_2}$ and the tangent space is also given therein. The result is given in Proposition~\ref{Mrmanifold} for completeness.
\begin{Proposition}\cite[Example 8.14]{lee2013smooth} \label{Mrmanifold}
	$\cM_r = \{\X\in \mathbb{R}^{p_1\times p_2}:\rank(\X) = r\}$ is a smooth embedded submanifold of dimension $(p_1 + p_2 -r)r$. Its tangent space $T_{\X}\cM_r$ at $\X \in \cM_r$ with the SVD decomposition $\X = \U\bSigma \V^\top$ ($\U\in \bbO_{p_1, r}$ and $\V\in \bbO_{p_2, r}$) is given by:
	\begin{equation}
	\label{eq:tangentS}
	T_\X \cM_r = \left\{
	[\U\quad \U_{\perp}] \begin{bmatrix}
	\bbR^{r\times r} &\bbR^{r\times (p_2 - r)} \\[2pt]
	\bbR^{(p_1 - r)\times r} & \0_{(p_1 -r)\times (p_2 - r)}
	\end{bmatrix}
	[\V\quad \V_{\perp}]^\top
	\right\}.
	\end{equation}
\end{Proposition}
The Riemannian metric of $\mathcal{M}_r$ that we use throughout this paper is the Euclidean inner product, i.e., $\langle \U, \V \rangle = \mathrm{trace}(\U^\top \V)$.

In the Euclidean setting, the update formula in an iterative algorithm is $X^t + \alpha \eta^t$, where $\alpha$ is the stepsize and $\eta^t$ is a descent direction. However, in the framework of Riemannian optimization, $\X^t + \alpha \eta^t$ is generally neither well-defined nor lying in the manifold. To overcome this difficulty, the notion of retraction is used, see e.g.,~\cite{absil2009optimization}. Considering the manifold $\mathcal{M}_r$, we have the definition that a retraction $R$ is a smooth map from $T \mathcal{M}_r$ to $\mathcal{M}_r$ satisfying i) $R(\X, 0) = \X$ and ii) $\frac{d}{d t} R(\X, t \eta) \vert_{t = 0} = \eta$ for all $\X \in \mathcal{M}_r$ and $\eta \in T_\X \mathcal{M}_r$, where $T \mathcal{M}_r = \{(\X, T_\X \mathcal{M}_r) : \X \in \mathcal{M}_r\}$, is the tangent bundle of $\mathcal{M}_r$. The two conditions guarantee that $R(\X, t \eta)$ stays in $\mathcal{M}_r$ and $R(\X, t \eta)$ is a first-order approximation of $\X + t \eta$ at $t = 0$.

Next, we show that Step~\ref{alg:RISRO:st7} in Algorithm \ref{alg: recursive IS alg 2} performs the orthographic retraction on the manifold of fixed-rank matrices given in~\cite{absil2012projection}. Suppose at iteration $t+1$, $\B^{t+1}$ is invertible (this is true under the RIP framework, see Remark \ref{rm: proof-sketch} and Step 2 in the proof of Theorem \ref{th: local contraction general setting}). We can show by some algebraic calculations that the update $\X^{t+1}$ in Step~\ref{alg:RISRO:st7} can be rewritten as
\begin{equation} \label{eq: retraction step verification}
		\begin{split}
		\X^{t+1} &= \X^{t+1}_U \left(\B^{t+1}\right)^{-1} \X_V^{t+1\top} = [\U^t \quad \U^{t}_\perp]
\begin{bmatrix}
\B^{t+1} & \D_2^{t+1\top} \\
\D_1^{t+1} & \D_1^{t+1} {(\B^{t+1})}^{-1} \D_2^{t+1\top}
\end{bmatrix} [\V^t \quad \V^t_\perp]^{\top}.
	\end{split}
\end{equation}
Let $\eta^t \in T_{\X^t} \cM_r$ be the update direction and $\X^t + \eta^t$ has the following representation, 
\begin{equation} \label{eq: update before retraction}
	\X^t + \eta^t = [\U^t \quad \U^t_\perp] \begin{bmatrix}
\B^{t+1}   & \D_2^{t+1\top} \\
\D^{t+1}_1  & \0
\end{bmatrix} [\V^t \quad \V^t_{\perp}]^\top.
\end{equation}
Comparing \eqref{eq: retraction step verification} and \eqref{eq: update before retraction}, we can view the update of $\X^{t+1}$ from $\X^t + \eta^t$ as simply completing the $\0$ matrix in $\begin{bmatrix}
\B^{t+1}   & \D_2^{t+1\top} \\
\D^{t+1}_1  & \0
\end{bmatrix}$ by $\D_1^{t+1} (\B^{t+1})^{-1} \D_2^{t+1\top}$. This operation maps the tangent vector on $T_{\X^t} \cM_r$ back to the manifold $\cM_r$ and it turns out that it coincides with the orthographic retraction
\begin{equation}\label{eq: orthographic retraction}
    R(\X^t, \eta^t) = [\U^t \quad \U^{t}_\perp]
\begin{bmatrix}
\B^{t+1} & \D_2^{t+1\top} \\
\D_1^{t+1} & \D_1^{t+1} {(\B^{t+1})}^{-1} \D_2^{t+1\top}
\end{bmatrix} [\V^t \quad \V^t_\perp]^{\top}
\end{equation}
on the set of fixed-rank matrices~\citep{absil2012projection}. Therefore, we have $\X^{t+1} = R(\X^t, \eta^t)$.
\begin{Remark}
    Although the orthographic retraction defined in~\cite{absil2012projection} requires that $\U^t$ and $\V^t$ are left and right singular vectors of $\X^t$, one can verify that even if the $\U^t$ and $\V^t$ are not exactly the left and right singular vectors but satisfy $\U^t = \widetilde \U^t \mathbf{O}$, $\V^t = \widetilde \V^t \mathbf{Q}$, then the mapping~\eqref{eq: orthographic retraction} is equivalent to the orthographic retraction in~\cite{absil2012projection}. Here, $\mathbf{O}, \mathbf{Q} \in \mathbb{O}_{r, r}$, and $\widetilde \U^t$ and $\widetilde \V^t$ are left and right singular vectors of $\X^t$.
\end{Remark}

The Riemannian gradient of a smooth function $f:\cM_r \to \bbR$ at $\X\in \cM_r$ is defined as the unique tangent vector ${\rm grad}\, f(\X) \in T_\X \cM_r$ such that $\inprod{{\rm grad}\, f(\X)}{\Z} = {\rm D}\, f(\X)[\Z], \forall\, \Z \in T_\X \cM_r,$
where ${\rm D} f(\X)[\Z]$ denotes the directional derivative of $f$ at point $\X$ along the direction $\Z$. Since $\cM_r$ is an embedded submanifold of $\bbR^{p_1\times p_2}$ and the Euclidean metric is used, from \cite[(3.37)]{absil2009optimization}, we know in our problem,
\begin{equation}\label{eq: Riemannian gradient}
{\rm grad}\, f(\X) = P_{T_\X}(\cA^*(\cA(\X) - \y)), 
\end{equation}
and here $P_{T_\X}$ is the orthogonal projector onto the tangent space at $\X$ defined as follows
\begin{equation} \label{eq: projection onto tangent space}
	P_{T_\X}(\Z) = P_\U \Z P_\V + P_{\U_\perp} \Z P_\V + P_\U \Z P_{\V_\perp},\quad \forall\, \Z\in\bbR^{p_1\times p_2},
\end{equation}
where $\U \in \bbO_{p_1, r}, \V \in \bbO_{p_2, r}$ are the left and right singular vectors of $\X $.

Next, we introduce the Riemannian Hessian. The Riemannian Hessian of $f$ at $\X\in\cM_r$ is the linear map ${\rm Hess} \,f(\X)$ of $T_\X\cM_r$ onto itself defined as ${\rm Hess}\, f(\X)[\Z] = \widebar{\nabla}_{\Z} {\rm grad}\, f, \quad \, \forall \Z \in T_\X \cM_r,$
where $\widebar{\nabla}$ is the Riemannian connection on $\cM_r$ \cite[Section 5.3]{absil2009optimization}. Lemma \ref{lm: Riemannian Hessian} gives an explicit formula for Riemannian Hessian in our problem.
\begin{Lemma}[Riemannian Hessian] \label{lm: Riemannian Hessian} 
	Consider $f(\X)$ in \eqref{eq:minimization}. If $\X \in \cM_r$ has singular value decomposition $\U \bSigma \V^\top$ and $\Z \in T_{\X} \cM_r$ has representation 
\begin{equation*}
   \Z = [\U \quad \U_\perp] \left[\begin{array}{c c}
       \Z_B & \Z_{D_2}^\top  \\
       \Z_{D_1}    & 0 
   \end{array} \right] [\V \quad \V_\perp]^\top, 
\end{equation*} 
then the Hessian operator in this setting satisfies
\begin{equation} \label{eq: Riemannian Hessian in least square}
	\begin{split}
	\Hess f(\X)[\Z] = & P_{T_\X} \left( \cA^* (\cA (\Z)) \right) + P_{\U_\perp} \cA^* (\cA (\X) - \y) \V_p \bSigma^{-1} \V^\top P_{\V} \\
	& + P_{\U} \U \bSigma^{-1} \U_p^\top \cA^* (\cA (\X) - \y) P_{\V_\perp},
	\end{split}
\end{equation}
where $\U_p = \U_\perp \Z_{D_1}, \V_p = \V_\perp \Z_{D_2}$.
\end{Lemma}

Next, we show that the update direction $\eta^t$, implicitly encoded in \eqref{eq: update before retraction}, finds the Riemannian Gauss-Newton direction in the manifold optimization of $\cM_r$.
Similar to the classic Newton's method, at $t$-th iteration, the Riemannian Newton method aims to find the Riemannian Newton direction $\eta^{t}_{\Newton}$ in $T_{\X^t} \cM_r$ that solves the following Newton equation
\begin{equation}\label{eq: Newton equation}
-\grad f(\X^t)= \Hess f(\X^t)[\eta_{\Newton}^t]. 
\end{equation}
If the residual $(\y - \cA(\X^t))$ is small, the last two terms in $\Hess f(\X^t)[\eta]$ of \eqref{eq: Riemannian Hessian in least square} are expected to be small, which means we can approximately solve the Riemannian Newton direction via
\begin{equation} \label{eq: Riemannian Gauss-newton equation}
	-\grad f(\X^t)= P_{T_{\X^t}} \left( \cA^* (\cA (\eta)) \right), \quad \eta \in T_{\X^t} \cM_r.
\end{equation}

In fact, Equation \eqref{eq: Riemannian Gauss-newton equation} has an interpretation from the \emph{Fisher scoring} algorithm. Consider the statistical setting $\y = \cA(\X)+\bepsilon$, where $\X$ is a fixed low-rank matrix and $\bepsilon_i \overset{i.i.d.}\sim N(0,\sigma^2)$. Then for any $\eta$,
\begin{equation*}
	\left\{\bbE (\Hess f(\X)[\eta] )\right\} |_{\X = \X^t} =P_{T_{\X^t}} \left( \cA^* (\cA (\eta)) \right),
\end{equation*}
where on the left hand side, the expression is evaluated at $\X^t$ after taking expectation. In the literature, the \emph{Fisher Scoring} algorithm computes the update direction via solving the modified Newton equation which replaces the Hessian with its expected value \citep{lange2010numerical}, 
i.e., $$ \left\{\bbE (\Hess f(\X)[\eta] )\right\} |_{\X = \X^t}  = - \grad f(\X^t), \quad \eta \in T_{\X^t} \cM_r,$$
which exactly becomes \eqref{eq: Riemannian Gauss-newton equation} in our setting. Meanwhile, it is not difficult to show that the Fisher Scoring algorithm here is equivalent to the Riemannian Gauss-Newton method for solving nonlinear least squares, see \cite[Section 14.6]{lange2010numerical} and \cite[Section 8.4]{absil2009optimization}. Thus, $\eta$ that solves the equation \eqref{eq: Riemannian Gauss-newton equation} is also the {\it Riemannian Gauss-Newton direction}.

It turns out that the update direction $\eta^t$ \eqref{eq: update before retraction} of RISRO solves the Fisher Scoring or Riemannian Gauss-Newton equation \eqref{eq: Riemannian Gauss-newton equation}: 
\begin{Theorem} \label{th: Riemannian Gauss-Newton of Simul Mini}
Let $\{\X^t\}$ be the sequence generated by RISRO under the same assumptions as in Theorem \ref{th: local contraction general setting}. Then, for all $t\ge 0$, the implicitly encoded update direction $\eta^t$ in \eqref{eq: update before retraction} solves the Riemannian Gauss-Newton equation \eqref{eq: Riemannian Gauss-newton equation}. 
\end{Theorem}

Theorem \ref{th: Riemannian Gauss-Newton of Simul Mini} together with the retraction explanation in \eqref{eq: orthographic retraction} establishes the connection of RISRO and Riemannian manifold optimization. Following this connection, we further show that each $\eta_t$ is always a decent direction in the next Proposition \ref{prop: descent direction}. This fact will be useful in boosting the local convergence of RISRO to the global convergence to be discussed in Remark \ref{rem: globalization}.
\begin{Proposition} \label{prop: descent direction}
For all $t\ge 0$, the update direction $\eta^t \in T_{\X^t} \cM_r $ in \eqref{eq: update before retraction} satisfies $\inprod{{\rm grad}f(\X^t)}{\eta^t} < 0,$
    i.e., $\eta_t$ is a descent direction.
    If  $\cA$ satisfies the 2r-RIP, then the direction sequence $\{\eta^t\}$ is gradient related.
\end{Proposition}
\begin{Remark}
The convergence of Riemannian Gauss-Newton was studied in a recent work \cite{breiding2018convergence}. Our results are significantly different from and offer improvements to \cite{breiding2018convergence} in the following ways. First, \cite{breiding2018convergence} considered a more general Riemannian Gauss-Newton setting, while their convergence results are established for a local minimum, which is a stronger and less practical requirement than the stationary point assumption we need. Second, the convergence rate in \cite{breiding2018convergence} includes several unspecified constants while we manage to work out all constants explicitly in our statement. Third, the local convergence radius in \cite{breiding2018convergence} does not specify the dependence on the $r$-th singular value of the target matrix while our result does. Fourth, our recursive importance sketching framework provides new sketching interpretations for several classical algorithms for rank constrained least squares. Finally, in Section \ref{sec: statistics applications} we also apply RISRO in popular statistical models and show RISRO achieves quadratic convergence in terms of {\rm estimation}. It is however not immediately clear how to utilize the results in \cite{breiding2018convergence} in these statistical settings.
\end{Remark}

\begin{Remark} \label{rem: value-of-importance-sketching-perspective} In addition to providing an interpretation of the superiority of RISRO, 
the Riemannian Gauss-Newton perspective developed in this section can inspire algorithmic developments in more general settings. For example, consider a general constrained optimization programming: $\min_{\X \in \cM} f(\X)$, where $\cM$ is an embedded submanifold of $R^N$ and $f$ is the restriction of a general twice differentiable objective in the ambient space to $\cM$. Although importance sketching is hard to define for this setting, Riemannian Gauss-Newton equation inspires to compute $\eta \in T_{\X^t} \cM$ by solving $P_{T_{\X^t}} \nabla^2 f(\X^t)[\eta] = - \grad f(\X^t)$,
then update the iterate as $\X^{t+1} = R({\X^t},\eta)$, where $R(\cdot, \cdot)$ is a retraction operator onto $\cM$. It is interesting to investigate the behavior of this algorithm from both optimization and statistical perspectives.

 Meanwhile, the recursive sketching perspective also provides solutions to a wider range of constrained optimization problems. For example, one can replace the $l_2$ loss, i.e., the least squares in Eq. \eqref{eq: alg2 least square} by other loss functions, such as the $l_1$ loss, Huber loss, or logistic loss, to handle different types of error corruptions and develop more robust algorithms. 
\end{Remark}

\begin{Remark}[Global Convergence of RISRO] \label{rem: globalization}
		By the classic theory of Riemannian optimization, the established connection of RISRO and Riemannian Gauss-Newton implies that vanilla RISRO may not converge when the RIP or the initialization condition fails. On the other hand, \cite[Section 8.4]{absil2009optimization} suggested that by adding or modifying the algorithm with certain line search or trust-region schemes, global convergence of Riemannian Gauss-Newton from any initialization to a stationary point can be guaranteed under proper assumptions. To be more specific, based on the Riemannian Gauss-Newton equation in \eqref{eq: Riemannian Gauss-newton equation} and Theorem \ref{th: Riemannian Gauss-Newton of Simul Mini}, the Riemannian Gauss-Newton direction at iteration $t$ satisfies
		\begin{equation} \label{eq: cha-RGN-direction}
			\begin{split}
				\eta^t = \arg\min_{\eta \in T_{\X^t} \cM_r } \|\cA P_{T_{\X^t}} (\X^t + \eta) - \y \|_2^2 = ( P_{T_{\X^t}} \cA^* \cA P_{T_{\X^t}} )^{-1} P_{T_{\X^t}} \cA^* (\y - \cA^*(\X^t) ).
			\end{split}
		\end{equation} After calculating $\eta^t$, we can update $\X^t$ to $\X^t + \eta^t$. 
		
		We can equip the algorithm with line search and update $\X^t$ to $\X^t + \alpha_t \eta^t$, where $\alpha_t$ is determined by some line search scheme, such as the Armijo method \cite[Section 4.3]{absil2009optimization}. Since the update direction $\eta^t$ is gradient related as shown in Proposition \ref{prop: descent direction} under RIP condition, this modified line search method has guaranteed global convergence property as shown in \cite[Theorem 4.3.1]{absil2009optimization}.
		
		We can also apply the trust region method to achieve global convergence. Specifically, we calculate the update direction as
		\begin{equation*}
			\begin{split}
				\tilde{\eta}^t = \arg\min_{\eta \in T_{\X^t} \cM_r, \eta \leq \Delta_t } \|\cA P_{T_{\X^t}} (\X^t + \eta) - \y \|_2^2
 			\end{split}
		\end{equation*} 
		for some radius $\Delta_t > 0$. Then if $\Delta_t$ is properly chosen such that $\tilde{\eta}^t$ guarantees sufficient decrease, the global convergence of this trust region method can be achieved under proper assumptions \cite[Theorem 7.4.2] {absil2009optimization}.
\end{Remark}

\section{Computational Complexity of RISRO} \label{sec: computation complexity}

In this section, we discuss the computational complexity of RISRO. Suppose $p_1 = p_2 = p$, the computational complexity of RISRO per iteration is $O(np^2r^2 + (pr)^3)$ in the general setting. A comparison of the computational complexity of RISRO and other common algorithms is provided in Table \ref{tab: time complexity compare}. Here the main complexity of RISRO and Alter Mini is from solving the least squares. The main complexity of the singular value projection (SVP) \citep{jain2010guaranteed} and gradient descent \citep{tu2016low} is from computing the gradient. From Table \ref{tab: time complexity compare}, we can see RISRO has the same per-iteration complexity as Alter Mini and comparable complexity with SVP and GD when $n \geq pr$ and $r$ is much less than $n$ and $p$. On the other hand, RISRO and Alter Mini are tuning-free, while a proper step size is crucial for SVP and GD to have fast convergence: the convergence theory of SVP and GD were often established when the step size is chosen to be smaller than a hard-to-find threshold; there are several practical ways to determine this step size and one needs to select the best one based on the data \citep{zheng2015convergent}, which may cost extra time. Finally, RISRO enjoys a high-order convergence as we have shown in Section~\ref{sec:theory}, and the convergence rates of all other algorithms are limited to being linear.
 \begin{table}
	\centering
	\begin{tabular}{c | c | c | c |c}
	\hline
	&  Alter Mini &  SVP & GD & RISRO (this work) \\
	\hline
	Complexity per iteration & $O(np^2 r^2 + (pr)^3 )$  & $O(np^2)$ & $O(np^2)$ & $O(np^2 r^2 + (pr)^3)$\\
	\hline
	Convergence rate & Linear & Linear & Linear & Quadratic-(linear)\\
	\hline
	\end{tabular}
	\caption{Computational complexity per iteration and convergence rate for Alternating Minimization (Alter Mini) \citep{jain2013low}, singular value projection (SVP) \citep{jain2010guaranteed}, gradient descent (GD) \citep{tu2016low}, and RISRO} \label{tab: time complexity compare}
\end{table}

The main computational bottleneck of RISRO is solving the least squares, which can be alleviated by using iterative linear system solvers, such as the (preconditioned) conjugate gradient method when the linear operator $\mathcal{A}$ has special structures. Such special structures occur, for example, in matrix completion problem ($\mathcal{A}$ is sparse)~\citep{vandereycken2013low}, phase retrieval for X-ray crystallography
imaging ($\mathcal{A}$ involves fast Fourier transforms)~\citep{huang2017solving}, and blind deconvolution for imaging deblurring ($\mathcal{A}$ involves fast Fourier transforms and Haar wavelet transforms)~\citep{huang2018blind}.

To utilize these structures, we introduce an intrinsic representation of tangent vectors in $\mathcal{M}_r$: if $\U, \V$ are the left and right singular vectors of a rank-$r$ matrix $\X$, an orthonormal basis of $T_\X \mathcal{M}_r$ can be
\begin{align*}
&\left\{ 
	[\U\quad \U_{\perp}] \begin{bmatrix}
	\e_i \e_j^\top & \0_{r \times (p - r)} \\[2pt]
	\0_{(p - r) \times r} & \0_{(p -r)\times (p - r)}
	\end{bmatrix}
	[\V\quad \V_{\perp}]^\top,
i = 1, \ldots, r, j = 1, \ldots, r
 \right\}
\cup \\
&\left\{ 
	[\U\quad \U_{\perp}] \begin{bmatrix}
	\0_{r \times r} & \e_i \tilde{\e}_j^\top \\[2pt]
	\0_{(p - r) \times r} & \0_{(p -r)\times (p - r)}
	\end{bmatrix}
	[\V\quad \V_{\perp}]^\top,
i = 1, \ldots, r, j = 1, \ldots, p - r
 \right\}
\cup \\
&\left\{ 
	[\U\quad \U_{\perp}] \begin{bmatrix}
	\0_{r \times r} & \0_{r \times (p - r)} \\[2pt]
	\tilde{\e}_i \e_j^\top & \0_{(p -r)\times (p - r)}
	\end{bmatrix}
	[\V\quad \V_{\perp}]^\top,
i = 1, \ldots, p - r, j = 1, \ldots, r
 \right\},
\end{align*}
where $\e_i$ and $\tilde{\e}_i$ denote the $i$-th canonical basis of $\mathbb{R}^r$ and $\mathbb{R}^{p - r}$, respectively. It follows that any tangent vector in $T_{\X} \mathcal{M}_r$ can be uniquely represented by a coefficient vector in $\mathbb{R}^{(2 p - r) r}$ via the basis above. This representation is called the intrinsic representation \citep{huang2017intrinsic}. Computing the intrinsic representations of a Riemannian gradient can be computationally efficient. For example, the complexity of computing the Riemannian gradient in matrix completion is $O(n r + p r^2)$ and its intrinsic representation can be computed by an additional $O(p r^2)$ operations \citep{vandereycken2013low}. The complexities of computing intrinsic representations of the Riemannian gradients of the phase retrieval and the blind deconvolution are both $O(n \log(n) r + p r^2)$ \citep{huang2017solving,huang2018blind}. 

By Theorem~\ref{th: Riemannian Gauss-Newton of Simul Mini}, the least squares problem~\eqref{eq: alg2 least square} of RISRO is equivalent to solve $\eta \in T_{\X^{t}}\cM_r$ such that $P_{T_{\X^t}} \cA^* (\cA (\eta)) = -\grad f(\X^t)$. Reformulating this equation by intrinsic representation yields
\begin{equation} \label{eq:intrModNewtonsystem}
	- \grad f(\X^t) = P_{T_{\X^t}} \cA^* ( \cA(\eta)) \quad
\Longrightarrow -u = \mathcal{B}_{\mathbf{X}}^{*} (\mathcal{A}^* (\mathcal{A}( \mathcal{B}_{\mathbf{X}} v ))), 
\end{equation}
where $u, v$ are the intrinsic representations of $\grad f(\X^t)$ and $\eta$, the mapping $\mathcal{B}_{\mathbf{X}}: \mathbb{R}^{(2 p - r)r} \rightarrow T_\X \mathcal{M}_r \subset \mathbb{R}^{p \times p}$ converts an intrinsic representation to the corresponding tangent vector, and $\mathcal{B}_{\mathbf{X}}^*: \mathbb{R}^{p \times p} \rightarrow \mathbb{R}^{(2 p - r)r}$ is the adjoint operator of $\mathcal{B}_{\mathbf{X}}$. The computational complexity of using conjugate gradient method to solve \eqref{eq:intrModNewtonsystem} is determined by the complexity of evaluating the operator $\mathcal{B}_{\mathbf{X}}^{*} \circ (\mathcal{A}^* \mathcal{A}) \circ \mathcal{B}_{\mathbf{X}}$ on a given vector. With the intrinsic representation, it can be shown that this evaluation costs $O(n r + p r^2)$ in matrix completion and $O(n \log(n) r + p r^2)$ in the phase retrieval and blind deconvolution. Thus, when solving \eqref{eq:intrModNewtonsystem} via the conjugate gradient method, the complexity is $O(k (n r + p r^2))$ in the matrix completion and $O(k (n \log(n) r + p r^2))$ in the phase retrieval and the blind deconvolution, where $k$ is the number of conjugate gradient iterations and is provably at most $(2 p - r)r$.  Hence, for special applications such as matrix completion, phase retrieval and blind deconvolution, by using the conjugate gradient method with the intrinsic representation, the per iteration complexity of RISRO can be greatly reduced. This point will be further exploited in our future research.

\section{Recursive Importance Sketching under Statistical Models}\label{sec: statistics applications}

In this section, we study the applications of RISRO in machine learning and statistics. We specifically investigate the low-rank matrix trace regression and phase retrieval, while our key ideas can be applied to more problems. For the execution of RISRO, we assume that some estimate for the rank of the target parameter matrix, denoted by $r$, is available. 
In many statistical applications such as phase retrieval and blind deconvolution, this assumption trivially holds as the parameter matrix is known to be rank-1. In other applications, while the rank of the parameter is unknown, it is generally not difficult to obtain a rough estimate given the domain knowledge. Then, we can optimize over the set of fixed rank matrices using the formulation of~\eqref{eq:minimization} and dynamically update the selected rank (see, e.g.,~\cite{vandereycken2010riemannian,zhou2016riemannian}). 

\subsection{Low-Rank Matrix Trace Regression} \label{sec: low rank matrix regression}

Consider the low-rank matrix trace regression model:
\begin{equation} \label{eq: low rank matrix recovery model}
	\y_i = \langle \A_i, \X^* \rangle + \bepsilon_i, \quad \text{ for } 1 \leq i \leq n,
\end{equation} 
where $\X^* \in \bbR^{p_1 \times p_2}$ is the true model parameter to be estimated. We estimate $\X^*$ by solving~\eqref{eq:minimization} where $r$ in the rank constraint satisfies $r \leq {\rm rank}(\X^*)$, i.e., $r$ is an estimate of ${\rm rank}(\X^*)$. 

The following Theorem \ref{th: local convergence in local rank matrix recovery} shows RISRO converges quadratically to the best rank $r$ approximation of $\X^*$, i.e., $\X^*_{\max(r)}$, up to some statistical error given a proper initialization. Under the Gaussian ensemble design, RISRO with spectral initialization achieves the minimax optimal estimation error rate. 

\begin{Theorem}\label{th: local convergence in local rank matrix recovery}
    {\bf (RISRO in Matrix Trace Regression)}
	Consider the low-rank matrix trace regression problem \eqref{eq: low rank matrix recovery model}. Define $\tilde{\bepsilon}_i := \bepsilon_i + \langle \A_i, \X^* - \X^*_{\max(r)} \rangle$ for $i = 1,\ldots, n$. Suppose that $\cA$ satisfies the $2r$-RIP, the initialization of RISRO satisfies
\begin{equation} \label{ineq: low rank matrix recovery initialization condition}
	\|\X^0 - \X^*_{\max(r)}\|_\F \leq \left( \frac{1}{4} \wedge \frac{1-R_{2r}}{2\sqrt{5}R_{3r}}  \right) \sigma_r(\X^*), 
\end{equation} 
and
\begin{equation} \label{ineq: LRMR sigma condition}
	\sigma_r(\X^*) \geq \left(16\sqrt{5} \vee  \frac{40 \sqrt{2} R_{3r}}{1-R_{2r}} \right) \frac{\|(\cA^*(\tilde{\bepsilon}))_{\max(r)}\|_\F}{1-R_{2r}}.
\end{equation} 
Then the iterations of RISRO converge as follows $\forall \, t\ge 0$:
	\begin{equation} \label{ineq: low rank matrix recovry quadratic-linear converge}
	\begin{split}
		& \|\X^{t+1} - \X^*_{\max(r)}\|_\F^2  \leq 10\frac{R_{3r}^2 \|\X^t - \X^*_{\max(r)} \|_\F^4 }{(1-R_{2r})^2\sigma_r^2(\X^*)} + \frac{20 \| ( \cA^*(\tilde{\bepsilon} ) )_{\max(r)} \|_\F^2}{(1-R_{2r})^2}\\
		\leq & 10\frac{R_{3r}^2 \|\X^t - \X^* \|_\F^4 }{(1-R_{2r})^2\sigma_r^2(\X^*)} +  \frac{20 \left(\| ( \cA^*(\bepsilon ) )_{\max(r)} \|_\F + \| (\cA^*\cA( \X^* - \X^*_{\max(r)} ))_{\max(r)} \|_\F \right)^2}{(1-R_{2r})^2}. 
	\end{split}
	\end{equation}
 The overall convergence of RISRO shows two phases:
\begin{itemize}
	\item (Phase I) When $\|\X^{t} - \X^*_{\max(r)}\|_\F^2 \geq \frac{\sqrt{2}}{R_{3r}} \| ( \cA^*(\tilde{\bepsilon} ) )_{\max(r)} \|_\F  \sigma_r(\X^*)$,
	\begin{equation*}
		\|\X^{t+1} - \X^*_{\max(r)}\|_\F \leq 2\sqrt{5}\frac{R_{3r} \|\X^t - \X^*_{\max(r)} \|_\F^2 }{(1-R_{2r})\sigma_r(\X^*)}, \quad \|\X^{t+1} - \X^*\|_\F \leq 2\sqrt{5}\frac{R_{3r} \|\X^t - \X^*_{\max(r)} \|_\F^2 }{(1-R_{2r})\sigma_r(\X^*)} + \|\X^*_{-\max(r)}\|_\F,
	\end{equation*} where $\X^*_{-\max(r)} = \X^* - \X^*_{\max(r)}$.
	\item (Phase II) When $\|\X^{t} - \X^*_{\max(r)}\|_\F^2 \leq \frac{\sqrt{2}}{R_{3r}} \| ( \cA^*(\tilde{\bepsilon} ) )_{\max(r)} \|_\F  \sigma_r(\X^*)$,
	\begin{equation*}
		\|\X^{t+1} - \X^*_{\max(r)}\|_\F \leq  \frac{2\sqrt{10} \| ( \cA^*(\tilde{\bepsilon} ) )_{\max(r)} \|_\F}{1-R_{2r}}, \quad 		\|\X^{t+1} - \X^*\|_\F \leq  \frac{2\sqrt{10} \| ( \cA^*(\tilde{\bepsilon} ) )_{\max(r)} \|_\F}{1-R_{2r}} + \|\X^*_{-\max(r)}\|_F.
	\end{equation*}
\end{itemize}

Moreover, we assume $\rank(\X^*) = r$, $(\A_i)_{[j,k]}$ are independent sub-Gaussian random variables with mean zero and variance $1/n$ and $\bepsilon_i$ are independent sub-Gaussian random variables with mean zero and variance $\sigma^2/n$ (i.e., $\mathbb{E}(\A_i)_{[j,k]}=\mathbb{E}(\bepsilon_i)=0$, $\textrm{Var}((\A_i)_{[j,k]}) = 1/n$, $\textrm{Var}(\bepsilon_i) = \sigma^2/n$, $\sup_{q\geq 1} (n/q)^{1/2}(\mathbb{E}|(\A_i)_{[j,k]}|^q)^{1/q}\leq C$, $\sup_{q\geq 1} (n/(q\sigma^2))^{1/2}(\mathbb{E}|\bepsilon_i|^q)^{1/q}\leq C$ for some fixed $C > 0$). Then there exist universal constants $C_1, C_2, C', c>0$ such that as long as $n \geq C_1 (p_1 + p_2)r(\frac{\sigma^2}{\sigma^2_r(\X^*)} \vee r\kappa^2)$ (here $\kappa = \frac{\sigma_{1}(\X^*)}{\sigma_r(\X^*)}$ is the condition number of $\X^*$) and  $t_{\max} \geq C_2 \log \log (\frac{\sigma^2_r(\X^*) n}{r(p_1 + p_2)\sigma^2}) \vee 1$, the output of RISRO with spectral initialization $\X^0 = (\cA^*(\y))_{\max(r)}$ satisfies $\|\X^{t_{\max}} - \X^*\|_\F^2 \leq c\frac{r(p_1+p_2)}{n}\sigma^2$
	with probability at least $1 - \exp(-C'(p_1+ p_2))$.
\end{Theorem}

\begin{Remark}\label{rem: trace regression remark 1}
	{\bf (Quadratic Convergence, Two-phases Convergence, Statistical Error, and Robustness)}
	The upper bound of $\|\X^t - \X^*_{\max(r)}\|^2_\F$ in \eqref{ineq: low rank matrix recovry quadratic-linear converge} includes two terms: the optimization error term $O(\|\X^t-\X^*_{\max(r)}\|_F^4)$ quadratically decreases over iteration $t$, and the statistical error term $O(\|(\cA^*(\tilde{\bepsilon}))_{\max(r)} \|_\F^2)$ is static through iterations. Moreover, RISRO includes two phases in its convergence. In Phase I with large $\|\X^{t} - \X^*_{\max(r)}\|_\F^2 $, RISRO converges quadratically towards $\X^*_{\max(r)}$; in Phase II with moderate $\|\X^{t} - \X^*_{\max(r)}\|_\F^2$, the estimator returned by one more iteration of RISRO achieves the best possible statistical error rate $O(\| (\cA^*(\tilde{\bepsilon}) )_{\max(r)} \|_\F^2)$ as suggested by the $d=2$ case in \cite[Theorem 2]{luo2021low}. Therefore, although the convergence rate of RISRO may decelerate to be linear in Phase II, Theorem \ref{th: local convergence in local rank matrix recovery} suggests there is no need to run further iterations as the estimator is already statistically optimal after one additional iteration. Such performance of ``quadratic convergence + one-iteration optimality'' is unique, which does not appear in common first-order methods.
	
	Finally, \eqref{ineq: low rank matrix recovry quadratic-linear converge} shows the error contraction factor is independent of the condition number $\kappa$, which demonstrates the robustness of RISRO to the ill-conditioning of the underlying low-rank matrix. We will further demonstrate this point by simulation studies in Section \ref{sec: numerical comparison}.
\end{Remark}

\begin{Remark}[Optimal Statistical Error] \label{rem: regression setting, statistical error}
	Under the Gaussian ensemble design and when $\rank(\X^*) = r$, RISRO with spectral initialization achieves the rate of estimation error $cr(p_1 + p_2)\sigma^2/n$ after {\it double-logarithmic} number of iterations when $n\geq C_1(p_1 + p_2 )r( \frac{\sigma^2}{\sigma^2_r(\X^*)} \vee r\kappa^2)$. Compared with the lower bound of the estimation error
	    $$\min_{\widehat{\X}} \max_{\rank(\X^*) \leq r} \bbE \|\widehat{\X} - \X^*\|_\F^2 \geq c' \frac{r(p_1 + p_2)\sigma^2}{n} $$ for some $c' >0$
	in \cite{candes2011tight}, RISRO achieves the minimax optimal estimation error with near-optimal sample complexity. To the best of our knowledge, RISRO is the first provable algorithm that achieves the minimax rate-optimal estimation error with only a double-logarithmic number of iterations and this is an exponential improvement over common first-order methods where a logarithmic number of iterations are needed.
\end{Remark}

\subsection{Phase Retrieval} \label{sec: phase retrieval}

In this section, we consider RISRO for solving the following quadratic equation system
\begin{equation} \label{eq: phase retrieval model}
	\y_i = |\langle\a_i, \x^*\rangle|^2 \quad \text{for} \quad 1 \leq i \leq n,
\end{equation} 
where $\y \in \bbR^n$ and covariates $\{\a_i \}_{i=1}^n  \in \bbR^{p}$ (or $\mathbb{C}^p$) are known whereas $\x^* \in \bbR^p$ (or $\mathbb{C}^p$) are unknown. The goal is to recover $\x^*$ based on $\{\y_i, \a_i\}_{i=1}^n$. One important application is known as \emph{phase retrieval} arising from physical science due to the nature of optical sensors \citep{fienup1982phase}. In the literature, various approaches have been proposed for phase retrieval with provable guarantees, such as convex relaxation \citep{candes2013phaselift,huang2017solving,waldspurger2015phase} and non-convex approaches \citep{candes2015phase,chen2017solving,gao2017phaseless,ma2019implicit,netrapalli2013phase,sanghavi2017local,wang2017solving,duchi2019solving}.

For ease of exposition, we focus on the real-value model, i.e., $\x^* \in \bbR^n$ and $\a_i \in \bbR^n$, while a simple trick in \cite{sanghavi2017local} can recast the problem \eqref{eq: phase retrieval model} in the complex model into a rank-2 real value matrix recovery problem, then our approach still applies. In the real-valued setting, we can rewrite model \eqref{eq: phase retrieval model} into a low-rank matrix recovery model 
\begin{equation}
\label{eq: matrix model pha}
\y = \cA(\X^*) \mbox{ with } \X^* = \x^* \x^{*\top} \mbox{ and } [\cA(\X^*)]_i = \langle \a_i \a_i^\top, \x^* \x^{*\top} \rangle.
\end{equation} There are two challenges in phase retrieval compared to the low-rank matrix trace regression considered previously. First, due to the symmetry of sensing matrices $\a_i \a_i^\top$ and $\x^* \x^{*\top}$ in phase retrieval, the importance covariates $\cA_{D_1}$ and $\cA_{D_2}$ in \eqref{eq: importance sketches} are exactly the same and an adaptation of Algorithm \ref{alg: recursive IS alg 2} is thus needed. Second, in phase retrieval, the mapping $\cA$ no longer satisfies a proper RIP condition in general \citep{cai2015rop,candes2013phaselift}, so a new theory is needed. 
To this end, we introduce a modified RISRO for phase retrieval in Algorithm \ref{alg: recursive IS alg for phase retrieval}. Particularly in Step 4 of Algorithm \ref{alg: recursive IS alg for phase retrieval}, we multiply the importance covariates $\A_2$ by an extra factor $2$ to account for the duplicate importance covariates due to symmetry. 

\begin{algorithm}[h]
\caption{RISRO for Phase Retrieval}
	\begin{algorithmic}[1]
		\STATE Input: design vectors $\{\a_i \}_{i=1}^n \in \bbR^p$, $\y\in \mathbb{R}^n$, initialization $\X^0$ that admits eigenvalue decomposition $\sigma_1^0 \u^0 \u^{0\top}$ 
		\FOR{$t=0, 1,\ldots, $} 
		\STATE Perform importance sketching on $\a_i$ and construct the covariates $\A_1 \in \bbR^n$, $\A_2 \in \bbR^{n \times (p-1)}$, where for $1 \leq i \leq n$, $(\A_1)_i = (\a_i^\top \u^t)^2, \quad (\A_{2})_{[i,:]} = \u_{\perp}^{t\top} \a_i \a_i^\top \u^{t}.$
		\STATE \label{alg: PR least square} Solve the unconstrained least squares problem $( b^{t+1}, \d^{t+1}) = \argmin_{b\in \bbR, \d \in \bbR^{(p -1)}} \left\|\y - \A_1 b - 2\A_{2} \d \right\|_2^2.$
		\STATE Compute the eigenvalue decomposition of 
		$ [\u^t \,\u^t_\perp ] \left[  \begin{array}{c c}
				b^{t+1} & \d^{t+1 \top}\\
				\d^{t+1} & \0
			\end{array} \right] [\u^t\, \u^t_\perp]^\top$, and denote it as
		$
				 [\v_1 \, \v_2] \left[\begin{array}{c c}
				\lambda_1 & 0\\
				0 & \lambda_2
			\end{array}\right] [\v_1 \, \v_2]^\top$
		 with $\lambda_1 \geq \lambda_2$. 
		\STATE Update $\u^{t+1} = \v_1$ and $\X^{t+1} = \lambda_1 \u^{t+1} \u^{t+1 \top}$.
		\ENDFOR
	\end{algorithmic} \label{alg: recursive IS alg for phase retrieval}
\end{algorithm}

Next, we show under Gaussian ensemble design, given the sample number $n = O(p \log p)$ and proper initialization, the sequence $\{\X^t \}$ generated by Algorithm \ref{alg: recursive IS alg for phase retrieval} converges quadratically to $\X^*$.

\begin{Theorem}[Local Quadratic Convergence of RISRO for Phase Retrieval] \label{thm: quadratic convergence for phase retrieval}
	In the phase retrieval problem \eqref{eq: phase retrieval model}, assume that 
	$\{\a_i\}_{i=1}^n$ are independently generated from  $N(0,\I_p)$. Then for any $\delta_1, \delta_2 \in (0,1)$, there exist $c,C(\delta_1), C' > 0$ such that when $p \geq c \log n, n \geq C(\delta_1) p \log p$, if $\|\X^0 - \X^*\|_\F \leq \frac{(1-\delta_1)}{C' (1+\delta_2) p} \|\X^*\|_\F $, with probability at least $1 - C_1\exp(-C_2(\delta_1, \delta_2)n) - C_3 n^{-p}$, the sequence $\{\X^t\}$ generated by Algorithm \ref{alg: recursive IS alg for phase retrieval} satisfies
	\begin{equation} \label{ineq: quadra conver in PR}
		\|\X^{t+1} - \X^*\|_\F \leq \frac{C'(1+\delta_2) p}{(1-\delta_1)\|\X^*\|_\F} \|\X^t - \X^*\|_\F^2, \quad \forall \, t\ge 0
	\end{equation}
	for some $C_1,C_2(\delta_1, \delta_2), C_3 >0$.
\end{Theorem}

	To overcome the technical difficulties in establishing quadratic convergence without RIP for phase retrieval, Theorem \ref{thm: quadratic convergence for phase retrieval} is established under the assumption $\|\X^0 - \X^*\|_\F \leq O(\|\X^*\|_\F/p)$. 
	Although it is difficult to prove that the spectral initializer meets this assumption under the near-optimal sample size (e.g., $n=Cp\log p$), 
	we find by simulation that the spectral initialization yields quadratic convergence for RISRO (Section \ref{sec:numerics}). On the other hand, we can also run a few iterations of factorized gradient descent to achieve the initialization condition in Theorem \ref{thm: quadratic convergence for phase retrieval} with near-optimal sample complexity guarantee \citep{candes2015phase,chen2017solving,ma2019implicit} and then switch to RISRO. Specifically, the initialization algorithm for RISRO in phase retrieval via factorized gradient descent is provided in Algorithm \ref{alg: initialization phase retrieval} and its guarantee is given in Proposition \ref{prop: init-phase-retrieval}.
\begin{algorithm}[h]
\caption{
RISRO for Phase Retrieval with Gradient Descent Initialization}
	\begin{algorithmic}[1]
		\STATE Input: design vectors $\{\a_i \}_{i=1}^n \in \bbR^p$ and $\y\in \mathbb{R}^n$. 
		\STATE Let $\lambda_1(\Y)$ and $\v_1$ be the leading eigenvalue and eigenvector of $\Y = \frac{1}{n} \sum_{j=1}^n \y_j \a_j \a_j^\top $, respectively, and set $\widetilde{\x}^0 = \sqrt{\lambda_1(\Y)/3} \v_1$. 
		\FOR{$t=0, 1,\ldots, T_0-1$} 
		\STATE Update $\widetilde{\x}^{t+1} = \widetilde{\x}^t - \eta_t \nabla g(\widetilde{\x}^t) $, where $g(\x) = \frac{1}{4n} \sum_{j=1}^n \left(  (\a_j^\top \x)^2 - \y_j  \right)^2$.
		\ENDFOR
		\STATE Apply Algorithm \ref{alg: recursive IS alg for phase retrieval} with initialization $\widetilde{\x}^{T_0} \widetilde{\x}^{T_0\top}$.
	\end{algorithmic}\label{alg: initialization phase retrieval}
\end{algorithm}

\begin{Proposition} \label{prop: init-phase-retrieval}
	In phase retrieval \eqref{eq: phase retrieval model}, suppose $\{\a_i\}_{i=1}^n$ are independently drawn from $N(0,\I_p)$ and $n \geq C p \log p$ for some sufficient large constant $C > 0$. Assume the step size in Algorithm \ref{alg: initialization phase retrieval} obeys $\eta_t \equiv \eta = c_1/(\log p \cdot \|\widetilde{\x}^0 \|_2^2 )$ for constant $c_1 > 0$, where $\widetilde{\x}^0$ is given in the algorithm. Then there exist absolute constants $c_2, c_3 > 0$ such that when $T_0 \geq c_2 \log p \cdot \log ( \|\x^*\|_2 p)$, the initialization $\X^0:= \widetilde{\x}^{T_0} \widetilde{\x}^{T_0 \top} $ in Algorithm \ref{alg: initialization phase retrieval} satisfies the initialization condition in Theorem \ref{thm: quadratic convergence for phase retrieval} and the conclusion of Theorem \ref{thm: quadratic convergence for phase retrieval} holds with probability at least $1- c_3 n p^{-5}$.
\end{Proposition}

\section{Numerical Studies}\label{sec:numerics}

In this section, we conduct simulation studies to investigate the numerical performance of RISRO. We specifically consider two settings: \begin{itemize}
\item Matrix trace regression. Let $p = p_1 = p_2$ and $\y_i = \langle \X^*, \A_i \rangle + \bepsilon_i $, where $\A_i$s are constructed with independent standard normal entries and $\bepsilon_i\overset{i.i.d.}\sim N(0, \sigma^2)$. $\X^* = \U^* \bSigma^* \V^{*\top}$ where $\U^*, \V^* \in \mathbb{O}_{p,r}$ are randomly generated, $\bSigma = \diag(\lambda_1,\ldots, \lambda_r)$. Also, we set $\lambda_1 = 3$ and $\lambda_i = \frac{\lambda_1}{\kappa^{i/r}}$ for $i= 2, \ldots, r$, so the condition number of $\X^*$ is $\kappa$. We initialize $\X^0$ via $(\cA^*(\y))_{\max(r)}$.
\item Phase retrieval. Let $\y_i = \langle \a_i , \x^* \rangle^2$, where $\x^* \in \bbR^p$ is a randomly generated unit vector, $\a_i \overset{i.i.d.}\sim N(0, \I_p)$. We initialize $\X^0$ via truncated spectral initialization \citep{chen2017solving}.
\end{itemize}
Throughout the simulation studies, we consider errors in two metrics: (1) $\|\X^{t} - \X^{t_{\max}}\|_{\F}/ \|\X^{t_{\max}}\|_{\F}$, which measures the convergence error; (2) $\|\X^{t} - \X^*\|_{\F}/ \|\X^*\|_{\F}$, which is the relative root mean-squared error (Relative RMSE) that measures the estimation error for $\X^*$. The algorithm is terminated when it reaches the maximum number of iterations $t_{\max} = 300$ or the corresponding error metric is less than $10^{-12}$. Unless otherwise noted, the reported results are based on the averages of 50 simulations and on a computer with Intel Xeon E5-2680 2.5GHz CPU. Additional development and simulation results for RISRO in matrix completion and robust PCA can be found in Appendix \ref{sec: RISRO-matrix-completion-robust-PCA}.

\subsection{Properties of RISRO} \label{sec: numerical Simul Mini property}

We first study the convergence rate of RISRO. Specifically, set $p = 100, r = 3, n \in \{1200, 1500, 1800, 2100, 2400\}, \kappa =1, \sigma = 0$ for low-rank matrix trace regression and $p = 1200, n \in \{4800, 6000, 7200, 8400, 9600\}$ for phase retrieval. The convergence performance of RISRO (Algorithm \ref{alg: recursive IS alg 2} in low-rank matrix trace regression and Algorithm \ref{alg: recursive IS alg for phase retrieval} in phase retrieval) is plotted in Figure \ref{fig: Simul Mini performance illustration}. We can see RISRO with the (truncated) spectral initialization converges quadratically to the true parameter $\X^*$ in both problems, which is in line with the theory developed in previous sections. Although our theory on phase retrieval in Theorem \ref{thm: quadratic convergence for phase retrieval} is based on a stronger initialization assumption, the truncated spectral initialization achieves great empirical performance. 

\begin{figure}[ht]
	\centering
	\subfigure{\includegraphics[width=0.7\textwidth]{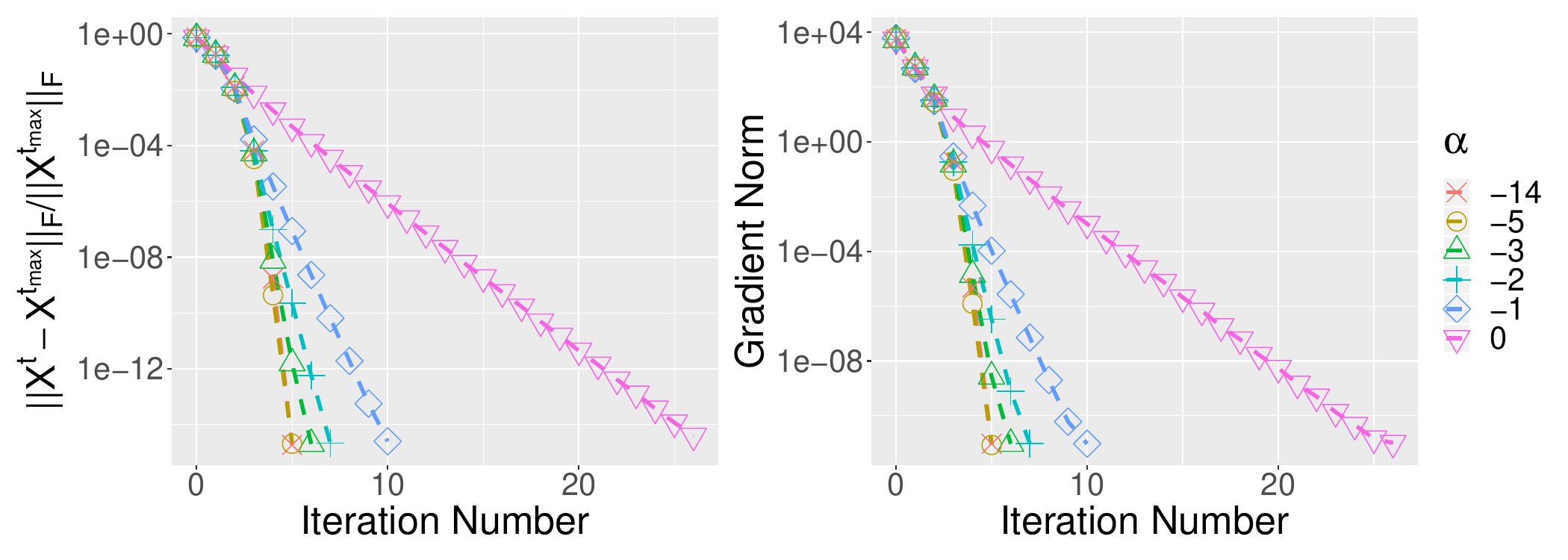}}
	\caption{Convergence plot of RISRO in matrix trace regression. $p = 100, r = 3, n = 1500, \kappa = 1$, $\sigma = 10^{\alpha}$ with varying $\alpha$} \label{fig: Simul Mini convergence prop}
\end{figure}

In another setting, we examine the quadratic-linear convergence for RISRO under the noisy setting. Consider the matrix trace regression problem, where $\sigma = 10^{\alpha}$, $\alpha \in \{0,-1,-2,-3,-5,-14 \}$, $n = 1500$, and $p, r, \kappa$ are the same as the previous setting. The simulation results in Figure \ref{fig: Simul Mini convergence prop} show the gradient norm $\|\grad f(\X^t)\|$ of the iterates converges to zero, which demonstrates the convergence of the algorithm. Meanwhile, since the observations are noisy, RISRO exhibits the quadratic-linear convergence as we discussed in Remark \ref{rem: quadratic-linear convergence}: when $\alpha = 0$, i.e., $\sigma = 1$, RISRO converges quadratically in the first 2-3 steps and then reduces to linear convergence afterward; as $\sigma$ gets smaller, we can see RISRO enjoys a longer path of quadratic convergence, which matches our theoretical prediction in Remark \ref{rem: quadratic-linear convergence}. 

\begin{figure}[ht]
	\centering
	\subfigure{\includegraphics[width=0.7\textwidth]{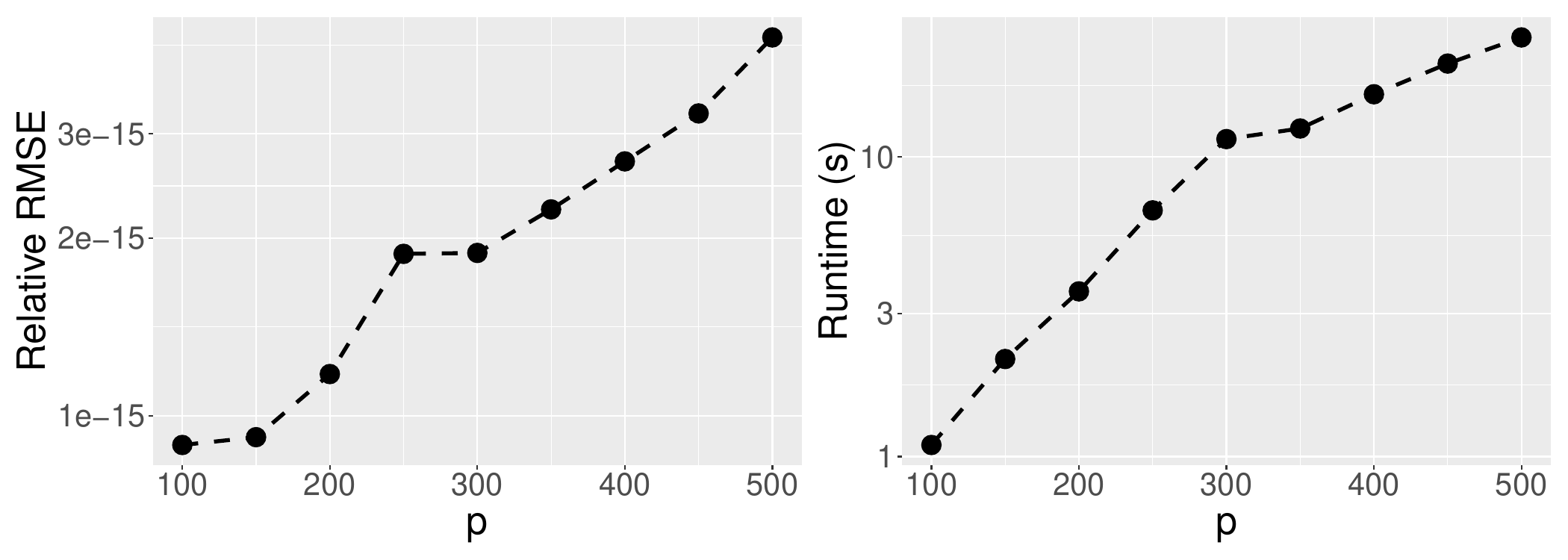}}
	\caption{Relative RMSE and runtime of RISRO in matrix trace regression. $p \in [100,500], r = 3, n = 7000, \kappa = 1$, $\sigma = 0$} \label{fig: Simul Mini convergence prop grow p}
\end{figure}

Finally, we study the performance of RISRO under the large-scale setting of the matrix trace regression. Fix $n = 7000, r = 3, \kappa = 1, \sigma = 0$ and let dimension $p$ grow from $100$ to $500$. For the largest case, the space cost of storing $\cA$ reaches $7000\cdot 500\cdot 500\cdot 8$B $= 13.04$GB. Figure \ref{fig: Simul Mini convergence prop grow p} shows the relative RMSE of the output of RISRO and runtime versus the dimension. We can clearly see the relative RMSE of the output is stable and the runtime scales reasonably well as the dimension $p$ grows.

\subsection{Comparison of RISRO with Other Algorithms in Literature } \label{sec: numerical comparison}

In this subsection, we further compare RISRO with existing algorithms in the literature. In the matrix trace regression, we compare our algorithm with singular value projection (SVP) \citep{goldfarb2011convergence,jain2010guaranteed}, Alternating Minimization (Alter Mini) \citep{jain2013low,zhao2015nonconvex}, gradient descent (GD) \citep{park2018finding,tu2016low,zheng2015convergent}, and convex nuclear norm minimization (NNM) \eqref{eq: convex relaxation} \citep{toh2010accelerated}. We consider the setting with $p = 100, r = 3, n = 1500$, $\kappa \in \{ 1, 50, 500\}$, $\sigma = 0$ (noiseless case) or $\sigma = 10^{-6}$ (noisy case). Following \cite{zheng2015convergent}, in the implementation of GD and SVP, we evaluate three choices of step size, $\{5\times 10^{-3}, 10^{-3}, 5 \times 10^{-4} \}$, then choose the best one. In phase retrieval, we compare Algorithm \ref{alg: recursive IS alg for phase retrieval} with Wirtinger Flow (WF) \citep{candes2015phase} and Truncated Wirtinger Flow (TWF) \citep{chen2017solving} with $p = 1200, n = 6000$. We use the codes of the accelerated proximal gradient for NNM, WF and TWF from the corresponding authors' websites and implement the other algorithms by ourselves. The stopping criteria of all procedures are the same as RISRO mentioned in the previous simulation settings.

\begin{figure}[ht!]
	\centering
	\subfigure[$\kappa = 1$]{\includegraphics[width=0.7\textwidth]{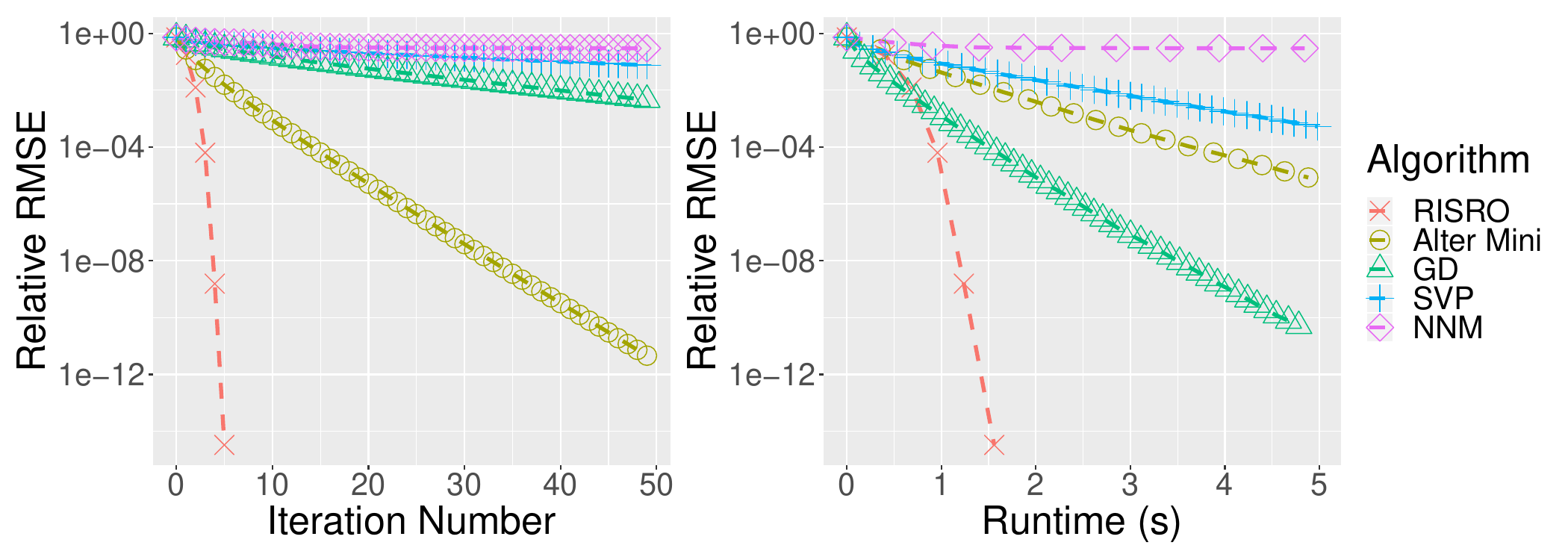}}
	\subfigure[$\kappa = 50$]{\includegraphics[width=0.7\textwidth]{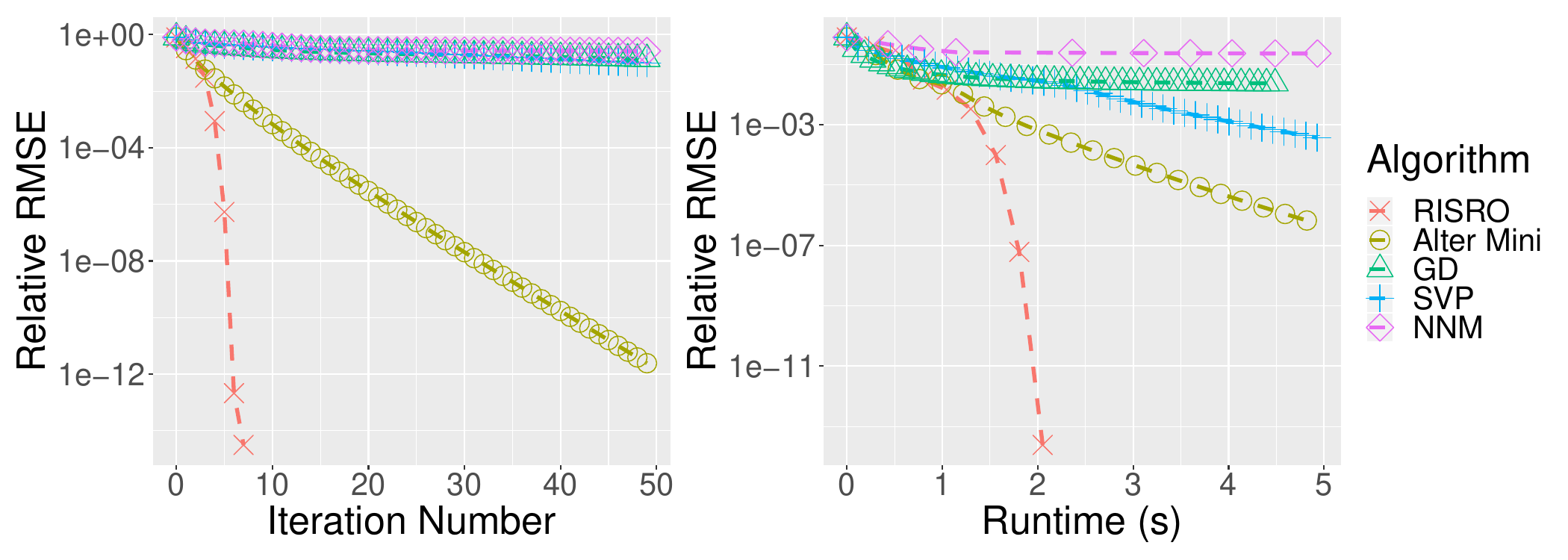}}
	\subfigure[$\kappa = 500$]{\includegraphics[width=0.7\textwidth]{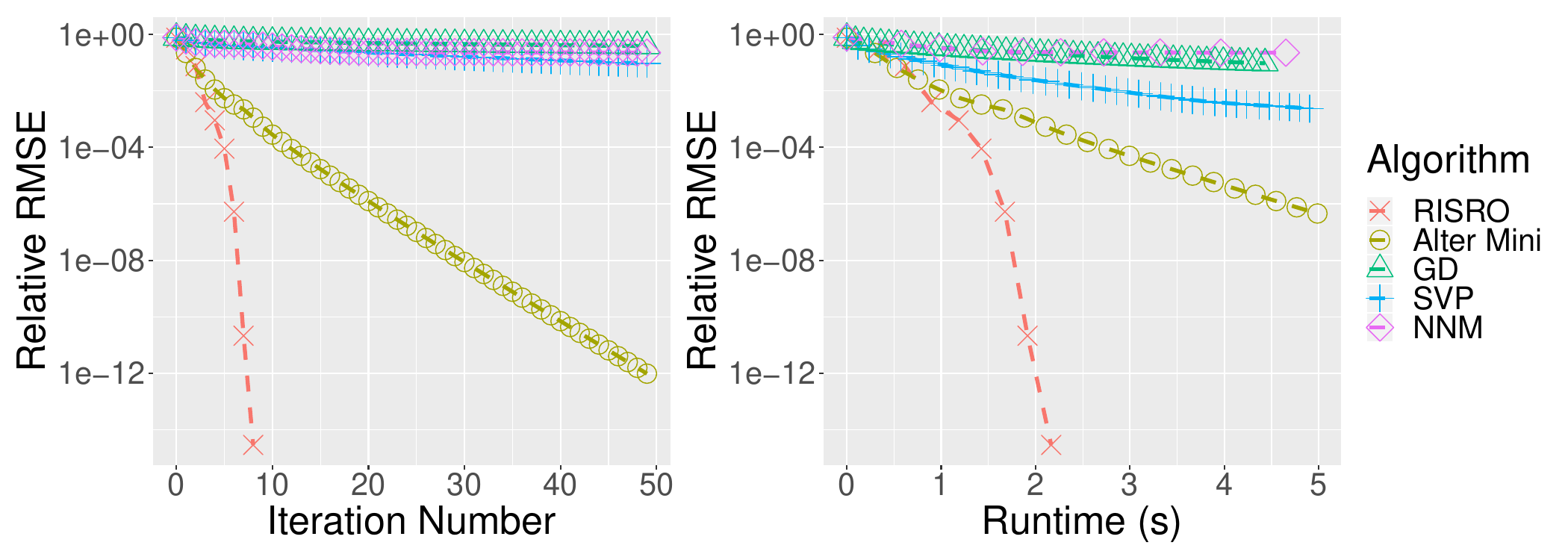}}
	\caption{Relative RMSE of RISRO, singular value projection (SVP), Alternating Minimization (Alter Mini), gradient descent (GD), and Nuclear Norm Minimization (NNM) in low-rank matrix trace regression. Here, $p = 100, r = 3, n = 1500, \sigma = 0, \kappa \in \{1,50,500 \}$.} \label{fig: Simul Mini comparison RMSE regression}
\end{figure}

We compare the performance of various procedures on noiseless matrix trace regression in Figure \ref{fig: Simul Mini comparison RMSE regression}. For all different choices of $\kappa$, RISRO converges quadratically to $\X^*$ in $7$ iterations with high accuracy, while the other baseline algorithms converge much slower at a linear rate. When $\kappa$ (condition number of $\X^*$) increases from $1$ to $50$ and $500$ so that the problem becomes more ill-conditioned, RISRO, Alter Mini, and SVP perform robustly, while GD converges more slowly. In Theorem \ref{th: local convergence in local rank matrix recovery}, we have shown the quadratic convergence rate of RISRO is robust to the condition number (see Remark \ref{rem: trace regression remark 1}). As we expect, the non-convex optimization methods converge much faster than the convex relaxation method. Moreover, to achieve a relative RMSE of $10^{-10}$, RISRO only takes about five iterations and $1/5$ runtime compared to other algorithms if $\kappa = 1$ and this factor is even smaller in the ill-conditioned cases that $\kappa=50$ and $500$. 

\begin{figure}[ht]
	\centering
	\subfigure{\includegraphics[width=0.7\textwidth]{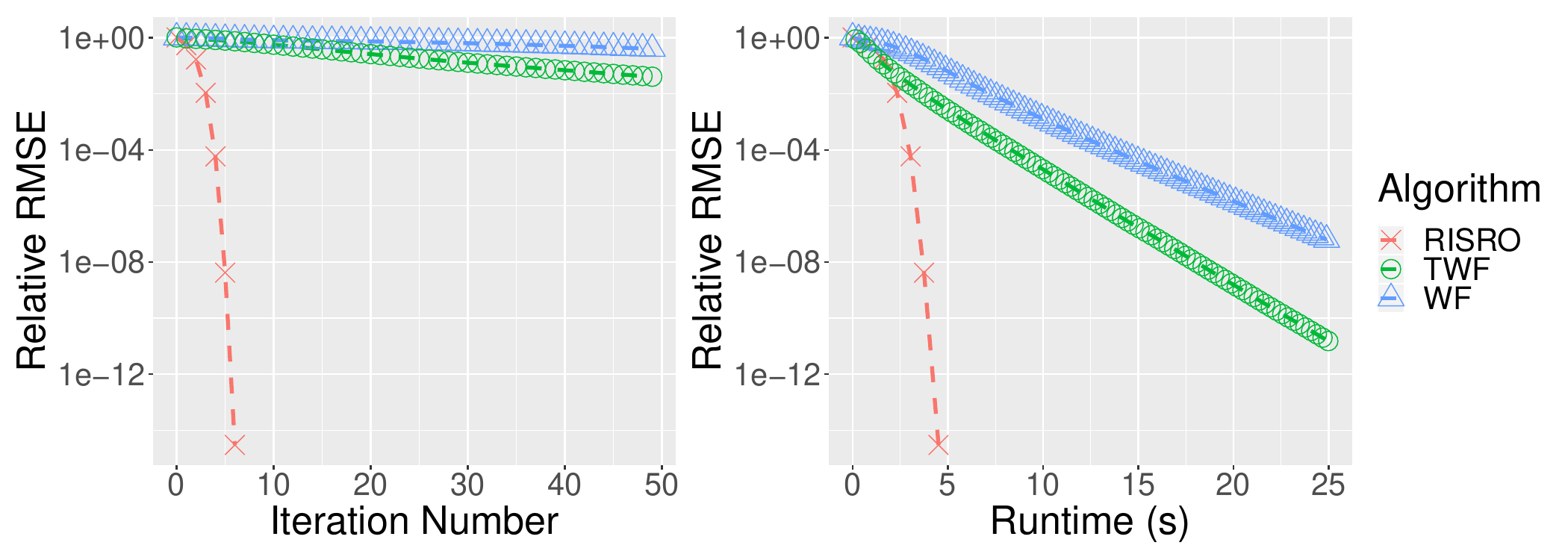}}
	\caption{Relative RMSE of RISRO, Wirtinger Flow (WF), Truncated Wirtinger Flow (TWF) in phase retrieval. Here, $p = 1200, n = 6000$} \label{fig: Simul Mini comparison RMSE phase retrieval}
\end{figure}

The comparison of RISRO, WF, and TWF in phase retrieval is plotted in Figure \ref{fig: Simul Mini comparison RMSE phase retrieval}. We can also see that RISRO can recover the underlying true signal with high accuracy in much less time than the other baseline methods. 

\begin{figure}[ht]
	\centering
	\subfigure{\includegraphics[width=0.7\textwidth]{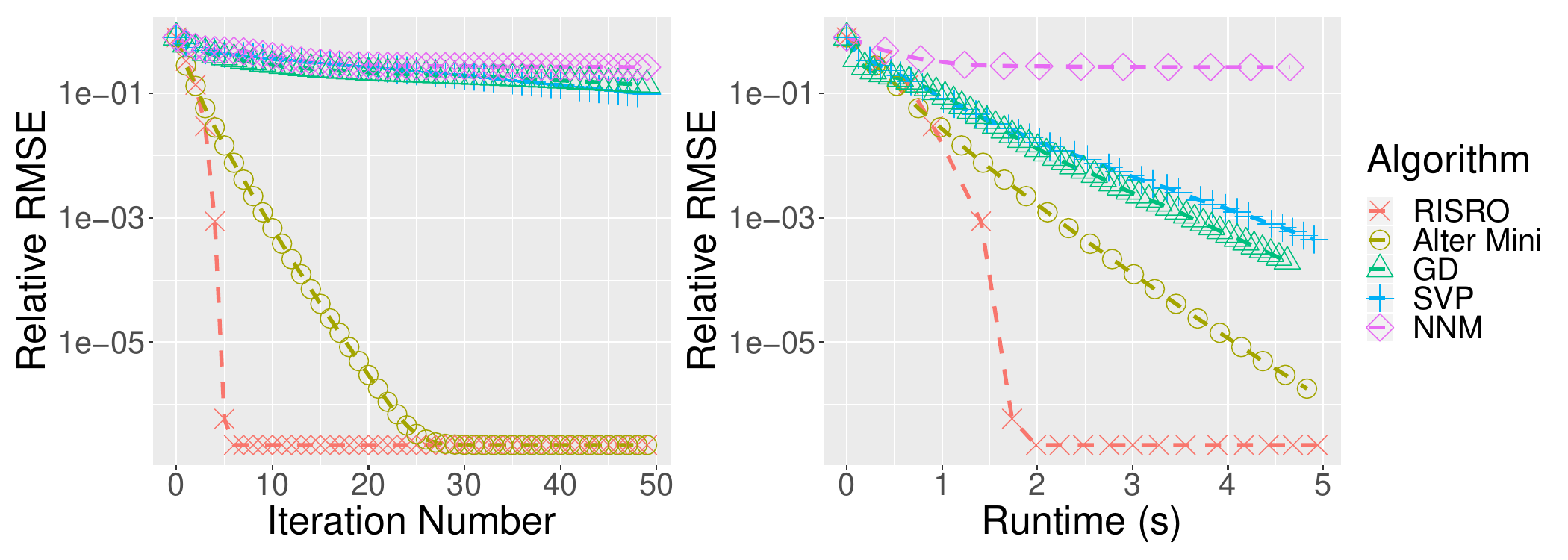}}
	\caption{Relative RMSE of RISRO, singular value projection (SVP), Alternating Minimization (Alter Mini), gradient descent (GD), and Nuclear Norm Minimization (NNM) in low-rank matrix trace regression. Here, $p = 100, r = 3, n = 1500, \kappa = 5, \sigma = 10^{-6}$} \label{fig: Simul Mini comparison RMSE noisy regression}
\end{figure}

Next, we compare the performance of RISRO with other algorithms in the noisy setting, $\sigma = 10^{-6}$, in the low-rank matrix trace regression. We can see from the results in Figure \ref{fig: Simul Mini comparison RMSE noisy regression} that due to the noise, the estimation error first decreases and then stabilizes after reaching a certain level. Meanwhile, we can also find RISRO converges at a much faster quadratic rate before reaching the stable level compared to all other algorithms. 

\begin{figure}[ht]
	\centering
	\subfigure[Matrix Trace Regression ($p = 100, r = 3, \sigma = 0, \kappa = 5$)]{\includegraphics[height=0.23\textwidth]{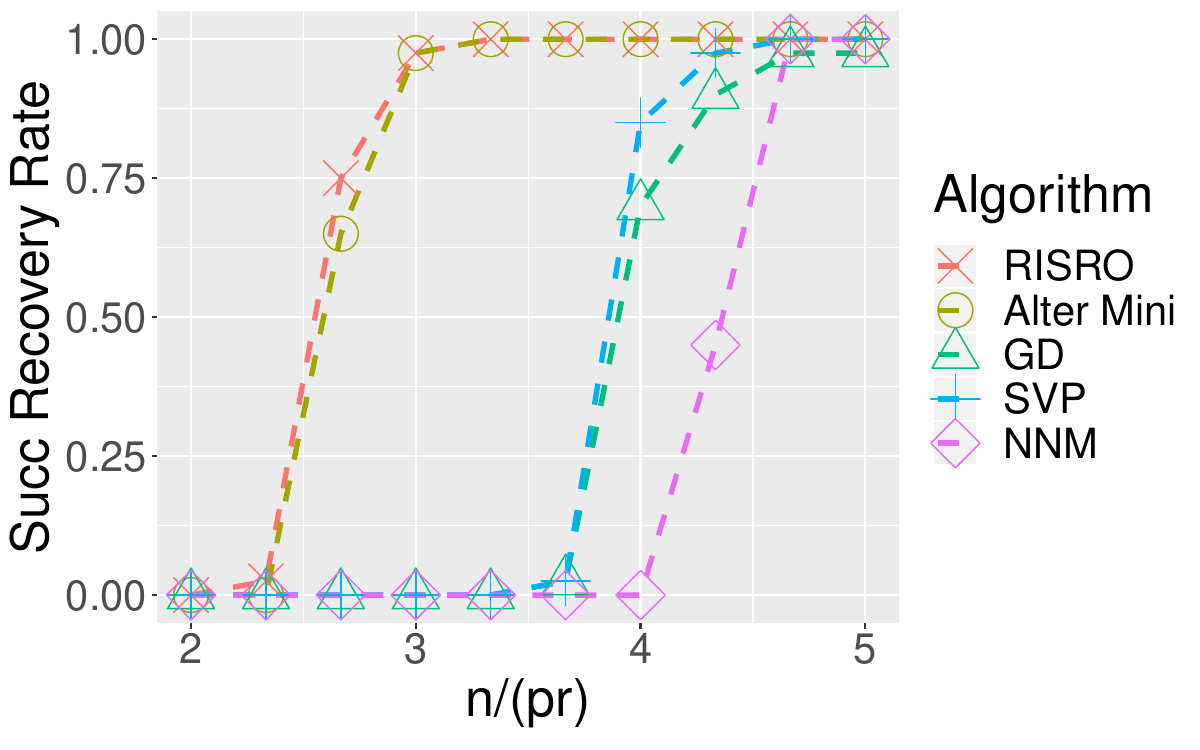}}
	\subfigure[Phase Retrieval ($p=1200$)]{\includegraphics[height=0.23\textwidth]{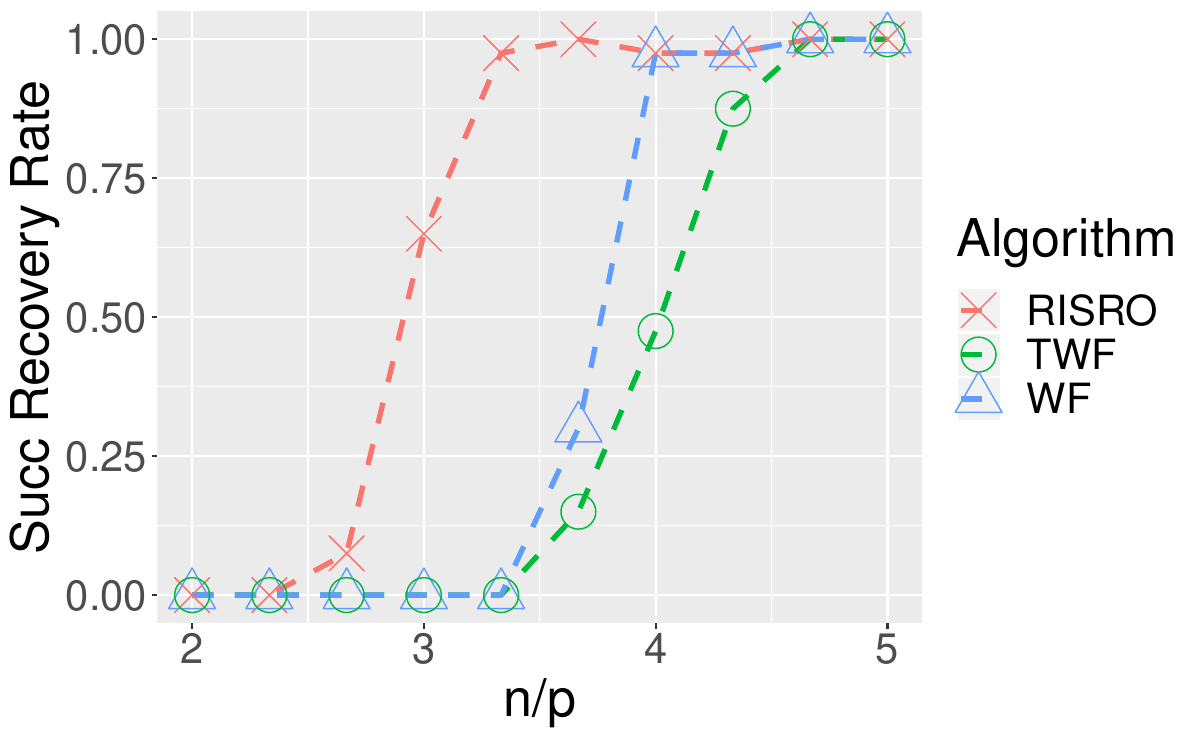}}
	\caption{Successful recovery rate comparison} 
	\label{fig: Simul Mini comparison succ recovery}
\end{figure}

Finally, we study the required sample size to guarantee successful recovery by RISRO and other algorithms. We set $p = 100, r = 3$, $\kappa = 5$, $n \in [600, 1500]$ in the noiseless matrix trace regression and $p = 1200$, $n \in [2400, 6000]$ in phase retrieval. We say the algorithm achieves successful recovery if the relative RMSE is less than $10^{-2}$ when the algorithm terminates. The simulation results in Figure \ref{fig: Simul Mini comparison succ recovery} show RISRO requires the minimum sample size to achieve a successful recovery in both matrix trace regression and phase retrieval; Alter Mini has similar performance to RISRO; and both RISRO and Alter Mini require smaller sample size than the rest of algorithms for successful recovery.

\section{Conclusion and Discussion} \label{sec: conclusion}

In this paper, we propose a new algorithm, RISRO, for solving rank constrained least squares. RISRO is based on a novel algorithmic framework, recursive importance sketching, which also provides new sketching interpretations for several existing algorithms for rank constrained least squares. RISRO is easy to implement and computationally efficient. Under some reasonable assumptions, local quadratic-linear and quadratic convergence are established for RISRO. Simulation studies demonstrate the superior performance of RISRO. 

The connection of recursive importance sketching and Riemannian Gauss-Newton discovered in this paper can be leveraged to other settings, such as in the low-rank tensor estimation problems (see a follow-up work in \cite{luo2021low} after the first preprint of this paper).

There are many interesting extensions to the results in this paper to be explored in the future. First, our current convergence theory on RISRO relies on the RIP assumption, which may not hold in many scenarios, such as phase retrieval, matrix completion, and robust PCA. In this paper, we give some theoretical guarantees of RISRO in phase retrieval with a strong initialization assumption. However, such an initialization requirement may be unnecessary and spectral initialization is good enough to guarantee quadratic convergence as we observe in the simulation studies. Empirically, we also observe RISRO achieves quadratic convergence in the matrix completion and robust PCA examples, see their development in Appendix \ref{sec: RISRO-matrix-completion-robust-PCA}. To improve and establish theoretical guarantees for RISRO in phase retrieval and matrix completion or robust PCA, we think more sophisticated analysis tools such as the ``leave-one-out'' method and some extra properties such as ``implicit regularization" \citep{ma2019implicit} need to be incorporated into the analysis and it will be interesting future work. Also, this paper focuses on the squared error loss in \eqref{eq:minimization}, while the other loss functions may be of interest in different settings, such as the $\ell_1$ loss in robust low-rank matrix recovery \citep{charisopoulos2019low,li2020nonconvex,li2020non}, which is worth exploring.

\section*{Acknowledgments.} We thank the editors and two anonymous reviewers for their suggestions and comments, which help significantly improve the presentation of this paper.

\bibliographystyle{apalike}
\bibliography{reference.bib}

%

%
%
\newpage

\appendix
\section{RISRO for Matrix Completion and Robust PCA} \label{sec: RISRO-matrix-completion-robust-PCA}
In this section, we provide implementation details and simulation results of RISRO in another two prominent examples:  
matrix completion \citep{candes2010power} and robust PCA \citep{candes2011robust}.

\vskip.3cm

\noindent{\bf RISRO for Matrix Completion}. In matrix completion, we observe a fraction of entries of an unknown low-rank matrix $\X^*$ (denote the index set of the observed entries by $\Omega$) corrupted by measurement error, $\Y_{\Omega}=\X^*_{\Omega} + \bvarepsilon_{\Omega}$, and aim to recover $\X^*$ from $\Y_{\Omega}$. Suppose the number of observed entries is $n$. Let $\{(i_k,j_k)^n_{k=1}: (i_k,j_k) \in \Omega\}$ be the ordered index pairs of the observed entries with the same order as the one in the vectorization of $\Y_{\Omega}$. Since sensing matrices in matrix completion are all binary with only one entry of $1$, RISRO for matrix completion can be simplified to Algorithm \ref{alg: RISRO-matrix-completion}. 
\begin{algorithm}[h]
\caption{RISRO for Matrix Completion}
	\begin{algorithmic}[1]
		\STATE Input: $\Y_{\Omega}$, ordered index pairs of the observed entries $\{ (i_k,j_k)^n_{k=1}: (i_k,j_k) \in \Omega \}$, rank $r$, initialization $\X^0$ that admits SVD $\U^0 \bSigma^0 \V^{0\top}$, where $\U^0 \in \bbO_{p_1,r}, \V^0 \in \bbO_{p_2, r}, \bSigma^0 \in \bbR^{r \times r}$
		\FOR{$t=0, 1, \ldots$} 
		\STATE Construct the importance covariates matrix $\A^t \in \bbR^{n \times (p_1 + p_2 -r)r}$:  
		\begin{equation*}
			\A^t_{[k,:]} = \left[ \V^t_{[j_k,:]} \otimes \U^t_{[i_k,:]}, \V^t_{[j_k,:]} \otimes \U^t_{\perp [i_k,:]}, \V^t_{\perp [j_k,:]} \otimes \U^t_{[i_k,:]} \right], \text{ for } k =1,\ldots,n = |\Omega|.
		\end{equation*} 
		\STATE Solve $( \rmvec(\B^{t+1}), \rmvec(\D_1^{t+1}), \rmvec(\D_2^{t+1\top})) = (\A^{t\top} \A^t)^{\dagger} \A^{t\top} \rmvec(\Y_{\Omega})$.
		\STATE Compute $\X^{t+1}_{U} = \left(\U^{t}\B^{t+1} + \U^{t}_{\perp}\D_1^{t+1}\right)$ and $\X^{t+1}_V = \left(\V^{t}\B^{t+1\top} + \V^{t}_{\perp}\D_2^{t+1}\right)$.
		\STATE  Perform QR orthogonalization: $\U^{t+1} = \QR(\X^{t+1}_U),\quad \V^{t+1} = \QR(\X^{t+1}_V).$
		\STATE  Update $\X^{t+1} = \X^{t+1}_U \left(\B^{t+1}\right)^{\dagger} \X_V^{t+1\top}$.
		\ENDFOR
	\end{algorithmic} \label{alg: RISRO-matrix-completion}
\end{algorithm}

\vskip.3cm

\noindent{\bf RISRO for Robust Principal Component Analysis}. The basic model of robust PCA is $\Y = \X^* + \S^* \in \mathbb{R}^{p_1\times p_2}$, where $\X^*$ is an unknown low-rank matrix of interest and $\S^*$ an unknown sparse corruption matrix. We consider the setting where a part of entries of $\Y$, whose indices are denoted by $\Omega \subseteq \{(i,j): 1\leq i\leq p_1, 1\leq j \leq p_2\}$, are observed. Our goal is to recover $\X^*$ based on $\Y_{\Omega}$. Define the following truncation operator $F: \bbR^{p_1 \times p_2} \to \bbR^{p_1 \times p_2}$ \citep{yi2016fast,zhang2018robust},
\begin{equation} \label{eq: threshold-operator}
	(F(\A))_{[i,j]} = \left\{ \begin{array}{ll}
		0, & \text{ if } |\A_{[i,j]}| > |\A_{[i,:]}|^{[\gamma, \Omega] }\text{ and } |\A_{[i,j]}| > |\A_{[:,j]}|^{[\gamma, \Omega] };\\
		\A_{[i,j]}, & \text{ otherwise}.
	\end{array} \right.
\end{equation} 
Here, $|\A_{[i,:]}|^{[\gamma, \Omega] }$ and $|\A_{[:,j]}|^{[\gamma, \Omega] }$ represent the $(1-\gamma)$-th percentile of the absolute values of the observed entries of $\A_{[i,:]}$ and $\A_{[:,j]}$ of the matrix $\A$, respectively. We provide an implementation of RISRO for robust PCA in Algorithm \ref{alg: RISRO-robust-PCA}.
\begin{algorithm}[h]
\caption{RISRO for Robust PCA}
	\begin{algorithmic}[1]
		\STATE Input: $\Y_{\Omega}$, rank $r$, initialization $\X^0$ which admits SVD $\U^0 \bSigma^0 \V^{0\top}$, $\U^0 \in \bbO_{p_1,r}, \V^0 \in \bbO_{p_2, r}, \bSigma^0 \in \bbR^{r \times r}$
		\FOR{$t=0, 1, \ldots$} 
		\STATE Denote the support of $F(\Y_{\Omega} - \X_{\Omega}^t )$ as $\Phi_t$. Construct the importance covariates matrix $\A^t \in \bbR^{n_t \times (p_1 + p_2 -r)r}$ where $n_t = |\Phi_t|$ is the cardinality of $\Phi_t$, and
		\begin{equation*}
			\A^t_{[k,:]} = \left[ \V^t_{[j_k,:]} \otimes \U^t_{[i_k,:]}, \V^t_{[j_k,:]} \otimes \U^t_{\perp [i_k,:]}, \V^t_{\perp [j_k,:]} \otimes \U^t_{[i_k,:]} \right], \text{ for } k =1,\ldots,n_t.
		\end{equation*} Here $\{ (i_k,j_k)^{n_t}_{k=1}: (i_k,j_k) \in \Phi_t \}$ are the ordered index pairs of $\Phi_t$.
		\STATE Solve $( \rmvec(\B^{t+1}), \rmvec(\D_1^{t+1}), \rmvec(\D_2^{t+1\top})) = (\A^{t\top} \A^t)^{\dagger} \A^{t\top} \rmvec(\Y_{\Phi_t})$.
		\STATE Compute $\X^{t+1}_{U} = \left(\U^{t}\B^{t+1} + \U^{t}_{\perp}\D_1^{t+1}\right)$ and $\X^{t+1}_V = \left(\V^{t}\B^{t+1\top} + \V^{t}_{\perp}\D_2^{t+1}\right)$.
		\STATE  Perform QR orthogonalization: $\U^{t+1} = \QR(\X^{t+1}_U),\quad \V^{t+1} = \QR(\X^{t+1}_V).$
		\STATE  Update $\X^{t+1} = \X^{t+1}_U \left(\B^{t+1}\right)^{\dagger} \X_V^{t+1\top}$.
		\ENDFOR
	\end{algorithmic} \label{alg: RISRO-robust-PCA}
\end{algorithm}

Next, we investigate the numerical performance of RISRO in matrix completion and robust PCA, i.e., Algorithms \ref{alg: RISRO-matrix-completion} and \ref{alg: RISRO-robust-PCA}. 
In both settings, we generate $\X^* \in \bbR^{p \times p}$ as a random rank-3 matrix with condition number $\kappa$ in the same way as the one in matrix trace regression described in Section \ref{sec:numerics}. In matrix completion, we set $p = 500, \kappa = \{1,50,500 \}$ and assume $n = 8 pr$ noiseless entries of $\X^*$ are observed uniformly at random. 
In robust PCA, we set $p = 100$, $\kappa = \{1,50,100\}$. We observe the full matrix $\Y = \X^* + \S^*$, where each entry of $\S^*$ follows $N(0,100)$ with probability $q = 0.02$ and equals zero with probability $1 - q$. The thresholding ratio $\gamma$ in the operator $F$ in \eqref{eq: threshold-operator} is set to be $15q$. 

\begin{figure}[h]
	\centering
	\subfigure[Matrix Completion]{\includegraphics[height=0.3\textwidth]{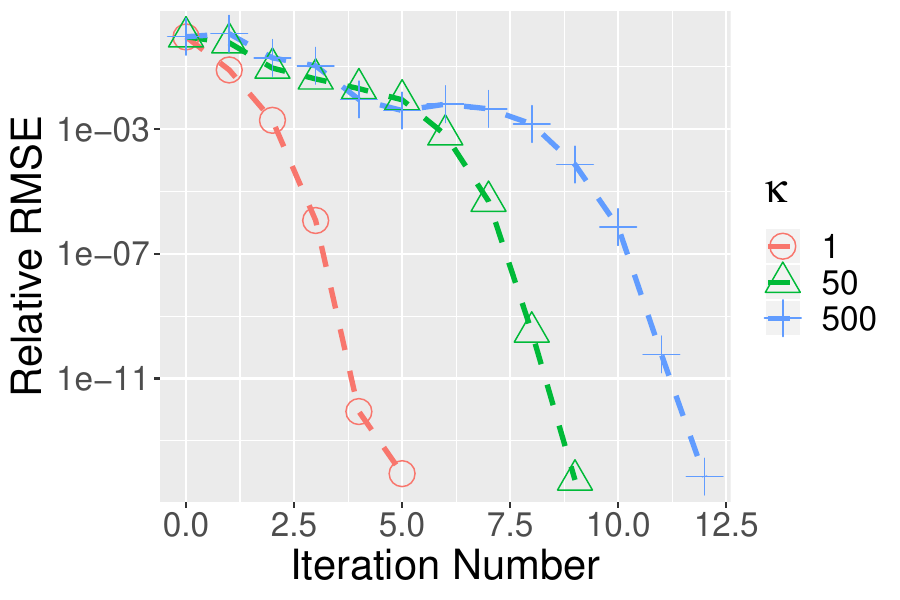}}
	\subfigure[Robust PCA]{\includegraphics[height=0.3\textwidth]{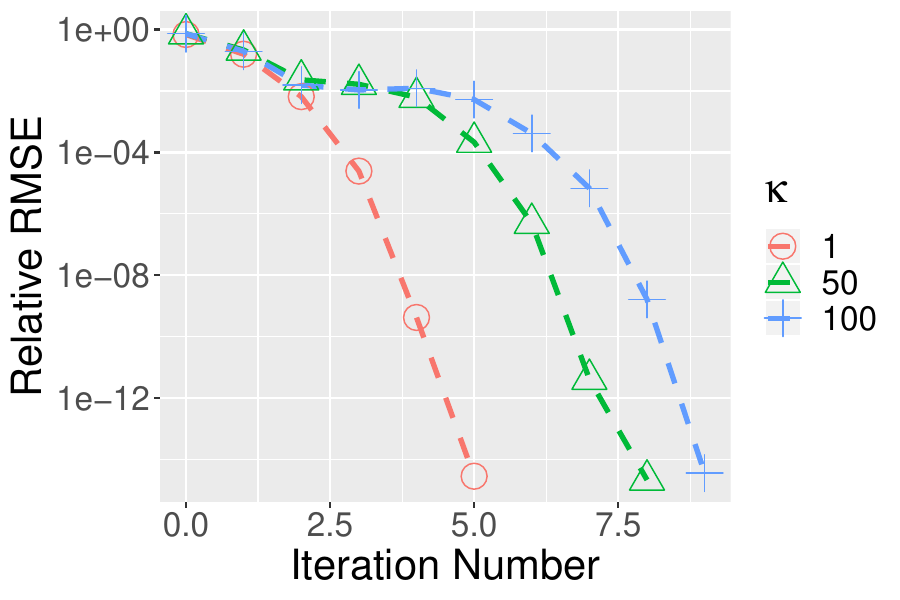}}
	\caption{Performance of RISRO in Matrix Completion and Robust PCA with Spectral Initialization.} 
	\label{fig: RISRO_matrix_completion_robust_PCA}
\end{figure}

The simulation results of RISRO with spectral initialization \citep{chi2019nonconvex} in matrix completion and robust PCA are provided in Figure \ref{fig: RISRO_matrix_completion_robust_PCA} (a) and (b), respectively. We observe that in both examples, RISRO converges quadratically to the true parameter $\X^*$ even though the RIP condition completely fails in matrix completion and robust PCA. Moreover, the performance of RISRO is robust even if $\X^*$ is ill-conditioned.

\section{Proof of the Main Results in the Paper}
 {\noindent \bf Proof of Lemma \ref{lm: partial least square for Simul Mini}}. 
First, by the decomposition of \eqref{eq: sketching interp of Simul Mini}, we have 
\begin{equation}
\label{eq:ydecomp}
\y = \cA \cL_t \left(\begin{bmatrix}
	\widetilde{\B}^t  & \widetilde{\D}_2^{t\top} \\
	\widetilde{\D}_1^t  & \0
	\end{bmatrix}\right) + \bepsilon^t.
\end{equation} In view of \eqref{eq: compact repre of least square}, if the operator $\cL_t^* \cA^* \cA \cL_t$ is invertible, the output of least squares in \eqref{eq: alg2 least square} satisfies
\begin{equation*}
	\begin{bmatrix}
		\B^{t+1}  & \D_2^{t+1 \top} \\
	\D_1^{t+1}  & \0
	\end{bmatrix} = (\cL_t^* \cA^* \cA \cL_t)^{-1} \cL_t^*\cA^* \y = \begin{bmatrix}
	\widetilde{\B}^t  & \widetilde{\D}_2^{t\top} \\
	\widetilde{\D}_1^t  & \0
	\end{bmatrix} + (\cL_t^* \cA^* \cA \cL_t)^{-1} \cL_t^*\cA^* \bepsilon^t,
\end{equation*} 
where the second equality is due to \eqref{eq:ydecomp}. This finishes the proof. \quad $\blacksquare$

{\noindent \bf Proof of Lemma \ref{lm: spectral norm bound of Atop A}}. 
Equation \eqref{eq: norm equal Range Lt} can be directly verified  from definitions of $\cL_t$ and $\cL_t^*$ in \eqref{eq: linear operator L_t}.

The second conclusion holds if $\M$ is a zero matrix. When $\M$ is not zero, to prove the claim it is equivalent to show the spectrum of $\cL_t^* \cA^* \cA \cL_t$ is upper and lower bounded by $1 + R_{2r}$ and $1-R_{2r}$, respectively, in the range of $\cL_t^*$. Since $\cL_t^* \cA^* \cA \cL_t$ is a symmetric operator, the upper and lower bounds of its spectrum are given below
\begin{equation*}
\begin{split}
	\sup_{\Z \in \Range(\cL^*_t): \|\Z\|_F = 1} \langle \Z, \cL_t^* \cA^* \cA \cL_t(\Z) \rangle = \sup_{\Z \in \Range(\cL^*_t):\|\Z\|_\F = 1} \|\cA \cL_t(\Z)\|_2^2 & \overset{(a)}\leq 1 + R_{2r} \\
	\inf_{\Z \in \Range(\cL^*_t): \|\Z\|_F = 1} \langle \Z, \cL_t^* \cA^* \cA \cL_t(\Z) \rangle = \inf_{\Z \in \Range(\cL^*_t):\|\Z\|_\F = 1} \|\cA \cL_t(\Z)\|_2^2 & \overset{(b)}\geq 1 - R_{2r}.
\end{split}
\end{equation*} Here (a) and (b) are due to the RIP condition of $\cA$, $\cL_t(\Z)$ is a at most rank $2r$ matrix by the definition of $\cL_t$ and \eqref{eq: norm equal Range Lt}. This finishes the proof.
\quad $\blacksquare$

{\noindent \bf Proof of Proposition \ref{prop: bound for iter approx error}}. Since $\cA$ satisfies 2r-RIP, $R_{2r}< 1$. Then, by Lemma \ref{lm: spectral norm bound of Atop A}, $\cL^*_t \cA^* \cA \cL_t$ is invertible over $\Range(\cL_t^*)$. 
With a slight abuse of
notation, define $P_{\X^t}$ as
\begin{equation} \label{eq: prop projection onto tangent space}
	P_{\X^t}(\Z) := \cL_t \cL_t^*(\Z) = P_{\U^t} \Z P_{\V^t} + P_{\U^t_\perp} \Z P_{\V^t} + P_{\U^t} \Z P_{\V^t_\perp},\quad \forall\, \Z\in\bbR^{p_1\times p_2},
\end{equation}
where $\U^t, \V^t$ are the updated sketching matrices at iteration $t$ defined in Step 7 of Algorithm \ref{alg: recursive IS alg 2}. We can verify $P_{\X^t}$ is an orthogonal projector. Let $P_{(\X^t)_\perp}(\Z) := \Z - P_{\X^t}(\Z) = P_{\U^t_\perp} \Z  P_{\V^t_\perp}$. Recall $\bepsilon^t = \cA (P_{\U^t_\perp} \widebar{\X} P_{\V^t_\perp} ) + \widebar{\bepsilon} = \cA P_{(\X^t)_\perp}(\widebar{\X}) + \widebar{\bepsilon} $ from \eqref{eq: sketching interp of Simul Mini}, we have
\begin{equation} \label{ineq: bound for the least square residual}
	\begin{split}
		 &\left\|(\cL^*_t \cA^* \cA \cL_t)^{-1} \cL_t^* \cA^* \bepsilon^t\right\|_\F^2\\
		 \overset{(a)}\leq &\frac{1}{(1 - R_{2r})^2} \|\cL_t^* \cA^* \left( \cA(P_{(\X^t)_\perp} \widebar{\X}) + \widebar{\bepsilon} \right)\|_\F^2\\
		 = &\frac{1}{(1 - R_{2r})^2} \left( \overbrace{\|\cL_t^* \cA^* \cA(P_{(\X^t)_\perp} \widebar{\X})  \|_\F^2}^{(1)} + \overbrace{\| \cL_t^* \cA^* (\widebar{\bepsilon})\|_\F^2}^{(2)} + \overbrace{2 \left\langle \cL_t^* \cA^* \cA P_{(\X^t)_\perp} \widebar{\X}, \cL_t^* \cA^*(\widebar{\bepsilon}) \right\rangle}^{(3)} \right),
	\end{split}
\end{equation} where (a) is due to Lemma \ref{lm: spectral norm bound of Atop A} and the definition of $\bepsilon^t$.

Notice term $(2)$ is the target term we want, next, we bound (1) and (3) at the right-hand side of \eqref{ineq: bound for the least square residual}.

{\noindent \bf Bound for (1)}.
\begin{equation} \label{ineq: least square error term 1 bound}
	\begin{split}
		\|\cL_t^* \cA^* \cA(P_{(\X^t)_\perp} \widebar{\X})  \|_\F^2 & = \langle \cL_t^* \cA^* \cA(P_{(\X^t)_\perp} \widebar{\X}), \cL_t^* \cA^* \cA(P_{(\X^t)_\perp} \widebar{\X}) \rangle\\
		& = \langle \cA(P_{(\X^t)_\perp} \widebar{\X}) , \cA P_{\X^t} \cA^* \cA(P_{(\X^t)_\perp} \widebar{\X})\rangle\\
		& \overset{(a)}\leq R_{3r} \| P_{(\X^t)_\perp} \widebar{\X} \|_\F \|P_{\X^t} \cA^* \cA(P_{(\X^t)_\perp} \widebar{\X})\|_\F\\
		& = R_{3r} \| P_{(\X^t)_\perp} \widebar{\X} \|_\F \|\cL_t^* \cA^* \cA(P_{(\X^t)_\perp} \widebar{\X})\|_\F,
	\end{split}
\end{equation} where (a) is due to the Lemma \ref{lm:retricted orthogonal property} and the fact that $\langle P_{(\X^t)_\perp} \widebar{\X}, P_{\X^t} \cA^* \cA P_{(\X^t)_\perp} \widebar{\X} \rangle = 0 $, $\rank(P_{(\X^t)_\perp} \widebar{\X}) \leq r$ and $\rank(P_{\X^t} \cA^* \cA P_{(\X^t)_\perp} \widebar{\X}) \leq 2r$. 
Note that
\begin{equation} \label{eq: decomposition of X-bar complement}
\begin{split}
	P_{(\X^t)_\perp} \widebar{\X} &= \widebar{\X} - P_{\X^t} \widebar{\X}\\
	& \overset{(a)}= P_{\widebar{\U}}\widebar{\X} + \widebar{\X} P_{\widebar{\V}} - P_{\widebar{\U}}\widebar{\X}P_{\widebar{\V}} - P_{\U^t} \widebar{\X} - \widebar{\X}P_{\V^t} + P_{\U^t} \widebar{\X} P_{\V^t}\\
	& = (P_{\widebar{\U}} - P_{\U^t})\widebar{\X} + \widebar{\X}(P_{\widebar{\V}} - P_{\V^t}) - P_{\widebar{\U}}\widebar{\X}P_{\widebar{\V}} + P_{\widebar{\U}}\widebar{\X}P_{\V^t}-P_{\widebar{\U}}\widebar{\X}P_{\V^t} + P_{\U^t} \widebar{\X} P_{\V^t}\\
	& = (P_{\widebar{\U}} - P_{\U^t})\widebar{\X} (\I - P_{\V^t}) + (\I - P_{\widebar{\U}}) \widebar{\X} (P_{\widebar{\V}} - P_{\V^t})\\
	& \overset{(b)}= (P_{\widebar{\U}} - P_{\U^t})\widebar{\X} (\I - P_{\V^t})\\
	& \overset{(c)} = (P_{\widebar{\U}} - P_{\U^t})(\widebar{\X}- \X^t) (\I - P_{\V^t}),
\end{split}
\end{equation} where $\widebar{\U}, \widebar{\V}$ are left and right singular vectors of $\widebar{\X}$, (a) is because $\widebar{\X} =P_{\widebar\U} \widebar\X +  \widebar\X P_{\widebar\V} - P_{\widebar\U} \widebar\X P_{\widebar\V}$, (b) is due to the fact that $(\I - P_{\widebar{\U}})\widebar{\X} (P_{\widebar{\V}} - P_{\V^t}) = 0$ and (c) is because $\X^t(\I - P_{\V^t} ) = 0$.
Now, from \eqref{ineq: least square error term 1 bound}, we get
\begin{equation}\label{ineq: bound for Thm 3 term (1)}
\begin{split}
	\|\cL_t^* \cA^* \cA(P_{(\X^t)_\perp} \widebar{\X})  \|_\F &\leq R_{3r}  \| P_{(\X^t)_\perp} \widebar{\X} \|_\F \\
	& \overset{(a)} \leq R_{3r} \|(P_{\widebar{\U}} - P_{\U^t})(\widebar{\X}- \X^t) (\I - P_{\V^t})\|_\F \\
	& \leq R_{3r} \|P_{\widebar{\U}} - P_{\U^t}\| \|\widebar{\X}- \X^t\|_\F  \\
	& \overset{(b)}\leq R_{3r} \frac{\|\X^t - \widebar{\X}\| \|\X^t - \widebar{\X}\|_\F }{\sigma_r(\widebar{\X})},
\end{split}
\end{equation} here (a) is due to \eqref{eq: decomposition of X-bar complement} and (b) is due to the singular subspace perturbation inequality $ \|P_{\widebar{\U}} - P_{\U^t}\| \leq \|\X^t- \widebar{\X}\|/\sigma_r(\widebar{\X}) $ obtained from Lemma \ref{lm: perturbation bound}.

{\noindent \bf Bound for (3)}. 
\begin{equation} \label{ineq: bound for Thm 2, term (3)}
	\begin{split}
		2 \left\langle \cL_t^* \cA^* \cA P_{(\X^t)_\perp} \widebar{\X}, \cL_t^* \cA^*(\widebar{\bepsilon}) \right\rangle &= 2 \langle \cA P_{(\X^t)_\perp} \widebar{\X} , \cA P_{\X^t} \cA^*(\widebar{\bepsilon}) \rangle\\
		& \overset{(a)}\leq2 R_{3r} \| P_{(\X^t)_\perp} \widebar{\X} \|_\F \|  \cL^*_t \cA^*(\widebar{\bepsilon})\|_\F\\
		& \overset{(b)} \leq 2 R_{3r}  \frac{\|\X^t - \widebar{\X}\| \|\X^t - \widebar{\X}\|_\F }{\sigma_r(\widebar{\X})} \|  \cL^*_t \cA^*(\widebar{\bepsilon})\|_\F,
	\end{split}
\end{equation} here (a) and (b) are by the same arguments in \eqref{ineq: least square error term 1 bound}  and  \eqref{ineq: bound for Thm 3 term (1)}, respectively.

Plugging \eqref{ineq: bound for Thm 2, term (3)} and \eqref{ineq: bound for Thm 3 term (1)} into 
\eqref{ineq: bound for the least square residual}, we get \eqref{ineq: prop least square residual bound}. This finishes the proof of this Proposition. \quad $\blacksquare$

{\noindent \bf Proof of Theorem \ref{th: local contraction general setting}}. 
 First notice that quadratic convergence result follows easily from \eqref{ineq: optimization theory quadratic and linear converge} when we set $\widebar{\bepsilon} = 0$. So the rest of the proof is devoted to proving \eqref{ineq: optimization theory quadratic and linear converge} and we also prove the linear convergence along the way. The proof can be divided into four steps. In Step 1, we use Proposition \ref{prop: bound for iter approx error} and give an upper bound for the approximation error in the case $\widebar{\X}$ is a stationary point. Then, we use induction to show the main results in Step 2,3,4.
To start, similar as \eqref{eq: prop projection onto tangent space}, define 
\[
P_{\widebar{\X}}(Z) = P_{\widebar \U} \Z P_{\widebar \V} + P_{\widebar \U_\perp} \Z P_{\widebar \V} + P_{\widebar \U} \Z P_{\widebar \V_\perp},\quad \forall\, \Z\in\bbR^{p_1\times p_2},
\]
where  $\widebar{\U}, \widebar{\V}$ are left and right singular vectors of $\widebar{\X}$.

{\noindent \bf Step 1}. In this step, we apply Proposition \ref{prop: bound for iter approx error} in the case $\widebar{\X}$ is a stationary point. In view of \eqref{ineq: prop least square residual bound}, the term that we can simplify is $\|\cL_t^* \cA^* (\widebar{\bepsilon})\|_\F$. Since $\widebar{\X}$ is a stationary point, we know $P_{\widebar{\X}} \big(\cA^*(\widebar{\bepsilon})\big) = P_{\widebar{\X}} \big(\cA^*(\y - \cA(\widebar \X)) \big) = 0$. Then
\begin{equation} \label{ineq: bound for Thm 2, term 2, first inequality}
	\begin{split}
		\| \cL_t^* \cA^* (\widebar{\bepsilon})\|_\F^2 = \| P_{\X^t} \cA^* (\widebar{\bepsilon}) \|_\F^2 = \| (P_{\X^t} - P_{\widebar{\X}}) \cA^*(\widebar{\bepsilon}) \|_\F^2 &\leq \| P_{\X^t} - P_{\widebar{\X}} \|^2 \| \cA^*(\widebar{\bepsilon}) \|_\F^2,
	\end{split}
\end{equation} where the first equality is by the definition of $P_{\X^t}$ in \eqref{eq: prop projection onto tangent space} and  equation \eqref{eq: norm equal Range Lt}. Meanwhile, it holds that
\begin{equation*}
	\begin{split}
		\| P_{\X^t} - P_{\widebar{\X}} \| &= \sup_{\|\Z\|_\F \leq 1} \| (P_{\X^t} - P_{\widebar{\X}})(\Z) \|_\F\\
		& \overset{(a)}= \sup_{ \|\Z\|_\F \leq 1} \|(P_{\U^t} - P_{\widebar{\U}} ) \Z (\I - P_{\widebar{\V}}) \|_\F + \| (\I - P_{\U^t}) \Z (P_{\V^t} - P_{\widebar{\V}})  \|_\F\\
		& \overset{(b)}\leq \sup_{ \|\Z\|_\F \leq 1} \frac{2 \|\X^t - \widebar{\X}\|}{\sigma_r(\widebar{\X})} \|\Z\|_\F \leq \frac{2 \|\X^t - \widebar{\X}\|}{\sigma_r(\widebar{\X})},
	\end{split}
\end{equation*} here (a) is because $(P_{\X^t} - P_{\widebar{\X}})(\Z) = (P_{\U^t}-P_{\widebar{\U}})\Z (\I - P_{\V^t}) + (\I - P_{\widebar{\U}}) \Z (P_{\V^t}-P_{\widebar{\V}} )$ by a similar argument in \eqref{eq: decomposition of X-bar complement} and (b) is by Lemma \ref{lm: perturbation bound}. Then, from \eqref{ineq: bound for Thm 2, term 2, first inequality}, we see that
\[
\| \cL_t^* \cA^* (\widebar{\bepsilon})\|_\F^2 \leq \frac{4\|\X^t - \widebar{\X} \|^2}{\sigma_r^2(\widebar{\X})} \|\cA^*(\widebar{\bepsilon})\|_\F^2,
\]
which, together with \eqref{ineq: prop least square residual bound}, implies that
\begin{equation} \label{ineq: least square residual upper bound}
	\begin{split}
	\|(\cL^*_t \cA^* & \cA \cL_t)^{-1} \cL_t^* \cA^* \bepsilon^t\|_\F^2 \leq  \frac{\|\X^t - \widebar{\X}\|^2}{(1-R_{2r})^2 \sigma^2_r(\widebar{\X})}\\
	&\quad \quad \quad \cdot \left( R_{3r}^2 \|\X^t - \widebar{\X}\|_\F^2 + 4 \|\cA^*(\widebar{\bepsilon})\|_\F^2 + 4 R_{3r}\|\cA^*(\widebar{\bepsilon})\|_\F \|\X^t - \widebar{\X}\|_\F \right).
	\end{split}
\end{equation}

{\noindent \bf Step 2.} In this step, we start using induction to show the convergence of $\X^t$ in \eqref{ineq: optimization theory quadratic and linear converge}. First, we introduce $\theta_{t}$ to measure the goodness of the current iterate in terms of subspace estimation:
\begin{equation} \label{ineq: sin theta bound}
\begin{split}
	\theta_t = \max\{ \left\| \sin \Theta(\U^{t}, \widebar{\U})  \right\|, \left\| \sin \Theta(\V^{t},\widebar{\V})  \right\| \},
	\end{split}	
\end{equation}where $\widebar{\U},\widebar{\V}$ are the left and right singular vectors of $\widebar{\X}$ and here $\|\sin\Theta(\U^t, \widebar{\U})\| := |\sin (\cos^{-1}(\sigma_r(\U^{t\top} \widebar{\U})))|$ measures the largest angle between subspaces $\U^t, \widebar{\U}$.

Recall the definition of $\widetilde{\B}^{t}, \widetilde{\D}_1^t $ and $\widetilde{\D}_2^t$ in \eqref{eq: tildeB,D}. The induction we want to show is: given $\theta_t \leq 1/2$, $\widetilde{\B}^t$ is invertible and $\|\X^{t} - \widebar{\X}\|_\F \leq \|\X^0 - \widebar{\X}\|_\F$, we prove  $\theta_{t+1} \leq 1/2$, $\widetilde{\B}^{t+1}$ is invertible, $\|\X^{t+1} - \widebar{\X}\|_\F \leq \|\X^0 - \widebar{\X}\|_\F$ as well as $\B^{t+1}$ is invertible and 
\begin{equation} \label{ineq: Xt+1 - barX bound}
\begin{split}
	\|\X^{t+1} - \widebar{\X}\|_\F^2& \leq 5 \left\|(\cL^*_t \cA^* \cA \cL_t)^{-1} \cL_t^* \cA^* \bepsilon^t\right\|_\F^2\\
	&\leq \frac{5\|\X^t - \widebar{\X}\|^2}{(1-R_{2r}^2) \sigma^2_r(\widebar{\X})} \left( R_{3r}^2 \|\X^t - \widebar{\X}\|_\F^2 + 4 R_{3r} \|\cA^*(\widebar{\bepsilon})\|_\F\|\X^t - \widebar{\X}\|_\F + 4\|\cA^*(\widebar{\bepsilon})\|^2_\F   \right)\\
	& \leq \frac{9}{16} \|\X^t - \widebar{\X} \|_\F^2.
\end{split}
\end{equation} 

For the rest of this step, we show when $t = 0$, the induction assumption holds, i.e. $\theta_0 \leq 1/2$ and $\widetilde{\B}^0$ is invertible. Since $\frac{\|\X^0 - \widebar{\X}\|}{\sigma_r(\widebar{\X})} \leq \frac{1}{4}$, it holds by Lemma \ref{lm: perturbation bound} that $\theta_0 = \max\{\|\sin \Theta(\U^0, \widebar{\U})\|,\|\sin \Theta(\V^0, \widebar{\V})\|\} \leq 2\frac{\|\X^0 - \widebar{\X}\|}{\sigma_r(\widebar{\X})} \leq \frac{1}{2}$. Then, we have
\begin{equation}
\label{eq: sigmarB0 g0}
	\sigma_r(\widetilde{\B}^0) \overset{(a)}\geq \sigma_r(\U^{0\top} \widebar{\U}) \sigma_r(\widebar{\X}) \sigma_r (\widebar{\V}^\top \V^0) \overset{(b)}= (1 - \theta_0^2) \sigma_r (\widebar{\X}) \geq 3/4 \cdot  \sigma_r (\widebar{\X}) > 0,
\end{equation} 
where (a) is because $\U^{0\top} \widebar{\U},\widebar{\X},\widebar{\V}^\top \V^0$ are all rank $r$ matrices and \cite[arXiv v1,Lemma 5]{luo2020schatten} and (b) is by the definition of $\| \sin \Theta (\U^t, \widebar{\U} ) \|$.

{\noindent \bf Step 3.} In this step, we first show given $\theta_t \leq 1/2$ and $\|\X^{t} - \widebar{\X}\|_\F \leq \|\X^0 - \widebar{\X}\|_\F$, $\B^{t+1}$ is invertible. Then we also do some preparation to show the main contraction \eqref{ineq: Xt+1 - barX bound}. Notice that
\begin{equation} \label{ineq: lower bound on Bt+1}
\begin{split}
	\sigma_r(\B^{t+1}) \geq \sigma_r(\widetilde{\B}^t) - \| \B^{t+1} - \widetilde{\B}^t \| &\geq \sigma_r(\U^{t\top} \widebar{\U} ) \sigma_r(\widebar{\X}) \sigma_r(\V^{t\top} \widebar{\V}) -  \| \B^{t+1} - \widetilde{\B}^t \|\\
	& \geq (1 - \theta_t^2) \sigma_r(\widebar{\X}) - \| \B^{t+1} - \widetilde{\B}^t \|,
\end{split}
\end{equation} and
\begin{equation} \label{ineq: B, D1, D2 spectral norm bound}
	\begin{split}
		\max\{ \| \B^{t+1} - \widetilde{\B}^t \|, \| \D_1^{t+1} - \widetilde{\D}_1^t \|, \|\D_2^{t+1} - \widetilde{\D}_2^t \|  \} &\overset{(a)}\leq \left\|(\cL^*_t \cA^* \cA \cL_t)^{-1} \cL_t^* \cA^* \bepsilon^t\right\|\\
		& \leq\left\|(\cL^*_t \cA^* \cA \cL_t)^{-1} \cL_t^* \cA^* \bepsilon^t\right\|_\F,
	\end{split}
\end{equation} 
 where (a) is due to \eqref{eq: Bt+1 - Btilde}.

Under the induction assumption $\|\X^{t} - \widebar{\X}\|_\F \leq \|\X^0 - \widebar{\X}\|_\F$, the initialization condition \eqref{ineq: optimization theory initialization condition} and $\|\cA^*(\widebar{\bepsilon})\|_\F \leq \frac{1-R_{2r}}{4\sqrt{5}}\sigma_r(\widebar{\X})$, we have
\begin{equation*}
	\begin{split}
		\frac{R_{3r}^2\|\X^t - \widebar{\X}\|_\F^2}{(1-R_{2r})^2 \sigma^2_r(\widebar{\X}) } \leq 1/80,\,
		\frac{4 \|\cA^*(\widebar{\bepsilon})\|_\F^2}{(1-R_{2r})^2 \sigma^2_r(\widebar{\X}) }  \leq 1/20, \, \frac{4 R_{3r}\|\X^t - \widebar{\X}\|_\F \|\cA^*(\widebar{\bepsilon})\|_\F  }{(1-R_{2r})^2 \sigma^2_r(\widebar{\X}) } \leq 1/20.
	\end{split}
\end{equation*} 
By combining these inequalities with \eqref{ineq: least square residual upper bound}, we obtain
\begin{equation*}
	\left\|(\cL^*_t \cA^* \cA \cL_t)^{-1} \cL_t^* \cA^* \bepsilon^t\right\|_\F \leq \frac{3}{4\sqrt{5}} \|\X^0 - \widebar{\X}\| \leq \frac{3}{16} \sigma_r(\widebar{\X}).
\end{equation*} 
Thus, from \eqref{ineq: B, D1, D2 spectral norm bound} we have  
\begin{equation} \label{ineq: Bt+1 - tildeBt}
	\max\{ \| \B^{t+1} - \widetilde{\B}^t \|, \| \D_1^{t+1} - \widetilde{\D}_1^t \|, \|\D_2^{t+1} - \widetilde{\D}_2^t \|  \} \leq \frac{3}{16} \sigma_r(\widebar{\X}),
\end{equation} and $ \sigma_r(\B^{t+1}) \geq \frac{9}{16} \sigma_r(\widebar{\X}) > 0 $ because of \eqref{ineq: lower bound on Bt+1}. This shows the invertibility of $\B^{t+1}$.

With the invertibility of $\B^{t+1}$, we also introduce $\rho_{t+1}$ to measure the goodness of the current iterate in the following way
\begin{equation} \label{def: rho t+1}
\begin{split}
		 \rho_{t+1} = \max \{ \left\| \D_1^{t+1} (\B^{t+1})^{-1} \right\|, \left\| (\B^{t+1})^{-1} \D_2^{t+1\top} \right\|  \} .
\end{split}	
\end{equation}
Notice $\widetilde{\D}_1^t (\widetilde{\B}^t)^{-1} = \U^{t\top}_\perp \widebar{\X} \V^t (\U^{t\top} \widebar{\X} \V^t )^{-1} = \U^{t\top}_\perp \widebar{\U} (\U^{t\top} \widebar{\U})^{-1}$, so
\begin{equation} \label{ineq: bound for D1 B(-1)}
	\begin{split}
		\| \widetilde{\D}_1^t (\widetilde{\B}^t)^{-1}\| \leq \|\U^{t\top}_\perp \widebar{\U}\| \|(\U^{t\top} \widebar{\U})^{-1} \| \overset{(a)}= \frac{\|\sin \Theta(\U^t, \widebar{\U} )\|}{ \sqrt{1 - \|\sin \Theta(\U^t, \widebar{\U} ) \|^2} } \overset{(b)}\leq \frac{
	\theta_t}{\sqrt{1 - \theta_t^2}},
	\end{split}
\end{equation} where $(a)$ is due to the $\sin \Theta$ property in Lemma 1 of \cite{cai2018rate} and (b) is due to \eqref{ineq: sin theta bound}. The same bound also holds for $\|  (\widetilde{\B}^t)^{-1}\widetilde{\D}_2^{t\top}\|$.
Meanwhile, it holds that 
\begin{equation} \label{ineq: rho t+1 bound}
	\begin{split}
		\|\D_1^{t+1} (\B^{t+1})^{-1}\| \leq& \| (\D_1^{t+1} - \widetilde{\D}_1^{t} )(\B^{t+1})^{-1} \| + \|\widetilde{\D}_1^t (\B^{t+1})^{-1}\|\\
		\overset{(a)} \leq & \frac{\| \D_1^{t+1} - \widetilde{\D}_1^{t} \|}{\sigma_r(\B^{t+1})} + \|\widetilde{\D}_1^t (\widetilde{\B}^{t})^{-1}\| + \|\widetilde{\D}_1^t (\widetilde{\B}^{t})^{-1} (\B^{t+1} - \widetilde{\B}^t ) (\B^{t+1})^{-1}  \|\\
		\overset{(b)} \leq & \frac{\| \D_1^{t+1} - \widetilde{\D}_1^{t} \|}{\sigma_r(\B^{t+1})} + \frac{\theta_t}{\sqrt{1 - \theta_t^2}} + \frac{\theta_t}{\sqrt{1 - \theta_t^2}} \frac{\|\B^{t+1} - \widetilde{\B}^t \|}{\sigma_r(\B^{t+1})},
	\end{split}
\end{equation} 
where (a) is because $(\B^{t+1})^{-1} = (\widetilde{\B}^t)^{-1} - (\widetilde{\B}^t)^{-1} (\B^{t+1} - \widetilde{\B}^t) (\B^{t+1})^{-1}$ and (b) is due to the bound for $\|\widetilde{\D}_1^t (\widetilde{\B}^{t})^{-1}\|$ in \eqref{ineq: bound for D1 B(-1)}. We can also bound $ \left\| (\B^{t+1})^{-1} \D_2^{t+1\top} \right\|$ in a similar way.

Thus, by plugging $\sigma_r(\B^{t+1}) \geq \frac{9}{16} \sigma_r(\widebar{\X})$, $\theta_t \leq 1/2$ and the upper bound in \eqref{ineq: Bt+1 - tildeBt} into \eqref{ineq: rho t+1 bound}, we have $$\rho_{t+1} \leq \frac{1}{3} + \frac{1}{\sqrt{3}} + \frac{1}{3\sqrt{3}} \leq \frac{4 + \sqrt{3}}{3\sqrt{3}}. $$

{\noindent \bf Step 4}. In the last step, we show the contraction inequality \eqref{ineq: Xt+1 - barX bound}. Note from Lemma \ref{lm: perturbation bound} that $\theta_t$ decreases as $\|\X^{t} - \widebar{\X}\|_\F$ decreases. So after we show the contraction of $\|\X^{t+1} - \widebar{\X}\|_\F$, it automatically implies that $\theta_{t+1} \leq 1/2$, which, together with similar arguments in \eqref{eq: sigmarB0 g0}, further guarantees the invertibility of  $\widetilde{\B}^{t+1}$.

Since $\widebar{\X}$ is a rank $r$ matrix and $\widetilde{\B}^t,\B^{t+1}$ are invertible, a quick calculation asserts $\widetilde{\D}_1^t (\widetilde{\B}^t)^{-1} \widetilde{\D}_2^{t\top} = \U^{t\top}_\perp \widebar{\X} \V^{t\top}_\perp $. Therefore,
\begin{equation*}
	\begin{split}
		\widebar{\X} &= [\U^t \quad \U^{t}_\perp] [\U^t \quad \U^{t}_\perp]^\top \widebar{\X} [\V^t \quad \V^t_\perp] [\V^t \quad \V^t_\perp]^\top\\
		& = [\U^t \quad \U^{t}_\perp] \begin{bmatrix}
\widetilde{\B}^{t}  & \widetilde{\D}_2^{t\top} \\
\widetilde{\D}_1^{t} & \widetilde{\D}_1^{t} (\widetilde{\B}^{t})^{-1} \widetilde{\D}_2^{t\top}
\end{bmatrix}  [\V^t \quad \V^t_\perp]^\top.
	\end{split}
\end{equation*}
At the same time, it is easy to check
\begin{equation*}
		\begin{split}
		\X^{t+1} &= \X^{t+1}_U \left(\B^{t+1}\right)^{-1} \X_V^{t+1\top} = [\U^t \quad \U^{t}_\perp]\begin{bmatrix}
\B^{t+1}  & \D_2^{t+1\top} \\
\D_1^{t+1} & \D_1^{t+1} (\B^{t+1})^{-1} \D_2^{t+1\top}
\end{bmatrix} [\V^t \quad \V^t_\perp]^{\top}.
	\end{split}
\end{equation*}  So
\begin{equation} \label{eq: X^t+1 - barX bound}
	\|\X^{t+1} - \widebar{\X}\|_\F^2 = \left\| \begin{array}{c c}
	\B^{t+1} - \widetilde{\B}^t & \D_2^{t+1 \top} - \widetilde{\D}_2^{t\top}\\
	\D_1^{t+1} - \widetilde{\D}_1^t & \bDelta^{t+1} 
	\end{array} \right\|_\F^2,
\end{equation}
where $\bDelta^{t+1} = \D_1^{t+1} (\B^{t+1})^{-1} \D_2^{t+1\top}- \widetilde{\D}_1^{t} (\widetilde{\B}^{t})^{-1} \widetilde{\D}_2^{t\top}$.
Recall that \eqref{eq: B and D part bound} gives a precise error characterization for $\|\B^{t+1} - \widetilde{\B}^t \|_\F^2 + \sum_{k=1}^2 \|\D_k^{t+1} - \widetilde{\D}_k^t\|_\F^2$. Hence, to bound $\|\X^{t+1} - \widebar{\X}\|_\F^2$ from \eqref{eq: X^t+1 - barX bound}, we only need to obtain an upper bound of $\|\bDelta^{t+1}\|_\F^2$. 

Combining \eqref{def: rho t+1}, \eqref{ineq: bound for D1 B(-1)} and Lemma \ref{lm:FGH}, we have
\begin{equation*}
\begin{split}
	\|\bDelta^{t+1}\|_\F 
	\leq \rho_{t+1} \| \D_1^{t+1} - \widetilde{\D}_1^t \|_\F + \frac{\theta_t}{\sqrt{1 - \theta_t^2}} \|\D_2^{t+1} - \widetilde{\D}_2^t\|_\F + \frac{\theta_t \rho_{t+1}}{\sqrt{1 - \theta_t^2}} \|\B^{t+1} - \widetilde{\B}^t\|_\F. 
\end{split}
\end{equation*} Moreover, since $(a + b + c)^2 \leq 3(a^2 + b^2 + c^2)$ and \eqref{eq: B and D part bound}, it holds that 
\begin{equation} \label{ineq: D1 B^{-1} D2^top bound}
	\|\bDelta^{t+1}\|_\F^2 \leq 3(\rho_{t+1} \vee \frac{\theta_t}{\sqrt{1 - \theta_t^2}} \vee \frac{\theta_t \rho_{t+1}}{\sqrt{1 - \theta_t^2}} )^2 \left\|(\cL^*_t \cA^* \cA \cL_t)^{-1} \cL_t^* \cA^* \bepsilon^t\right\|_\F^2.
\end{equation}
In summary, combining \eqref{eq: X^t+1 - barX bound}, \eqref{eq: B and D part bound} and \eqref{ineq: D1 B^{-1} D2^top bound}, we get
\begin{equation} 
	\|\X^{t+1} - \widebar{\X}\|_\F^2 \leq (1 + 3\big(\rho_{t+1} \vee \frac{\theta_t}{\sqrt{1 - \theta_t^2}} \vee \frac{\theta_t \rho_{t+1}}{\sqrt{1 - \theta_t^2}} )^2\big) \left\|(\cL^*_t \cA^* \cA \cL_t)^{-1} \cL_t^* \cA^* \bepsilon^t\right\|_\F^2.
\end{equation}
Plugging in the upper bounds for $\theta_t$ and $\rho_{t+1}$, we have
$ 1 + 3(\rho_{t+1} \vee \frac{\theta_t}{\sqrt{1 - \theta_t^2}} \vee \frac{\theta_t \rho_{t+1}}{\sqrt{1 - \theta_t^2}} )^2 \leq 5$, and thus
\begin{equation*}
\begin{split}
	\|\X^{t+1} - \widebar{\X}\|_\F^2 &\leq 5 \left\|(\cL^*_t \cA^* \cA \cL_t)^{-1} \cL_t^* \cA^* \bepsilon^t\right\|_\F^2\\
	& \overset{(a)}\leq \frac{5\|\X^t - \widebar{\X}\|^2}{(1-R_{2r})^2 \sigma^2_r(\widebar{\X}) } \left( R_{3r}^2 \|\X^t - \widebar{\X}\|_\F^2 + 4 \|\cA^*(\widebar{\bepsilon})\|_\F^2 + 4 R_{3r}\|\cA^*(\widebar{\bepsilon})\|_\F \|\X^t - \widebar{\X}\|_\F \right) \\
	& \overset{(b)}\leq \frac{9}{16} \|\X^t - \widebar{\X} \|_\F^2,
\end{split}
\end{equation*} 
where (a) is due to \eqref{ineq: least square residual upper bound} and (b) is because under the initialization condition \eqref{ineq: optimization theory initialization condition} and $\|\cA^*(\widebar{\bepsilon})\|_\F \leq \frac{1-R_{2r}}{4\sqrt{5}}\sigma_r(\widebar{\X})$, the following inequalities hold
\begin{equation*}
	\begin{split}
		\frac{5 R_{3r}^2\|\X^t - \widebar{\X}\|_\F^2}{(1-R_{2r})^2 \sigma^2_r(\widebar{\X}) } \leq 1/16,\quad 
		\frac{20 \|\cA^*(\widebar{\bepsilon})\|_\F^2}{(1-R_{2r})^2 \sigma^2_r(\widebar{\X}) }  \leq 1/4, \quad \frac{20 R_{3r}\|\X^t - \widebar{\X}\|_\F \|\cA^*(\widebar{\bepsilon})\|_\F  }{(1-R_{2r})^2 \sigma^2_r(\widebar{\X}) } \leq 1/4.
	\end{split}
\end{equation*} This completes the induction and finishes the proof of this theorem. \quad $\blacksquare$

{\noindent \bf Proof of Theorem \ref{th: Riemannian Gauss-Newton of Simul Mini}}.  
In view of the Riemannian Gauss-Newton equation in \eqref{eq: Riemannian Gauss-newton equation}, to prove the claim, we only need to show
\begin{equation} \label{eq: Thm Rie G-N need to proof}
	P_{T_{\X^t}}(\cA^*(\cA(\eta^t + \X^t) -\y)) = 0.
\end{equation} 
From the optimality condition of the least squares problem \eqref{eq: compact repre of least square}, we know that
$\B^{t+1}$, $\D_1^{t+1}$ and $\D_2^{t+1}$ obtained in \eqref{eq: alg2 least square} satisfy
\begin{equation}\label{eq:LkZk}
\begin{split}
	&{\cal L}_t^*\cA^*(\cA {\cal L}_t
\begin{bmatrix}
\B^{t+1} & (\D_2^{t+1})^\top \\[2pt]
\D_1^{t+1} & \0
\end{bmatrix}
 - \y) = 0.
\end{split}
\end{equation}
Then, the updating formula in \eqref{eq: update before retraction}, together with the definition of ${\cal L}_t$, implies that
\begin{equation}\label{eq: eta^t + X^t}
	\eta^t + \X^t = {\cal L}_t
\begin{bmatrix}
\B^{t+1} & (\D_2^{t+1})^\top \\[2pt]
\D_1^{t+1} & \0
\end{bmatrix}.
\end{equation}
Hence, \eqref{eq:LkZk} implies ${\cal L}_t^*{\cal A}^*({\cal A} (\eta^t + \X^t) - \y) = 0$. 
Note from the proof of Theorem \ref{th: local contraction general setting} that for all $t\ge 1$, $\B^{t}$ is invertible. Then, it is not difficult to verify that $\U^t,\V^t$ are orthonormal bases of the column and row spans of $\X^t$ and  $P_{T_{\X^t}} = {\cal L}_t {\cal L}_t^*$ for all $t\ge 0$. We thus proved \eqref{eq: Thm Rie G-N need to proof}. \quad $\blacksquare$

{\noindent \bf Proof of Proposition \ref{prop: descent direction}}. We compute the inner product between the update direction $\eta^t$ in \eqref{eq: update before retraction} and the Riemannian gradient:
\begin{equation*}
	\begin{split}
		\langle \grad f(\X^t) , \eta^t \rangle &= \langle P_{T_{\X^t}} \cA^* (\cA(\X^t) - \y) ,\eta^t \rangle\\
		& \overset{(a)}= \langle - P_{T_{\X^t}} \cA^*\cA \eta^t , \eta^t \rangle \\
		& \overset{(b)}= - \langle \cA^*\cA \eta^t , \eta^t\rangle = - \|\cA(\eta^t)\|_2^2,
	\end{split}
\end{equation*} here (a) is due to \eqref{eq: Thm Rie G-N need to proof} and (b) is because $\eta^t$ lies in ${T_{\X^t} \cM_r}$. With this, we conclude the update direction $\eta^t$ has a negative inner product with the Riemannian gradient unless it is $0$. Thus the update $\eta^t$ is a descent direction.

If $\cA$ satisfies the $2r$-RIP, by similar arguments as in Lemma \ref{lm: spectral norm bound of Atop A}, we see that $P_{T_{\X^t}}\cA^*\cA P_{T_{\X^t}}$ is symmetric positive definite over ${T_{\X^t}\cM_r}$ for all $t \ge 0$. 
Since $\eta^t$ solves the Riemannian Gauss-Newton equation \eqref{eq: Riemannian Gauss-newton equation}, we know that $\eta^t = - (P_{T_{{\X^t}}} \cA^*\cA P_{T_{{\X^t}}})^{-1} \grad f({\X^t})$ for all $t\ge 0$.
For any subsequence $\{\X^t\}_{t\in \cal K}$ that converges to a nonstationary point $\widetilde{\X}$, it is not difficult to show that
\[
\lim_{t \to \infty, t\in \cal K} P_{T_{\X^t}} \cA^*\cA P_{T_{\X^t}} = P_{T_{\widetilde{\X}}} \cA^*\cA P_{T_{\widetilde{\X}}} \quad \mbox{and} \quad \tilde \eta = \lim_{t \to \infty, t\in \cal K} \eta^t = - (P_{T_{\widetilde{\X}}} \cA^*\cA P_{T_{\widetilde{\X}}})^{-1} \grad f(\widetilde{\X}).
\]
Hence, $\{\eta^t\}_{t\in \cal K}$ is bounded, $\tilde \eta \neq 0$ and
\[
	\lim_{t\to \infty, t\in \cal K} \inprod{\grad f(\X^t)}{\eta^t} = \inprod{\grad f(\widetilde{\X})}{\tilde \eta} = -\|\cA(\tilde \eta)\|_2^2 \overset{\text{(RIP  condition)}} \leq  -  (1 - R_{2r}) \|\tilde \eta\|_2^2 <0,
\]
i.e., the direction sequence $\{\eta^t\}$ is gradient related by \cite[Definition 4.2.1]{absil2009optimization}.
\quad $\blacksquare$

{\noindent \bf Proof of Theorem \ref{th: local convergence in local rank matrix recovery}}. The proof of \eqref{ineq: low rank matrix recovry quadratic-linear converge} shares many similar ideas to the proof of Theorem \ref{th: local contraction general setting}. Hence, we point out the main difference first and then give the complete proof. Compared to Theorem \ref{th: local contraction general setting} where the target matrix is a stationary point, here the target matrix is $\X^*_{\max(r)}$. So when we apply Proposition \ref{prop: bound for iter approx error}, $\widebar{\X} = \X^*_{\max(r)}$, $\widebar{\bepsilon} = \tilde{\bepsilon}$ and we no longer have \eqref{ineq: bound for Thm 2, term 2, first inequality}. Due to this difference, here we have an unavoidable statistical error term in the upper bound. 

We begin by proving \eqref{ineq: low rank matrix recovry quadratic-linear converge}. We first apply Proposition \ref{prop: bound for iter approx error} to bound $\left\|(\cL^*_t \cA^* \cA \cL_t)^{-1} \cL_t^* \cA^* \bepsilon^t\right\|_\F^2$ in this setting. Set $\widebar{\X} = \X^*_{\max(r)}$ and $\widebar{\bepsilon} = \tilde{\bepsilon}$, by Proposition \ref{prop: bound for iter approx error}, we have
\begin{equation} \label{ineq: least square residual bound}
	\begin{split}
		&\left\|(\cL^*_t \cA^* \cA \cL_t)^{-1} \cL_t^* \cA^* \bepsilon^t\right\|_\F^2\\
		 \leq & \frac{R_{3r}^2 \|\X^t - \X^* \|^2 \|\X^t - \X^* \|_\F^2 }{(1-R_{2r})^2\sigma_r^2(\X^*)} + \frac{\| \cL_t^* \cA^* (\tilde{\bepsilon})\|_\F^2}{(1-R_{2r})^2} + \| \cL_t^* \cA^* (\tilde{\bepsilon})\|_\F \frac{2R_{3r} \|\X^t - \X^* \| \|\X^t - \X^* \|_\F }{\sigma_r(\X^*)(1-R_{2r})^2}\\
		 \overset{(a)}\leq & 2\frac{R_{3r}^2 \|\X^t - \X^* \|^2 \|\X^t - \X^* \|_\F^2 }{(1-R_{2r})^2\sigma_r^2(\X^*)} + 2\frac{\| \cL_t^* \cA^* (\tilde{\bepsilon})\|_\F^2}{(1-R_{2r})^2},
	\end{split}
\end{equation} where (a) is by Cauchy-Schwarz inequality.
Recall $P_{\X^t}$ in \eqref{eq: prop projection onto tangent space}, $P_{\X^t} (\cA^*(\tilde{\bepsilon}))$ is a at most rank $2r$ matrix and $\sigma_i(P_{\X^t} \cA^*(\tilde{\bepsilon})) \leq \sigma_i(\cA^*(\tilde{\bepsilon}))$ for $1 \leq i \leq p_1 \wedge p_2$ by the projection property of $P_{\X^t}$. Then we have
\begin{equation} \label{ineq: noise bound truncated to max r}
	\|\cL_t^* \cA^*(\tilde{\bepsilon})\|^2_\F = \|P_{\X^t} \cA^*(\tilde{\bepsilon})\|^2_\F \leq \| (\cA^*(\tilde{\bepsilon}))_{\max(2r)} \|^2_\F \leq 2 \| (\cA^*(\tilde{\bepsilon}))_{\max(r)} \|^2_\F.
\end{equation}

Recall in this setting, the target matrix is $\X^*_{\max(r)} $. We replace $\widebar \X$ in \eqref{eq: tildeB,D} by $\X^*_{\max(r)} $ and obtain $\widetilde{\B}^t := \U^{t\top} \X^*_{\max(r)}  \V^t, \widetilde{\D}_1^t := \U^{t\top}_\perp\X^*_{\max(r)}  \V^t, \widetilde{\D}_2^{t\top} := \U^{t\top} \X^*_{\max(r)}  \V^{t}_\perp$. Next, we use induction to prove the main results. Define $\theta_t, \rho_{t+1}$ in the same way as in the proof of Theorem \ref{th: local contraction general setting}. We aim to show: given $\theta_t \leq 1/2$, $\widetilde{\B}^t$ is invertible, $\|\X^t -\X^*_{\max(r)} \|_\F \leq \|\X^0 - \X^*_{\max(r)} \|_\F \vee \frac{2\sqrt{10}}{(1-R_{2r})} \| (\cA^*(\tilde{\bepsilon}))_{\max(r)} \|_\F$, then $\theta_{t+1} \leq 1/2$, $\widetilde{\B}^{t+1}$ is invertible, $\|\X^{t+1} -\X^*_{\max(r)} \|_\F \leq \|\X^0 - \X^*_{\max(r)} \|_\F \vee \frac{2\sqrt{10}}{(1-R_{2r})} \| (\cA^*(\tilde{\bepsilon}))_{\max(r)} \|_\F$, as well as $\B^{t+1}$ is invertible and \eqref{ineq: low rank matrix recovry quadratic-linear converge}.

First, we can easily check the assumption is true when $t = 0$ under the initialization condition. Now, suppose the induction assumption is true at iteration $t$. Under the conditions \eqref{ineq: low rank matrix recovery initialization condition}, \eqref{ineq: LRMR sigma condition}, from \eqref{ineq: least square residual bound} and \eqref{ineq: noise bound truncated to max r}, we have
\begin{equation*}
	\left\|(\cL^*_t \cA^* \cA \cL_t)^{-1} \cL_t^* \cA^* \bepsilon^t\right\|_\F \leq \frac{3}{16} \sigma_r(\X^*).
\end{equation*} Then following the same proof as the Step 3,4 of Theorem \ref{th: local contraction general setting}, we have $\B^{t+1}$ is invertible, $\rho_{t+1} \leq (4 + \sqrt{3})/3\sqrt{3}$ and
\begin{equation} \label{ineq: Xt+1 - barX bound in low rank matrix recovery}
\begin{split}
	\|\X^{t+1} - \X^*_{\max(r)} \|_\F^2 &\leq (1 + 3(\rho_{t+1} \vee \frac{\theta_t}{\sqrt{1 - \theta_t^2}} \vee \frac{\theta_t \rho_{t+1}}{\sqrt{1 - \theta_t^2}} )^2) \left\|(\cL^*_t \cA^* \cA \cL_t)^{-1} \cL_t^* \cA^* \bepsilon^t\right\|_\F^2\\
	& \leq 5 \left\|(\cL^*_t \cA^* \cA \cL_t)^{-1} \cL_t^* \cA^* \bepsilon^t\right\|_\F^2.
\end{split}
\end{equation}
Plugging \eqref{ineq: noise bound truncated to max r} and \eqref{ineq: least square residual bound} into \eqref{ineq: Xt+1 - barX bound in low rank matrix recovery}, we arrive at 
	\begin{equation} \label{ineq: proof low rank matrix recovry quadratic-linear converge}
	\begin{split}
		\|\X^{t+1} - \X^*_{\max(r)} \|_\F^2  &\leq 10\frac{R_{3r}^2 \|\X^t - \X^*_{\max(r)}  \|^2 \|\X^t - \X^*_{\max(r)}  \|_\F^2 }{(1-R_{2r})^2\sigma_r^2(\X^*)} + \frac{20 \| ( \cA^*(\tilde{\bepsilon} ) )_{\max(r)} \|_\F^2}{(1-R_{2r})^2}\\
		& \overset{(a)}\leq \frac{1}{2}\|\X^{t} - \X^*_{\max(r)} \|_\F^2 + \frac{20}{(1-R_{2r})^2} \| ( \cA^*(\tilde{\bepsilon} ) )_{\max(r)} \|_\F^2,
	\end{split}
	\end{equation} where (a) is because under conditions \eqref{ineq: low rank matrix recovery initialization condition}, \eqref{ineq: LRMR sigma condition} and induction assumption at iteration $t$, it holds that
	\begin{equation*}
		\frac{10 R_{3r}^2\|\X^t - \X^*_{\max(r)} \|_\F^2}{(1-R_{2r})^2 \sigma^2_r(\X^*)} \leq 1/2.
	\end{equation*}
	By \eqref{ineq: proof low rank matrix recovry quadratic-linear converge}, we get $\|\X^{t+1} -\X^*_{\max(r)} \|_\F \leq \|\X^0 - \X^*_{\max(r)} \|_\F \vee \frac{2\sqrt{10}}{(1-R_{2r})} \| (\cA^*(\tilde{\bepsilon}))_{\max(r)} \|_\F$. Under the initialization conditions, Lemma \ref{lm: perturbation bound} also implies $\theta_{t+1} \leq 1/2$ and $\widetilde{\B}^{t+1}$ is invertible. This finishes the proof of \eqref{ineq: low rank matrix recovry quadratic-linear converge}.

Next, we prove the guarantee of RISRO under the sub-Gaussian ensemble design with spectral initialization. Throughout the proof, we use various $c, C, C_1$ to denote constants and they may vary from line to line. Recall now $\X^*$ is a rank $r$ matrix, so we have $\tilde{\bepsilon} = \bepsilon$. First we give the guarantee for the initialization $\X^0 = (\cA^*(\y))_{\max(r)}$. Define $\Q_0 \in \bbR^{p_1 \times 2r}$ as an orthogonal matrix which spans the column subspaces of $\X^0$ and $\X^*$. Let $\Q_{0\perp}$ be the orthogonal complement of $\Q_0$. Since
\begin{equation*}
	\|\X^0 - \cA^*(\y)\|_\F^2 = \|\X^0 - P_{\Q_0} (\cA^*(\y)) \|_\F^2 + \| P_{\Q_{0\perp}} (\cA^*(\y)) \|_\F^2
\end{equation*} and 
\begin{equation*}
	\|\X^* - \cA^*(\y)\|_\F^2 = \|\X^* - P_{\Q_0} (\cA^*(\y)) \|_\F^2 + \| P_{\Q_{0\perp}} (\cA^*(\y)) \|_\F^2,
\end{equation*} the SVD property $\|\X^0 - \cA^*(\y)\|_\F^2 \leq \|\X^* - \cA^*(\y)\|_\F^2$ implies that
\begin{equation*}
	\|\X^0 - P_{\Q_0} (\cA^*(\y)) \|_\F^2 \leq \|\X^* - P_{\Q_0} (\cA^*(\y)) \|_\F^2.
\end{equation*} 
Note that
\begin{equation} \label{ineq: spec bnd in init}
	\begin{split}
		\|P_{\Q_0} - P_{\Q_0} \cA^*\cA P_{\Q_0}\| & \overset{(a)}= \sup_{\|\Z\|_\F \leq 1} \langle(P_{\Q_0} - P_{\Q_0} \cA^*\cA P_{\Q_0})\Z, \Z \rangle \\
		& = \sup_{\|\Z\|_\F \leq 1} \left| \|P_{\Q_0}\Z\|_\F^2 - \|\cA P_{\Q_0}\Z\|_\F^2  \right|\\
		& \overset{(b)}\leq \sup_{\|\Z\|_\F \leq 1} R_{2r} \|P_{\Q_0}\Z\|_\F^2 \leq R_{2r},
	\end{split}
\end{equation} 
where (a) is because $P_{\Q_0} - P_{\Q_0} \cA^*\cA P_{\Q_0}$ is  symmetric and (b) is by the  $2r$-RIP of $\cA$. Hence, 
\begin{equation} \label{ineq: bound for initialization}
\begin{split}
	\|\X^0 - \X^*\|_\F &\leq \|\X^0 - P_{\Q_0} (\cA^*(\y)) \|_\F + \|\X^* - P_{\Q_0} (\cA^*(\y)) \|_\F\\
	& \leq 2\|\X^* - P_{\Q_0} (\cA^*(\y)) \|_\F\\
	& \overset{(a)}= 2 \| \X^* - P_{\Q_0} (\cA^*(\cA(\X^*) + \bepsilon)) \|_\F\\
	& = 2 \| P_{\Q_0}\X^* - P_{\Q_0} \cA^*\cA(P_{\Q_0}\X^*) -  P_{\Q_0} (\cA^*(\bepsilon)) \|_\F\\
	& \leq 2 \left( \| (P_{\Q_0} - P_{\Q_0} \cA^*\cA P_{\Q_0})\X^*\|_\F + \|  P_{\Q_0} (\cA^*(\bepsilon)) \|_\F \right)\\
	& \overset{(b)}\leq 2 R_{2r} \|\X^*\|_\F + 2\sqrt{2} \|(\cA^*(\bepsilon))_{\max(r)}\|_\F\\
	& \leq 2R_{2r} \sqrt{r} \kappa \sigma_r(\X^*) + 2\sqrt{2} \|(\cA^*(\bepsilon))_{\max(r)}\|_\F,
\end{split}
\end{equation} where (a) is due to the model of $\y$ and (b) is due to that $P_{\Q_0} (\cA^*(\bepsilon))$ is a at most rank 2r matrix and the spectral norm bound for the operator $(P_{\Q_0} - P_{\Q_0} \cA^*\cA P_{\Q_0})$ in \eqref{ineq: spec bnd in init}.
Hence, there exists $c_1, c_2, C > 0$ such that when 
\begin{equation} \label{ineq: RIP and least singular value requirement}
	R_{2r} \leq c_1 \frac{1}{\kappa \sqrt{r} }, \, R_{3r} < \frac{1}{2}, \text{ and } \sigma_r(\X^*) \geq C \| (\cA^*(\bepsilon))_{\max(r)} \|_\F,
\end{equation} we have $\|\X^0 - \X^*\|_\F \leq c_2 \sigma_r(\X^*)$ by \eqref{ineq: bound for initialization} and the conditions in \eqref{ineq: low rank matrix recovery initialization condition} and \eqref{ineq: LRMR sigma condition} are satisfied.

Next, we show under the sample complexity indicated in the Theorem, \eqref{ineq: RIP and least singular value requirement} are satisfied with high probability. First by a similar argument of \cite[Lemma 6]{zhang2020islet}, for the sub-Gaussian ensemble design considered here, we have with probability at least $1 - \exp(-c(p_1 + p_2))$ for some $c > 0$ that $\|\cA^*(\bepsilon)\| \leq c'\sqrt{\frac{p_1 + p_2}{n}} \sigma$. So with the same high probability, we have 
\begin{equation} \label{ineq: high prob bound for A^*(epsilon)}
\| (\cA^*(\bepsilon))_{\max(r)} \|_\F \leq c'\sqrt{\frac{(p_1 + p_2)r}{n}} \sigma,
\end{equation}
and when $n \geq C (p_1 + p_2) r \frac{\sigma^2}{\sigma^2_r(\X^*)} $, we have $\sigma_r(\X^*) \geq C \| (\cA^*(\bepsilon))_{\max(r)} \|_\F$. At the same time, by \cite[Theorem 2.3]{candes2011tight}, there exists $C > 0$ when $n \geq C (p_1 + p_2) \kappa^2 r^2$, $R_{2r} \leq c_1 \frac{1}{\kappa \sqrt{r} }$ and $R_{3r} < \frac{1}{2}$ are satisfied with probability at least $1 - \exp(- c(p_1 + p_2))$ for some $c > 0$.

In summary, there exists $C > 0$ such that when $n \geq C (p_1 + p_2) r (\frac{\sigma^2}{\sigma_r^2(\X^*)} \vee r\kappa^2 )$,  \eqref{ineq: RIP and least singular value requirement} holds with probability at least $1 - \exp(-c(p_1 + p_2))$ for some $c > 0$. So by the first part of the Theorem, we have with the same high probability:
	\begin{equation*}
		\|\X^{t+1} - \X^*\|_\F^2 \leq 10\frac{R_{3r}^2 \|\X^t - \X^* \|_\F^4 }{(1-R_{2r})^2\sigma_r^2(\X^*)} + \frac{20 \| ( \cA^*(\bepsilon ) )_{\max(r)} \|_\F^2}{(1-R_{2r})^2}, \quad \forall\, t\ge 0.
	\end{equation*}
	
More specifically, the above convergence can be divided into two phases. Let
\begin{itemize}
	\item (Phase I) When $\|\X^{t} - \X^*\|_\F^2 \geq \frac{\sqrt{2}}{R_{3r}} \| ( \cA^*(\bepsilon ) )_{\max(r)} \|_\F  \sigma_r(\X^*)$,
	\begin{equation*}
		\|\X^{t+1} - \X^*\|_\F \leq 2\sqrt{5}\frac{R_{3r} \|\X^t - \X^* \|_\F^2 }{(1-R_{2r})\sigma_r(\X^*)}
	\end{equation*}
	\item (Phase II) When $\|\X^{t} - \X^*\|_\F^2 \leq \frac{\sqrt{2}}{R_{3r}} \| ( \cA^*(\bepsilon ) )_{\max(r)} \|_\F  \sigma_r(\X^*)$,
	\begin{equation*}
		\|\X^{t+1} - \X^*\|_\F \leq  \frac{2\sqrt{10} \| ( \cA^*(\bepsilon ) )_{\max(r)} \|_\F}{1-R_{2r}}
	\end{equation*}
\end{itemize}
	
Combining Phase I, II and \eqref{ineq: high prob bound for A^*(epsilon)}, by induction we have $\|\X^t - \X^*\|_\F \leq 2^{-2^t} \|\X^0 - \X^*\|_\F + c\sqrt{\frac{r(p_1 + p_2)\sigma^2}{n}}$ and this implies the desired error bound for $\|\X^t -\X^*\|_\F$ after double-logarithmic number of iterations.
\quad $\blacksquare$

{\noindent \bf Proof of Theorem \ref{thm: quadratic convergence for phase retrieval}}.
In the phase retrieval example, the mapping $\cA$ no longer satisfies a proper RIP condition and the strategy we use is to show the contraction of $\X^t - \X^*$ in terms of its nuclear norm and then transform it back to Frobenius norm. 

We also use induction to show the main results. Specifically, we show: given $|\u^{*\top}\u^t| > 0$ where $\u^* = \frac{\x^*}{\|\x^*\|_2}$ and $\|\X^t - \X^*\|_\F \leq \|\X^0 - \X^*\|_\F$, then $|\u^{*\top}\u^{t+1}| > 0$, $\|\X^{t+1} - \X^*\|_\F \leq \|\X^0 - \X^*\|_\F$ and \eqref{ineq: quadra conver in PR}. 

First, the induction assumption is true when $t=0$ by the initialization condition and the perturbation bound in Lemma \ref{lm: perturbation bound}. Assume it is also correct at iteration $t$. Let $\widetilde{b}^t = \u^{t\top}\X^*\u^t, \widetilde{\d}^t = (\u_\perp)^{t\top} \X^*\u^t$. It is easy to verify $\widetilde{\d}^t(\widetilde{b}^t)^{-1} \widetilde{\d}^t = \u_\perp^{t\top} \X^* \u_\perp^t$ and $\X^* = [\u^t \, \u_\perp^t] \left[ \begin{array}{c c}
	\widetilde{b}^t & \widetilde{\d}^{t\top}\\
	\widetilde{\d}^t & \widetilde{\d}^t(\widetilde{b}^t)^{-1} \widetilde{\d}^t
	\end{array}
	  \right] [\u^t \, \u^t_{\perp} ]^\top.$
Define the linear operator $\cL_t$ similar as \eqref{eq: linear operator L_t} in this setting in the following way	  
 \begin{equation} \label{eq: Lt in PR}
 	\cL_t: \W= \left[ \begin{array}{c c}
	w_0 \in \bbR & \w_1^\top \in \bbR^{1 \times (p-1)}\\
	\w_1 \in \bbR^{(p-1) \times 1} & \0
\end{array}
 \right] \to [\u^t \, \u_\perp^t] \left[ \begin{array}{c c}
	w_0 & \w_1^{\top}\\
	\w_1 & \0	
	\end{array}
	  \right] [\u^t \, \u^t_{\perp} ]^\top,
 \end{equation} 
and it is easy to compute its adjoint $\cL_t^*(\M) = \left[ \begin{array}{c c}
	\u^{t\top} \M \u^t & \u^{t\top} \M \u_\perp^t \\
	\u_{\perp}^{t\top} \M \u^t & 0
	\end{array}
	  \right]$, where $\M$ is a rank 2 symmetric matrix. Define operator $P_{\X^t}$ similar as \eqref{eq: prop projection onto tangent space} over the space of $p\times p$ symmetric matrices $$P_{\X^t}(\W) := \cL_t \cL_t^* (\W) = \u^t \u^{t\top} \W \u^t \u^{t\top} + \u_\perp^t \u_\perp^{t\top} \W \u^t \u^{t\top} + \u^t \u^{t\top} \W \u_\perp^t \u_\perp^{t\top}.$$ It is easy to verify that $P_{\X^t}$ is an orthogonal projector. Meanwhile, let $P_{(\X^t)_\perp}(\W) = \W - P_{\X^t}(\W)= \u_\perp^t \u_\perp^{t\top} \W \u_\perp^t \u_\perp^{t\top}$. 

By using the operator $\cL_t$, the least squares solution in Step \ref{alg: PR least square} can be rewritten in the following way
\begin{equation}\label{eq: least square solution in PR}
\begin{split}
		\left[ \begin{array}{c c}
		b^{t+1} & \d^{t+1\top}\\
		\d^{t+1} & \0
	\end{array}  \right] &= \argmin_{b \in \bbR, \d \in \bbR^{p-1}} \left\|\y - \cA \cL_t \left[ \begin{array}{c c}
		b & \d^\top\\
		\d & \0
	\end{array}  \right] \right\|_2^2\\
	&= (\cL_t^* \cA^* \cA \cL_t)^{-1} \cL_t^* \cA^* \y\\
	&  = (\cL_t^* \cA^* \cA \cL_t)^{-1} \cL_t^* \cA^* (\cA(P_{\X^t}\X^*  ) + \cA(P_{(\X^t)_\perp} (\X^*) ))\\
	& = \left[ \begin{array}{c c}
	\widetilde{b}^t & \widetilde{\d}^{t\top}\\
	\widetilde{\d}^t & \0
	\end{array}
	  \right] +  (\cL_t^* \cA^* \cA \cL_t)^{-1} \cL_t^* \cA^*  \cA(P_{(\X^t)_\perp} (\X^*) ).
\end{split}
\end{equation} 
Here, $\cL_t^* \cA^* \cA \cL_t$ is invertible is due to the lower bound of the spectrum of $\cL_t^* \cA^* \cA \cL_t$ in Lemma \ref{lm: nuclear norm bound of AtopA in phase retrieval}.
So 
\begin{equation} \label{ineq: PR iteration error bound}
\begin{split}
	\|\X^{t+1} - \X^*\|_* \leq & \left\| \X^{t+1} - [\u^t\, \u^t_{\perp}] \left[ \begin{array}{c c}
		b^{t+1} & \d^{t+1\top}\\
		\d^{t+1} & \0
	\end{array}  \right]  [\u^t\, \u^t_{\perp}]^\top \right\|_* \\
	    & +  \left\| [\u^t\, \u^t_{\perp}] \left[ \begin{array}{c c}
		b^{t+1} & \d^{t+1\top}\\
		\d^{t+1} & \0
	\end{array}  \right]  [\u^t\, \u^t_{\perp}]^\top - (\u^t\, \u^t_{\perp}) \left[ \begin{array}{c c}
		\widetilde{b}^t & \widetilde{\d}^{t\top}\\
		\widetilde{\d}^{t} & \widetilde{\d}^{t} (\widetilde{b}^t )^{-1}\widetilde{\d}^{t\top}
	\end{array}  \right]  [\u^t\, \u^t_{\perp}]^\top \right\|_*\\
	\overset{(a)} \leq & 2 \left\| [\u^t\, \u^t_{\perp}] \left[ \begin{array}{c c}
		b^{t+1} & \d^{t+1\top}\\
		\d^{t+1} & \0
	\end{array}  \right]  [\u^t\, \u^t_{\perp}]^\top - [\u^t\, \u^t_{\perp}] \left[ \begin{array}{c c}
		\widetilde{b}^t & \widetilde{\d}^{t\top}\\
		\widetilde{\d}^{t} & \widetilde{\d}^{t} (\widetilde{b}^t )^{-1}\widetilde{\d}^{t\top}
	\end{array}  \right]  [\u^t\, \u^t_{\perp}]^\top \right\|_*\\
	\leq & 2 \left\| \begin{array}{c c}
		b^{t+1} - \widetilde{b}^t & \d^{t+1\top} - \widetilde{\d}^{t\top} \\
		\d^{t+1} - \widetilde{\d}^t & \0
	\end{array}  \right\|_* + 2 \| \widetilde{\d}^{t} (\widetilde{b}^t )^{-1}\widetilde{\d}^{t\top} \|_*\\
	\overset{(b)} = & 2 \|(\cL_t^* \cA^* \cA \cL_t)^{-1} \cL_t^* \cA^*  \cA(P_{(\X^t)_\perp} \X^* )\|_* + 2 \| \widetilde{\d}^{t} (\widetilde{b}^t )^{-1}\widetilde{\d}^{t\top} \|_*,
\end{split}
\end{equation} here (a) is due to Lemma \ref{lm: PSD projection in nuclear norm} and (b) is due to \eqref{eq: least square solution in PR}.

First notice $\| \widetilde{\d}^{t} (\widetilde{b}^t )^{-1}\widetilde{\d}^{t\top} \|_* = \|P_{(\X^{t})_\perp} \X^* \|_*$. Next we give bound for $\|(\cL_t^* \cA^* \cA \cL_t)^{-1} \cL_t^* \cA^*  \cA(P_{(\X^t)_\perp} \X^* )\|_*$.  With probability at least $1 - C_1 \exp(-C_2(\delta_1, \delta_2) p) - C_3n^{-p}$ ($C_1, C_2, C_3 > 0$), we have
\begin{equation*}
	\begin{split}
		&\|(\cL_t^* \cA^* \cA \cL_t)^{-1} \cL_t^* \cA^*  \cA(P_{(\X^t)_\perp} \X^* )\|_*\\
		\overset{(a)}\leq & \frac{4}{(1-\delta_1)n} \|\cL_t^*\cA^* \cA P_{(\X^t)_\perp}(\X^*)  \|_*\\
		\overset{(b)} \leq & C \frac{4p}{(1-\delta_1)n} \|\cA P_{(\X^t)_\perp}(\X^*)\|_1\\
		\overset{(c)} \leq & C' \frac{1+\delta_2}{1-\delta_1} p \| P_{(\X^t)_\perp}(\X^*)\|_*
	\end{split}
\end{equation*} 
for some $ C,C' > 0$. Here $\|\cdot\|_1$ denotes the $\ell_1$ norm of a vector, (a) is due to Lemma \ref{lm: nuclear norm bound of AtopA in phase retrieval}, (b) is due to Lemma \ref{lm: nuclear norm bound for LtA} and (c) is due to \cite[Lemma 3.1]{candes2013phaselift} and $P_{(\X^t)_\perp}(\X^*)$ is a symmetric matrix.

Putting above results into \eqref{ineq: PR iteration error bound}, we have with probability at least $1 - C_1 \exp(-C_2(\delta_1,\delta_2) p) - C_3n^{-p}$,
\begin{equation*}
\begin{split}
 \|\X^{t+1} - \X^*\|_\F \leq & \|\X^{t+1} - \X^*\|_* \leq C\frac{1+\delta_2}{1-\delta_1}p \|P_{(\X^t)_\perp}(\X^*)\|_* \overset{(a)}= C\frac{1+\delta_2}{1-\delta_1}p \|P_{(\X^t)_\perp}(\X^*)\|_\F\\
 \overset{(b)}\leq & C\frac{1+\delta_2}{1-\delta_1}p \frac{\|\X^t - \X^*\|_\F^2}{\sigma_1(\X^*)},
\end{split}
\end{equation*} 
where (a) is because $P_{(\X^t)_\perp}(\X^*)$ is a symmetric rank $1$ matrix, (b) is due to the same argument as \eqref{ineq: bound for Thm 3 term (1)}. Since $\sigma_1(\X^*) = \|\X^*\|_\F$, when $\|\X^0 - \X^*\|_\F \leq \frac{(1-\delta_1)}{C (1+\delta_2) p} \|\X^*\|_\F $ for some large enough $C$, we have $\|\X^{t+1} - \X^*\|_\F \leq \|\X^t - \X^*\|_\F \leq \|\X^0 - \X^*\|_\F$. Also $|\u^{t+1\top}\u^*| > 0$ as it is a non-decreasing function of $\|\X^{t+1} - \X^*\|_\F$ by Lemma \ref{lm: perturbation bound}. This finishes the induction and the proof.
\quad $\blacksquare$

{\noindent \bf Proof of Proposition \ref{prop: init-phase-retrieval}}. First, under the assumptions assumed in the proposition and \cite[Theorem 1]{ma2019implicit}, we have 
\begin{equation*}
	\min( \| \widetilde{\x}^{t} - \x^* \|_2, \| \widetilde{\x}^t + \x^* \|_2 ) \leq \epsilon ( 1- \eta \|\x^*\|_2^2/2 )^t \|\x^*\|_2, \quad \forall t \geq 0,
\end{equation*} holds for some $ \epsilon \in (0,1)$ with probability at least $1- c_3 n p^{-5}$. So when $T_0 \geq c_2 \log p \cdot \log (\|\x^*\|_2 p)$, we have 
\begin{equation*}
	\|\widetilde{\x}^{T_0} \widetilde{\x}^{T_0 \top} - \X^*\|_\F \overset{(a)}\leq \frac{9}{4} \|\X^*\|_\F \min( \| \widetilde{\x}^{T_0} - \x^* \|_2, \| \widetilde{\x}^{T_0} + \x^* \|_2 ) \leq c \|\X^*\|_\F/p.
\end{equation*} Here (a) is by \cite[Lemma 5.3]{tu2016low}. This finishes the proof of initialization and the rest of the proof follows from Theorem \ref{thm: quadratic convergence for phase retrieval}. \quad $\blacksquare$

 \section{Additional Proofs and Technical Lemmas}
 We collect the additional proofs and technical lemmas that support the main technical results in this section.

{\noindent \bf Proof of Equation \eqref{eq: sketching in R2RILS}.}
First, we denote 
\begin{equation*} 
	\U^{t} = [\u_1, \ldots, \u_{p_1}]^\top,
	\V^{t} = [\v_1, \ldots, \v_{p_2}]^\top,
	\M = [\m_1,\ldots,\m_{p_1}]^\top,
	\N = [\n_1,\ldots, \n_{p_2}]^\top.
\end{equation*} Then
\begin{equation}\label{eq: R2RILE lemma proof equation1}
\begin{split}
		\argmin_{\substack{\M \in \bbR^{p_1 \times r},\\ \N \in \bbR^{p_2 \times r}}} \sum_{(i,j) \in \Omega} \left\{\left(\U^t \N^\top + \M \V^{t \top} - \X\right)_{[i,j]} \right\}^2
		 =  \argmin_{\substack{\M \in \bbR^{p_1 \times r},\\ \N \in \bbR^{p_2 \times r}}} \sum_{(i,j) \in \Omega} (\u_i^\top \n_j + \m_i^\top \v_j - \X_{[i,j]} )^2.
\end{split}
\end{equation} 
If $(i,j) \in \Omega$, then the corresponding design matrix $\A^{ij}$ has $1$ at location $(i,j)$ and $0$ at the rest of locations. Then
\begin{equation*}
	\A^{ij} \V^{t} = [0,\ldots, \overbrace{\v_j}^{i^{th}},\ldots,0 ]^\top, \quad \A^{ij\top} \U^{t} = [0,\ldots, \overbrace{\u_i}^{j^{th}},\ldots,0 ]^\top.
\end{equation*} So on the sketching perspective of R2RILS, we have 
\begin{equation*}
\begin{split}
&\argmin_{\M \in \bbR^{p_1 \times r}, \N \in \bbR^{p_2 \times r}} \sum_{(i,j) \in \Omega} \left( \langle \U^{t\top} \A^{ij} , \N^\top \rangle +  \langle \M, \A^{ij} \V^t \rangle  - \X_{[i,j]} \right)^2\\
	= &\argmin_{\M \in \bbR^{p_1 \times r}, \N \in \bbR^{p_2 \times r}} \sum_{(i,j) \in \Omega} ( \langle \N, [0,\ldots, \overbrace{\u_i}^{j^{th}},\ldots,0 ]^\top \rangle + \langle \M, [0,\ldots, \overbrace{\v_j}^{i^{th}},\ldots,0 ]^\top \rangle - \X_{ij} )^2 \\
	= &\argmin_{\M \in \bbR^{p_1 \times r}, \N \in \bbR^{p_2 \times r}} \sum_{(i,j) \in \Omega} (\u_i^\top \n_j + \m_i^\top \v_j - \X_{[i,j]} )^2,
\end{split}
\end{equation*} which is exactly the same as \eqref{eq: R2RILE lemma proof equation1} and this finishes the proof.
\quad $\blacksquare$
 
 {\noindent \bf Proof of Lemma \ref{lm: Riemannian Hessian}}. 
 The proof is the same as the proof of Proposition 2.3 \cite{vandereycken2013low}, except here we need to replace the gradient in the matrix completion setting to the gradient $\cA^*(\cA(\X) - \y)$ in our setting. \quad $\blacksquare$
 
\begin{Lemma}[Projection onto the Positive Semidefinite Cone in the Nuclear norm]\label{lm: PSD projection in nuclear norm}
Given any symmetric matrix $\A \in \bbR^{p \times p}$, and denotes its eigenvalue decomposition as $\sum_{i=1}^p \lambda_i \v_i \v_i^\top$ with $\lambda_1 \geq \cdots \geq \lambda_p$. Let $\A_0 = \sum_{i=1}^p (\lambda_i \vee 0 )\v_i \v_i^\top$, then
\begin{equation*}
	\A_0 = \argmin_{\X \in \S^p_+} \|\A - \X\|_*,
\end{equation*} here $\S^p_+$ is the set of $p \times p$ positive semidefinite (PSD) matrices.
\end{Lemma}
{\noindent \bf Proof of Lemma \ref{lm: PSD projection in nuclear norm}}
Here the main property we use is the variational representation of nuclear norm. Let $m = \max \{i: \lambda_i \geq 0\}$. For any PSD matrix $\X$,
\begin{equation*}
	\begin{split}
		\|\X - \A \|_* \geq \sum_{i=1}^{p-m} \sigma_i(\X - \A)&= \sup_{\U \in \bbO_{p , (p-m)}, \V \in \bbO_{p, (p-m)}} \tr(\U^\top (\X - \A) \V)\\
		&\geq \sup_{\U \in \bbO_{p, (p-m)}} \tr(\U^\top (\X - \A)\U) \geq 0 - \inf_{\U \in \bbO_{p, (p-m)} }\tr(\U^\top \A \U) \\
		& \geq - (\sum_{i=m+1}^p \lambda_i). 
	\end{split}
\end{equation*}
On the other hand, $\|\A_0 - \A\|_* =  - (\sum_{i=m+1}^p \lambda_i)$ and this finishes the proof.
\quad $\blacksquare$

\begin{Lemma}[Bounds for spectrum of $\cL^* \cA^* \cA \cL$ in Phase Retrieval] \label{lm: nuclear norm bound of AtopA in phase retrieval}
For any given unit vector $\u\in\bbR^p$, define the linear map 	  
 \begin{equation*} 
 	\cL:  \W= \left[ \begin{array}{c c}
	w_0 \in \bbR & \w_1^\top \in \bbR^{1 \times (p-1)}\\
	\w_1 \in \bbR^{(p-1) \times 1} & \0
\end{array}
 \right] \to [\u \, \u_\perp] \left[ \begin{array}{c c}
	w_0 & \w_1^{\top}\\
	\w_1 & \0	
	\end{array}
	  \right] [\u \, \u_{\perp} ]^\top.
 \end{equation*} It is easy to compute $\cL^*(\M) = \left[ \begin{array}{c c}
	\u^{\top} \M \u & \u^{\top} \M \u_\perp \\
	\u_{\perp}^{\top} \M \u & 0
	\end{array}
	  \right]$, where $\M\in \bbR^{p\times p}$ is a symmetric matrix. Suppose $\a_i \overset{i.i.d.}\sim N(0, \I_p)$. Then $\forall \delta \in (0,1)$, $\exists C(\delta) > 0$ such that when $n \geq C(\delta) p \log p$, with probability at least $1 - c_1\exp(-c_2(\delta) p) - c_3 n^{-p}$, we have for any $\u$ and $\M \in \Range(\cL^*)$
\begin{equation} \label{ineq: spectrum lower bound of AtopA in PR}
    \|\cL^* \cA^* \cA \cL(\M)\|_\F \geq \frac{1-\delta}{2} n \|\M\|_\F.
\end{equation}
where $\cA$ is the linear map in \eqref{eq: matrix model pha} generated by $\{\a_i\}_{i=1}^n$. 
Also for any $\u$ and matrix $\M \in \Range(\cL^*)$, with the same high probability, we have
\begin{equation} \label{ineq: PR AtopA inverse norm}
	\|(\cL^* \cA^* \cA \cL)^{-1}(\M)\|_* \leq \frac{4}{(1-\delta)n} \|\M\|_*,
\end{equation} 
\end{Lemma}
{\noindent \bf Proof of Lemma \ref{lm: nuclear norm bound of AtopA in phase retrieval} } Note that \eqref{ineq: spectrum lower bound of AtopA in PR} is true when $\M$ is a zero matrix. When $\M$ is non-zero, $\cL(\M) = \u \m^\top + \m \u^\top$ for some $\m$. Then $\frac{1}{n}\|\cA\cL(\M)\|_2^2 = \frac{1}{n} \sum_{i=1}^n |\a_i^\top \u |^2 |\a_i^\top \m |^2$ and $\|\M\|^2_\F = \|\cL(\M)\|^2_\F = 2(|\m^\top \u|^2 + \|\m\|_2^2 \|\u\|_2^2)$. For any $\cL$ and $\M$, with probability  $1 - c_1\exp(-c_2(\delta) p) - c_3 n^{-p}$, we have
\begin{equation} \label{ineq: PR spec lower bound}
    \|\cL^* \cA^* \cA \cL(\M)\|_\F = \|\cA \cL(\M)\|_2^2/\|\M\|_\F \overset{(a)}\geq \frac{1-\delta}{2} n \|\cL(\M)\|^2_\F/\|\M\|_\F = \frac{1-\delta}{2} n \|\cL(\M)\|_\F,
\end{equation} where (a) is due to the \cite[Lemma 6.4]{sun2018geometric}.

Next, we prove \eqref{ineq: PR AtopA inverse norm}. First suppose $\W \in \Range(\cL^*)$ satisfies $(\cL^* \cA^* \cA \cL)(\W) = \M$, then we have \begin{equation} \label{eq: operator inverse result}
\begin{split}
	\frac{\|(\cL^* \cA^* \cA \cL)^{-1}\M\|_*}{\|\M\|_*} &=  \frac{\|\W\|_*}{\|\cL_t^* \cA^* \cA \cL_t(\W)\|_*}.
\end{split}	
\end{equation}

Hence, to prove the desired result, we only need to obtain a lower bound of $\|\cL_t^* \cA^* \cA \cL_t(\W)\|_*$:
\begin{equation} \label{ineq: AtopA PR spec bound}
	\begin{split}
		\| \cL^* \cA^* \cA \cL(\W) \|_* = \sup_{\|\Z\| \leq 1 } \langle \cA \cL(\W) , \cA \cL (\Z) \rangle &\geq  \langle \cA \cL(\W) , \cA \cL (\frac{\W}{\|\W\|}) \rangle\\
		& = \|\cA \cL(\W)\|_2^2/\|\W\|\\
		& \overset{(a)} \geq   \frac{(1-\delta)}{2} n \|\cL(\W)\|_\F^2/ \|\W\| \overset{(b)}\geq \frac{1-\delta}{4}n\|\W\|_*,  
	\end{split}
\end{equation} where (a) holds for any $\cL, \W$ with probability  $1 - c_1\exp(-c_2(\delta) p) - c_3 n^{-p}$ by the same reason as (a) in \eqref{ineq: PR spec lower bound}; (b) is true because $\W \in \Range(\cL^*)$  and $\|\cL(\W)\|_\F = \|\W\|_\F \geq  \|\W\|_*/\sqrt{2}$.
\quad $\blacksquare$

\begin{Lemma}[Upper Bound for $\|\cL^* \cA^*(\z)\|_*$ in Phase Retrieval] \label{lm: nuclear norm bound for LtA}
Consider the same linear operator $\cL$ as in Lemma \ref{lm: nuclear norm bound of AtopA in phase retrieval}. Suppose $\a_i \overset{i.i.d.}\sim N(0, \I_p)$ and $\cA$ is the linear map in \eqref{eq: matrix model pha} generated by $\{\a_i\}_{i=1}^n$. Then there exists $C, c_1, c_2 > 0$ such that when $n \geq Cp$, with probability at least $ 1- c_1\exp(-c_2 p)$, for any $\cL$ and $\z$, we have $\|\cL^* \cA^*(\z)\|_* \leq c p\|\z\|_1$ for some $c > 0$, where $\|\z\|_1$ denotes the $\ell_1$ norm of $\z$.
\end{Lemma}
{\noindent \bf Proof of Lemma \ref{lm: nuclear norm bound for LtA}.}
The proof is based on the concentration of sub-exponential random variables. First, for fixed $\cL$, we have
\begin{equation*}
	\begin{split}
		&\sup_{\|\z\|_1 \leq 1} \|\cL^* \cA^* (\z)\|_* = \sup_{\|\z\|_1 \leq 1} \sup_{\W \in \S,\|\W\|\leq 1} \langle \cL^* \cA^*(\z), \W \rangle\\
		 =& \sup_{\|\z\|_1 \leq 1} \sup_{\W \in \S,\|\W\|\leq 1} \langle \z, \cA \cL(\W) \rangle = \sup_{\W \in \S,\|\W\|\leq 1}\sup_{\|\z\|_1 \leq 1} \langle \z, \cA \cL(\W) \rangle = \sup_{\W \in \S,\|\W\|\leq 1} \|\cA \cL(\W)\|_{\infty},
	\end{split}
\end{equation*} here $\S$ is the set of symmetric matrices and $\|\cA \cL(\W)\|_{\infty}$ denotes the largest absolute value in the vector $\cA \cL(\W)$.

Notice $\cL (\W)$ is a symmetric rank-2 matrix with spectral norm bounded by $1$, without loss of generality, we can consider the bound for $\|\cA (\M)\|_{\infty}$ for fixed rank 2 matrix $\M$ with eigenvalue decomposition $\u_1 \u_1^\top - t\u_2\u_2^\top$ and $t \in [-1,1]$. In this case $\|\cA(\M)\|_\infty = \max_i | |\a_i^\top \u_1|^2 - t|\a_i^\top \u_2|^2 |$ and $ | |\a_i^\top \u_1|^2 - t|\a_i^\top \u_2|^2 |$ is a subexponential random variable. By the concentration of subexponential random variable \cite{vershynin2010introduction}, we have
\begin{equation*}
	\bbP( | |\a_i^\top \u_1|^2 - t|\a_i^\top \u_2|^2 | -\xi > x  ) \leq \exp(-c \min(x^2, x)),
\end{equation*} where $\xi = \bbE | |\a_i^\top \u_1|^2 - t|\a_i^\top \u_2|^2 |$. And a union bound yields
\begin{equation} \label{ineq: union bound after maximum}
	\bbP( \max_i \{ | |\a_i^\top \u_1|^2 - t|\a_i^\top \u_2|^2 | -\xi \} > x ) \leq n\exp(-c \min(x^2, x)).
\end{equation}

Next, we use the $\epsilon$-net argument to extend the bound to hold for any symmetric rank 2 matrices $\M$ with spectral norm bounded by $1$. Notice that by proving that, we also prove the desired inequality for any $\cL$. Let $S_\epsilon$ be an $\epsilon$-net on the unit sphere, $T_\epsilon$ be an $\epsilon$ net on $[-1,1]$ and set 
\begin{equation*}
	N_\epsilon = \{ \M= \u_1 \u_1^\top - t\u_2 \u_2^\top: (\u_1, \u_2, t) \in S_\epsilon \times S_\epsilon \times T_\epsilon  \}.
\end{equation*} Since $|S_\epsilon| \leq (3/\epsilon)^p$, we have $|N_\epsilon| \leq (3/\epsilon)^{2p+1}$. A union bound yields 
\begin{equation} \label{ineq: union bound after epsilon net}
	\bbP( \forall \M \in N_\epsilon, \max_i \{ | |\a_i^\top \u_1|^2 - t|\a_i^\top \u_2|^2 | -\xi \}> x ) \leq n\exp(-c \min(x^2, x) + (2p+1) \log (3/\epsilon)).
\end{equation}

Now suppose $(\u_1^*, \u_2^*, t^*) = \argmax_{\u_1, \u_2, t\in [-1,1]} \|\cA(\M)\|_{\infty}$ and denote $\M^* = \u_1^* \u_1^{*\top} - t^* \u_2^* \u_2^{*\top}$, $\mu = \|\cA(\M^*)\|_\infty$. Then find the approximation $\M_0 = \u_0 \u_0^\top - t_0 \v_0 \v_0^\top \in N_\epsilon$ such that $\|\u_0 - \u^*\|_2, \|\v^*-\v_0\|_2, |t^*-t_0| $ are each at most $\epsilon$. First notice
\begin{equation*}
	\begin{split}
		\|\u^* \u^{*\top} - \u_0 \u_0^\top\| &= \sup_{\|\x\|_2 = 1} | |\u_0^\top \x|^2 - |\u^{*\top} \x|^2  |\\
		&= \sup_{\|\x\|_2 = 1}  |(\u_0-\u^*)^\top \x| |(\u_0+\u^*)^\top \x|\\
		& \leq \|\u^*-\u_0\|_2 \|\u^*+\u_0\|_2 \leq 2 \|\u^* - \u_0\|_2 \leq 2\epsilon.   
	\end{split}
\end{equation*} Using the above bound, we have
\begin{equation}\label{ineq: closeness in epsilon net}
	\|\M^* - \M_0\| \leq \|\u^* \u^{*\top} - \u_0 \u_0^\top\| + |t^*-t_0| \|\v^* \v^{*\top}\| + |t_0| \|\v^* \v^{*\top} - \v_0 \v_0^\top\| \leq 5\epsilon.
\end{equation}
Now take $x \geq cp$ for some $c > 0$, then from \eqref{ineq: union bound after epsilon net}, we have the event $\{ \|\cA(\M_0)\|_\infty \leq c p \}$ happens with probability at least $1 - \exp(-C p)$. And on this event, we have
\begin{equation*}
	\mu \leq \|\cA(\M^* - \M_0)\|_\infty + \|\cA(\M_0)\|_\infty \overset{(a)}\leq 5\epsilon \cdot 2\mu +  cp \quad  \Longrightarrow \quad  \mu \leq \frac{cp}{(1-10\epsilon)},
\end{equation*} where (a) is by triangle inequality and the fact that $\M^* - \M_0$ can be decomposed into the sum of two rank $2$ symmetric matrices with spectral norm bounded by $5\epsilon$ \eqref{ineq: closeness in epsilon net}. Take $\epsilon < 1/10$, we get $\|\cA(\M^*)\|_\infty < cp$ for some $c >0$ with probability at least $1 - \exp(-C p)$. This finishes the proof.
\quad $\blacksquare$

\begin{Lemma} \cite[Lemma 3 ]{zhang2020islet}\label{lm:FGH}
	Suppose $\bF, \widehat{\bF}\in \mathbb{R}^{p_1\times r}, \G, \widehat{\G}\in \mathbb{R}^{r\times r}, \H, \widehat{\H}\in \mathbb{R}^{r\times p_2}$. If $\G$ and $\widehat{\G}$ are invertible, $\|\bF\G^{-1}\|\leq \lambda_1$, and $\|\widehat{\G}^{-1}\widehat{\H}\| \leq \lambda_2$, we have
	\begin{equation}\label{ineq:FGH1}
	\left\|\widehat{\bF}\widehat{\G}^{-1}\widehat{\H} - \bF\G^{-1}\H \right\|_\F \leq \lambda_2\|\widehat{\bF}-\bF\|_F + \lambda_1\|\widehat{\H}-\H\|_F + \lambda_1\lambda_2\|\widehat{\G}-\G\|_\F.
	\end{equation}
\end{Lemma}

\begin{Lemma}\cite[Lemma 3.3]{candes2011tight}\label{lm:retricted orthogonal property}
    Let $\Z_1, \Z_2 \in \bbR^{p_1 \times p_2}$ be two low rank matrices with $r_1 = \rank(\Z_1), r_2 = \rank(\Z_2)$. Suppose $\langle \Z_1, \Z_2 \rangle = 0$ and $r_1 + r_2 \leq \min(p_1, p_2)$. 
    Then
    \begin{equation*}
        |\langle \cA(\Z_1), \cA(\Z_2) \rangle | \leq R_{r_1 + r_2} \|\Z_1\|_{\F} \|\Z_2\|_{\F}.
    \end{equation*}
\end{Lemma}

\begin{Lemma}\label{lm: perturbation bound}
    Let $\X_1 = \U_1 \bSigma_1 \V_1^\top$ and $\X_2 = \U_2 \bSigma_2 \V_2^{\top}$ be two rank $r$ matrices with corresponding singular value decompositions. Then
    \begin{equation*}
    \left\{ 
    \begin{aligned}
    &\|\U_1 \U_1^\top - \U_2 \U_2^\top\| \leq \frac{\|\X_1 - \X_2\|}{\sigma_r(\X_1)\vee \sigma_r(\X_2)}, \\[2pt]
    & \max\{\|\sin \Theta(\U_1, \U_2)\|, \|\sin \Theta(\V_1, \V_2)\| \} \leq \frac{2\|\X_1 - \X_2\|}{\sigma_r(\X_1)\vee \sigma_r(\X_2)}.
    \end{aligned}
    \right.
    \end{equation*}
\end{Lemma}
{\noindent \bf Proof.} See Lemma 4.2 of \cite{wei2016guarantees} and Theorem 5 of \cite{luo2020schatten}. \quad $\blacksquare$







\end{document}